\documentclass[a4paper,10pt,reqno]{amsart}
\usepackage[english]{babel}
\usepackage{amssymb,amsthm,amsmath, amssymb,amscd} 
\usepackage{mathtools}
\usepackage{mathrsfs}
\usepackage{enumerate}
\usepackage{cite}
\usepackage[symbol]{footmisc}
\usepackage{color}
\usepackage{geometry}
\usepackage{multirow}
\usepackage[colorlinks=true, linkcolor=red, citecolor=blue, urlcolor=blue,hyperfootnotes=false]{hyperref}
\theoremstyle{plain} 
\newtheorem{theorem}{Theorem}[section]
\newtheorem{lemma}{Lemma}[section]

\theoremstyle{definition}
\newtheorem{definition}{Definition}[section]

\theoremstyle{remark}
\newtheorem{remark}{Remark}[section]

\newcommand{\R}{\mathbb{R}}
\newcommand{\Z}{\mathbb{Z}}

\newcommand{\N}{\mathbb{N}}
\newcommand{\normeq}[1]{{\left\vert\kern-0.25ex\left\vert\kern-0.25ex\left\vert #1 
    \right\vert\kern-0.25ex\right\vert\kern-0.25ex\right\vert}}

{\left\lbrace\begin{array}{r@{\hspace{1mm}}ll}}%
{\end{array}\right.}

\newcommand{\be}{\begin{equation}}
\newcommand{\ee}{\end{equation}}
\newcommand{\ba}{\begin{aligned}}
\newcommand{\ea}{\end{aligned}}

\title[KdV-DNN]{Error bounds for Physics Informed Neural Networks in generalized KdV equations placed on unbounded domains}

\author{Ricardo Freire}  
	\address{Departamento de Ingenier\'{\i}a Matem\'atica, Universidad de Chile.}
	\email{rfreire@dim.uchile.cl}
	\thanks{R.F. was partially funded by ANID Fondecyt FONDECYT 3230256, INRIA Lille PANDA project, and ANID 2022 Exploration 13220060.}

\author[Mu\~noz]{Claudio Mu\~noz}
\address{CNRS and Departamento de Ingenier\'{\i}a Matem\'atica and Centro
de Modelamiento Matem\'atico (UMI 2807 CNRS), Universidad de Chile, Casilla
170 Correo 3, Santiago, Chile.}
\email{cmunoz@dim.uchile.cl}
\thanks{C.\ M.'s work was partly funded by Inria Lille PANDA and Chilean research grants FONDECYT 1231250, Centro de Modelamiento Matemático (CMM) BASAL fund FB210005 for center of excellence from ANID-Chile, and Grant PID2022-137228OB-I00 funded by the Spanish Ministerio de Ciencia, Innovaci\'on y Universidades, MICIU/AEI/10.13039/501100011033.}

\author{Nicol\'as Valenzuela}
\address{Departamento de Ingenier\'{\i}a Matem\'atica, Universidad de Chile.}
	\email{nvalenzuela@dim.uchile.cl}
	\thanks{N.V. was partially funded by Chilean research grants ANID 2022 Exploration 13220060, FONDECYT 1231250, INRIA Lille PANDA project, a Latin America PhD Google Fellowship, ANID-Subdirecci\'on de Capital Humano/Doctorado Nacional/2023-21231021 and Basal CMM FB210005.}

\numberwithin{equation}{section}

\usepackage{subcaption}
\begin{document}

\begin{abstract}
In this paper we study a rigorous setting for the numerical approximation via deep neural networks of the generalized Korteweg-de Vries (gKdV) model in one dimension, for subcritical and critical nonlinearities, and assuming that the domain is the unbounded real line. The fact that the model is posed on the real line makes the problem difficult from the point of view of learning techniques, since the setting required to model gKdV is structured on intricate oscillatory estimates dating from Kato, Bourgain and Kenig, Ponce and Vega, among others. Therefore, a first task is to adapt the setting of these techniques to the deep learning setting. We shall use a battery of Kenig--Ponce--Vega suitable norms and Physics Informed Neural Networks (PINNs) to describe this approximative scheme, proving rigorous bounds on the approximation for each critical and subcritical gKdV model. We shall use these results to provide clear approximation results in the case of several gKdV nonlinear patterns such as solitons, multi-solitons, breathers, among other solutions.
\end{abstract}

\subjclass[2010]{Primary: 65K10, 65M99, Secondary: 68T07}

\keywords{KdV, Physics Informed Neural Networks, Deep Neural Networks, Unbounded Domains, Solitons}

\maketitle

\section{Introduction}

\subsection{Setting} Consider the generalized Korteweg-de Vries (g-KdV) equation placed on the real line
\begin{equation}\label{eq:KdV}
	\partial_t u + \partial_{xxx} u = \mu \partial_x (u^k),\qquad (t,x)\in \R\times\R,
\end{equation}
where $k \in \N$, $k \geq 2$ and $\mu = \pm 1$. The case $\mu=-1$ corresponds to the \emph{focusing} case, which will be the main focus of this paper, while $\mu = 1$ is the \emph{defocusing} case. We assume initial data $u_0\in H^s_x(\R)$, where $s\geq 0$ is a different value for each $k$, but always above the critical regularity setting, which depends on the value of $k$ described below in \eqref{s_k}. The cases $k=2$ and $k=3$, denoted as KdV and mKdV, are particularly important since they are obtained from the Boussinesq approximation of water waves in the shallow water and long wave regime, under further approximations of the dispersion coefficient \cite{Bous}, and from movement of curves in the plane \cite{GoPe}, respectively. Additionally, they are integrable in the sense of Liouville \cite{AC,AS,Sch}. The quartic case $k=4$ is recognized as the only subcritical integer pure power that is not integrable but has globally defined solutions. The case $k=5$ is important since it is $L^2$ critical, making it a domain where one studies blow up in finite time.  

The generalized KdV model \eqref{eq:KdV} has been extensively studied during the past decades. For a detailed account, the reader may consult the monograph by Linares and Ponce \cite{LP}. It is well-known that the initial value problem associated to \eqref{eq:KdV} is locally well-posed in suitable $H^s$ spaces, and globally well-posed if $2\leq k\leq 4$ \cite{Ka,KPV1}. In the focusing case and $k=5$, there are blow up solutions \cite{MtMe1}. The Cauchy theory has been improved to minimal requirements, see \cite{Bourgain,KPV1991.0,CKSTT,Guo,Kishimoto,KiVi19}. However, for the purposes of approximation techniques, in this work we shall only require regularity at the level of Kenig, Ponce and Vega estimates. The behavior of solutions in the case of regions which are not of solitonic type, in the nonintegrable setting, has been described in \cite{MMPP1}, see also \cite{FLMP} for the extension to the closely related Benjamin-Ono model. In the supercritical regime, see the recent work \cite{FrMu}.

The model \eqref{eq:KdV} is also well-known by the existence of suitable localized nonlinear structures such as solitons, multi-solitons, solitary waves, kinks and breathers. These structures have been studied in great details during the past 50 years, starting from the original work by Fermi, Pasta, Ulam and Tsingou \cite{FPU}, and the posterior work of Zabusky and Kruskal \cite{ZaKr}. First of all, we concentrate ourselves on solitons and solitary waves, in the case $\mu=-1$,
\[
u_{c,k}(t,x) = \left(\frac{k+1}{2} c ~\hbox{sech}^2\left(\frac {k-1}2 \sqrt c (x-ct)\right)\right)^{\frac 1{k-1}}, \quad k\in [2,\infty)\cap \mathbb N, \quad c>0,
\]
which are known from the work of Korteweg and de Vries \cite{KdV}, Boussinesq \cite{Bous} and Russell \cite{Russell} in the KdV case, but holds for any focusing gKdV model. A solution is said to be orbitally stable if any initial data starting sufficiently close to the given initial profile generates a solution that remains close, for all time, to the family of states obtained from the original solution through the natural symmetries of the equation (e.g., spatial translations), that is, small perturbations of the initial condition do not cause the evolution to drift far from the orbit of the reference solution for any later time. The stability of these waves have been studied by Benjamin \cite{Benj}, Bona, Souganidis and Strauss \cite{BSS} and Weinstein \cite{Weinstein3}, in the cases $k=1,2,3,4$. The critical case $k=5$ required more work. Martel and Merle \cite{MtMe1} proved the instability of these waves, improving previous results by Bona, Souganidis and Strauss \cite{BSS}. The blow-up was proved in \cite{MM2,MM3}. Later, Martel \cite{Martel} constructed a solution asymptotic to the sum of well-decoupled solitons for gKdV with any power between 2 and 5. The $H^1$ stability of these waves was studied by Martel, Merle and Tsai \cite{MMT}, and the asymptotic stability of these waves have been considered in \cite{PW,MM1,MMnon}. Even better, Alejo, Mu\~noz and Vega \cite{AMV} have studied the sum of $N$ solitons in the KdV case and shown their stability in $L^2$.  In the KdV $(k=2)$ and the mKdV $(k=3)$ cases, Hirota \cite{Hir71,Hir72} was able to construct exact solutions for the multiple collision of solitons.

Another special structure for the mKdV $(k=3)$ case are the breather solutions. They were originally introduced by Wadati \cite{Wad}. More precisely, for $\alpha,\beta>0$ and $x_1,x_2 \in \R$, the real breather solution of the mKdV equation \eqref{eq:KdV} is explicitly given by
\begin{equation}\label{breather}
    B_{\alpha,\beta}(t,x)
:= 2\sqrt{2}\,\partial_x
\bigg[
\arctan\!\bigg(
\frac{\beta}{\alpha}
\frac{\sin(\alpha y_1)}{\cosh(\beta y_2)}
\bigg)
\bigg],
\end{equation}
where $y_1 = x + \delta t + x_1,$  $ 
y_2 = x + \gamma t + x_2$, and $ \delta := \alpha^{2} - 3\beta^{2}$ , $\gamma := 3\alpha^{2} - \beta^{2}$. Then, Alejo and Muñoz proved in \cite{AM} that mKdV breathers are orbitally stable in their natural $H^2(\R)$-topology, see also Chen and Liu \cite{ChenLiu}. The nonexistence of breathers have been considered in many works, see e.g. \cite{MuPo1,Munoz1} and references therein. 

%
%
%

\subsection{g-KdV models and deep learning techniques} Recently, deep learning techniques and the use of artificial neural networks have been proved to be able to quickly represent and reproduce many scientific problems, including linear and semilinear parabolic equations \cite{Han,Han2,Hure}, linear time-independent elliptic equations \cite{LLP,Grohs}, and both time-dependent and time-independent non-local models \cite{PLK,Raissi0,Val22,Val23}, among others. The primary scientific basis is rigorously set on the capability of these methods to approximate continuous functions on compact sets \cite{Hornik,Dmitry}, and being able to produce a quick effective answer in short time, compared with classical methods. We can mention the original works by Raissi, Perniadakis and Karniadakis \cite{RK,RPK,Raissi}, where the PINNs method was introduced and where numerical approximation of solutions for several equations were given. More recently, this methodology has shown a great capability on the approximation of solutions with singularities \cite{JGS1}. Specifically, the authors give numerical approximations of self-similar solutions, such that there exist blow-up at finite time. These solutions were addressed over the two dimensional Boussinesq equation, the one dimensional Burgers equation, and the one dimensional Córdoba-Córdoba-Fontelos equation.

The KdV model placed on a bounded domain has been extensively studied using deep learning techniques. In the works of Raissi, Perniadakis and Karniadakis, numerical approximations of the KdV soliton solutions were given. Mishra and Molinaro \cite{MM22} introduced rigorous error estimates in several models of diffusive type in the setting of approximation of solutions via deep neural networks of PINNs type. In the case of dispersive models, the situation is more complicated. Bai, Mishra and Molinaro  \cite{BKM22} considered the case of bounded domains for several dispersive equations (including KdV models), proving the first rigorous error estimates. In \cite{MuVa24} the authors proved approximation bounds in the case of nonlinear waves in light cones of unbounded domains, and recently \cite{ACFMV24} the case of the nonlinear Schr\"odinger model on unbounded domains was also described in great detail. The gKdV case has been less treated in the context of Deep Learning techniques. Ortiz, Marín and Ortiz \cite{OMO} studied the accuracy of the standard PINN methodology on gKdV models, for the nonlinearities $k=2$ and $k=4$ without boundary conditions. They made a comparison between the PINN approximation and a more classical numerical scheme, in the case of Solitons and 2-Solitons. On the other hand, Chen, Shi, He and Fang \cite{CSHF} performed the PINN methodology for KdV and mKdV equations, as well as high-order KdV models, all of them over bounded domains.

Although it is relatively direct to simulate gKdV solitons assuming a finite interval, as the previous literature shows, these approaches are not suitable if one wants to study the long time behavior of solutions. For this reason, we have followed a new direction that is seeking to provide the best setting for a future long time dynamics description of gKdV solitons and solitary waves via DNN techniques. This approach must first take into account the deep dispersive nature of the gKdV model, coming from fundamental works by Kenig--Ponce--Vega and later deeper Fourier based developments by many researchers cited above. It is then sensible to state this point as the basis of our numerically based work. Then, there is an important forthcoming step, which will be the mixing between Cauchy-theory techniques, DNNs, and (orbital, asymptotic) stability analysis adapted to the two previous settings. For the moment, we provide a suitable description of the first steps. 

\subsection{Main results} In this work we prove quantitative estimates on the approximation of g-KdV solutions by deep learning approximations of PINNs type. This approach is based on the fact that abstract stability estimates inherent to dispersive models, which consider residual and other error estimates, are suitable for the numerical approximation of solutions in the gKdV setting. To prove this result, we need first a series of assumptions that are reasonable in the literature: a good integration scheme, a suitable PINNs method, and very importantly, a deep understanding of the local and global well-posedness theory in terms of Kenig, Ponce and Vega norms. We first consider the numerical hypotheses required for our main results.

\subsection*{Numerical hypotheses} Now we proceed to introduce the numerical approximations of the classical norms (quadratures) needed in this paper. Following \cite{ACFMV24}, we start with the following definition.

\begin{definition}[Quadratures] Let $I$ be a bounded time interval containing zero. Let $g = g(t,x)$ be a bounded continuous function. Let $N,M \geq 1$ be integers and let $(w_{\ell,j})_{\ell,j=1}^{M,N} \subseteq \R_+$ be matrix weights. For any fixed sequence of collocation point $(x_j)_{j=1}^N\subseteq \R$ and times $(t_\ell)_{\ell=1}^M\subseteq I$, consider the approximative norms 
\begin{equation}\label{Loss1}
\begin{aligned}
\mathcal J_{L_x^2,N_1}[f]:= &~{} \left(  \frac1{N_1}\sum_{j=1}^{N_1} w_{1,j} |f(x_j)|^2 \right)^{1/2},\\
 \mathcal J_{p,q,N,M}[g]:=&~{} \left( \frac1{N}\sum_{j=1}^{N} \left( \frac1{M} \sum_{\ell=1}^{M}   w_{j,\ell} |g(t_{\ell},x_{j})|^q \right)^{p/q} \right)^{1/p}, \\
\mathcal K_{q,p,M,N}[g]:=&~{} \left( \frac1{M}\sum_{\ell=1}^{M} \left( \frac1{N} \sum_{j=1}^{N}   w_{j,\ell} |g(t_{\ell},x_{j})|^q \right)^{p/q} \right)^{1/p}.
 \end{aligned}
\end{equation}
The values of $p=\infty$ or $q=\infty$ are valid by changing the summation over $(x_j)_{j=1}^N$ (resp. $(t_\ell)_{\ell=1}^M$) by a maximum over the same sequences, and by neglecting the exponent $p$ (resp. $q$). As an example, if $q = \infty$,
\[
\mathcal J_{p,\infty,N,M}[g]=~{} \left( \frac1{N}\sum_{j=1}^{N} \left( \max_{\ell=1,\ldots,M}   w_{j,\ell} |g(t_{\ell},x_{j})| \right)^{p} \right)^{1/p}.
\]
\end{definition}
The classical Riemann's sums are recovered with weights $w_{n,j}=1$ and $w_{n,\ell,j}=1$. Usually, in numerical computations the collocation points and times will be chosen uniformly separated, emulating the simplest classical Riemann's sum.

\medskip

\noindent
{\bf Hypotheses on integration schemes}. Following \cite{MM22}, our first requirement is a suitable way to approach integration.
\begin{enumerate}
\item[(H1)] Efficient integration. There exist efficient approximative rules for the computations of the $L^2_x(\mathbb R)$, $L^p_x L^q_t(I\times \mathbb R)$ and $L^q_t L^{p}_x(I\times \mathbb R)$ norms, in terms of the quantities defined in \eqref{Loss1}, in the following sense: 
\begin{itemize}
\item For any $\delta>0$, and all $f\in L^2_x(\mathbb R)$, there exist $N_1 \in\mathbb N$, $R_1>0$, points $(x_{1,j})_{j=1}^{N_1} \subseteq [-R_1,R_1]$ and weights $(w_{1,j})_{j=1}^{N_1} \subseteq [0,1]$ such that
\[
 \left| \|f\|_{L^2_x(\mathbb R)} -\mathcal J_{L^2_x,N_1}(f) \right| <\delta.
\]
\item  For any $\delta>0$, $p,q\in[1,\infty]$, and all $g\in L^p_x L^q_t(I\times \R)$, there are $N_2,M_2 \in\mathbb N$, $R_2>0$, points $(t_{2,\ell}, x_{2,j})_{\ell,j=1}^{M_2,N_2} \subseteq I\times [-R_2,R_2]$ and weights $(w_{2,\ell,j})_{\ell,j=1}^{M,N} \subseteq \R_+$ such that
\[
 \left| \|g\|_{L^q_x L^p_t(I\times \R)} -\mathcal J_{p,q,N_2,M_2}[g] \right| <\delta.
\]
\item  For any $\delta>0$, $p,q \in[1,\infty]$, and all $g\in L^{p}_t L^{q}_x(I\times \R)$, there are $N_3,M_3 \in\mathbb N$, $R_3>0$, points $(t_{3,\ell}, x_{3,j})_{j=1}^{M_3,N_3} \subseteq I\times [-R_3,R_3]$ and weights $(w_{3,\ell,j})_{\ell,j=1}^{M,N} \subseteq \R_+$ such that
\[
 \left| \|g\|_{L^p_t L^{q}_x(I\times \R)} -\mathcal K_{q,p,M_3,N_3}[g] \right| <\delta.
\]
\end{itemize}
\end{enumerate}
Notice that these hypotheses have been already considered in the literature in several cases. Mishra and Molinaro \cite{MM22} used similar hypotheses (actually more quantitative than ours) to describe rigorous approximations of semilinear heat models on bounded domains. In \cite{ACFMV24}, similar hypotheses were used to describe the approximations of nonlinear Schr\"odinger models in the case of unbounded domains, where $L^q_tL^r_x$ mixed Strichartz estimates are needed. In our situation, we shall use the norms cited above but with derivatives, assuming the classical Kenig--Ponce--Vega \cite{KPV1} theory.

\subsection*{Functional setting} More precisely, we consider for $k\in\{2,3,4,5\}$ and $s\geq 0$ the norm $\|\cdot\|_{Y_{k,s}}$ as follows:for the time-space region $A = I\times\R \subseteq \R_t\times\R_x$ or $A = \R\times I \subseteq \R_x\times\R_t$ (with $I\subseteq\R$ a time interval containing zero) depending on the order of integration, as follows:
\begin{align}
 &\begin{aligned}
    \|u\|_{Y_{2,s}(A)} &:= \|\partial_x u\|_{L_t^4L_x^{\infty}(A)} 
+ \|D_x^s\partial_x u\|_{L_x^{\infty}L_t^2(A)} \\
&~{} \quad + \|u\|_{L_x^2L_t^\infty(A)} 
+  \|u\|_{L_t^\infty H^s_x(A)};
\end{aligned} \label{eq:Y2s_norm}\\
& \begin{aligned}
\|u\|_{Y_{3,s}(A)} &:= \|\partial_x u\|_{L_x^\infty L_t^2(A)} 
+ \|u\|_{L_x^4L_t^\infty(A)} 
+ \|D^s_x \partial_x u\|_{L_x^\infty L_t^2(A)} \\
&~{} \quad + \|D_x^s u\|_{L_x^5L_t^{10}(A)}
+ \|\partial_x u\|_{L_x^{20}L_t^{5/2}(A)}; 
\end{aligned} \label{eq:Y3s_norm}\\
  &\begin{aligned}
 \|u\|_{Y_{4,s}(A)} &:= \max_{t\in I} \|D^s_x u(t)\|_{L_x^2(A)}+ \|u\|_{L_t^\infty L^2_x(A)}
+ \|u\|_{L_x^{42/13}L_t^{21/4}(A)}
+ \|u\|_{L_x^{60/13}L_t^{15}(A)} \\
&~{} \quad + \|u\|_{W_x^{s,10/3}L_t^{30/7}(A)}
+ \|\partial_x u\|_{L_x^{\infty} L_t^2(A)}+ \|D^s_x\partial_x u\|_{L_x^{\infty} L_t^2(A)}; 
\end{aligned} \label{eq:Y4s_norm}\\
 &\|u\|_{Y_{5,s}(A)} := \|u\|_{L_x^5L_t^{10}(A)}. \label{eq:Y5s_norm}
\end{align}
are the norms that we will work in the paper. These norms \eqref{eq:Y2s_norm}-\eqref{eq:Y5s_norm} arise naturally from the Kenig--Ponce--Vega theory. In addition we will say that $g \in Y_{k,s}(A)$ if $\|g\|_{Y_{k,s}(A)}<\infty$. Notice that under hypothesis (H1), the norms $\|\cdot\|_{Y_{k,s}}$ can be suitable approximated by a functional $\mathcal Y_{k,s,N,M}$ (cf. Definitions \eqref{eq:Y2s_appr}-\eqref{eq:Y5s_appr}).

\noindent\emph{Notational convention.} Throughout the paper, unless explicitly stated otherwise, all spatial norms are taken over $\R$ and all space--time norms over $I\times\R$, where $I=[-T,T]$ is a fixed time interval containing zero. Accordingly, we write $L^2_x:=L^2_x(\R)$, $L^q_TL^p_x:=L^q_t(I;L^p_x(\R))$, $L^p_xL^q_t:=L^p_x(\R;L^q_t(I))$, and $Y_{k,s}:=Y_{k,s}(I\times\R)$, and we omit the underlying domain whenever it is clear from the context. Mixed space--time norms are always written with the spatial variable listed first, i.e. $L^p_xL^q_t$ denotes $\big(\int_\R(\int_I|f|^q\,dt)^{p/q}dx\big)^{1/p}$, in agreement with the Kenig--Ponce--Vega local-smoothing and maximal-function norms used below.

\medskip

\noindent
{\bf Hypotheses on PINNs}.
 Let $u_{\text{DNN},\#}= u_{\text{DNN},\#} (t,x)$ be a smooth bounded function constructed by means of an algorithmic procedure (either SGD or any other ML optimization procedure) and realization of a suitable PINNs, in the following sense: under the framework \eqref{Loss1}, and with the hypothesis (H1) satisfied on a given quadrature, one has:
\begin{enumerate}
\item[(H2)] \emph{Uniformly bounded large time $L^\infty_tH^s_x$, initial time $H^s_x$ and large time $Y_{k,s}$ control for $u_0$ and $u_{\text{DNN},\#}$}.  There are $A,\widetilde A,B>0$ such that the following holds. For any $N_{1,0},N_{2,0},N_{3,0}\geq 1$ and $M_{2,0},M_{3,0}\geq 1$, there are $N_j\geq N_{j,0}$ and $M_j\geq M_{j,0}$,  such that each approximate norm satisfies
\[
\begin{aligned}
 \mathcal J_{L^2_x,N_1} [(u_0-u_{\text{DNN},\#}(0))]+\mathcal J_{L^2_x,N_1} [D_x^s(u_0-u_{\text{DNN},\#}(0))] &\leq \widetilde{A}, \\
 \mathcal K_{2,\infty,N_2,M_2} [u_{\text{DNN},\#}] + \mathcal K_{2,\infty,N_2,M_2} [D_x^su_{\text{DNN},\#}] &\leq A, \\
 \mathcal Y_{k,s,N_3,M_3} [u_{\text{DNN},\#}]  &\leq L.
\end{aligned}
\]
\item[(H3)] Small linear $Y_{k,s}$ and nonlinear $L_{x}^pL_t^2$ control. Given $\varepsilon>0$, and given  $N_{4,0},N_{5,0}\geq 1$ and $M_{4,0},M_{5,0}\geq 1$, there are $N_j\geq N_{j,0}$ and $M_j\geq M_{j,0}$ and corresponding approximative norms such that, if $k \leq 4$,
\[
\begin{aligned}
  &\mathcal J_{2,2,N_4,M_4} [\mathcal E_k[u_{\text{DNN},\#}]] + \mathcal J_{2,2,N_4,M_4} [D_x^s\mathcal E_k[u_{\text{DNN},\#}]] \\
  &\hspace{4cm} + \mathcal Y_{k,s,N_5,M_5} [ e^{-t\partial_x^3} (u_0-u_{\text{DNN},\#}(0))] < \varepsilon,
\end{aligned}
\]
and if $k=5$,
\[
  \mathcal J_{1,2,N_4,M_4} [D_x^{-1}\mathcal E_k[u_{\text{DNN},\#}]] + \mathcal Y_{k,s,N_5,M_5} [ e^{-t\partial_x^3} (u_0-u_{\text{DNN},\#}(0))] < \varepsilon,
\]
with 
\begin{equation}\label{def_E}
\mathcal E_k[u_\#]:= \partial_t u_\# + \partial_{xxx} u_\# +\partial_x(u_\#^k).
\end{equation}
Here $e^{-t \partial_x^3}f(x)$ represents the standard linear KdV solution issued of initial data $f$ at time $t$ and position $x$. 
\end{enumerate}
The LWP theory in \cite{KPV1} is addressed for $H^s(\R)$ spaces, for exponents $s>0$ satisfying $s\geq s_k$, $k=3,4,5$ or $s>s_k$, $k=2$, where
\begin{equation}\label{s_k}
s_k := \begin{cases}
   \frac34, &~{} k=2; \\
    \frac14,&~{} k = 3;\\
    \frac1{12},&~{} k = 4; \\
     0, &~{} k = 5,
    \end{cases}
\end{equation}
The choice of exponents is directly related to the technique that is used in this work, and we believe that they can be taken even below by using highly Fourier-analysis motivated norms. In order to provide the first statement of these DNN approximations in the gKdV setting, close to the classical regularity that DNNs have, we have chosen the KPV method as the most adapted to the use of DNNs. Under this regime, our main result reads as follows:

\begin{theorem}\label{MT}
Let $k \in \{2,3,4,5\}$ and let $s \geq s_k$ for $k \in\{3,4,5\}$ and $s>s_k$ for $k=2$. Assume now the following 
\begin{enumerate}
\item Let $u_0 \in H^s(\R)$, and assume \emph{(H1)}. 
\item Let $A,\widetilde A, L >0$ be fixed numbers, and let $0<\varepsilon<\varepsilon_1$ sufficiently small. 
\item Finally, let $u_{\emph{DNN},\#}$ be a DNN satisfying \emph{(H2)-(H3)}. 
\end{enumerate}
Then there exists a solution $u \in C(I;H^s(\R))\cap Y_{k,s}(I\times\R)$ to \eqref{eq:KdV} on $I\times\R$ with initial datum $u_0 \in H^s(\R)$ such that for all $R>0$ sufficiently large and constants $C$, there is smallness in the $Y_{k,s}$ norm:
\begin{align}
    \|u-u_{\emph{DNN},\#}\|_{Y_{k,s}(I\times[-R,R])}& \leq C(A, \widetilde{A},L)\varepsilon. \label{eq:MT2}
\end{align}
Additionally, one has bounded control of the energy norm:
\begin{align}
\|\partial_x(u^k-u^k_{\emph{DNN},\#})\|_{L_x^2L_t^2(I\times[-R,R])} & \leq C(A, \widetilde{A},L)\varepsilon, \label{eq:MT1}\\
 \|u-u_{\emph{DNN},\#}\|_{L_t^{\infty}H_x^s(I\times[-R,R])} & \leq C(A, \widetilde{A},L). \label{eq:MT3}
\end{align}

\end{theorem}

\begin{figure}[htbp]
   \centering
 \includegraphics[width=0.7\textwidth]{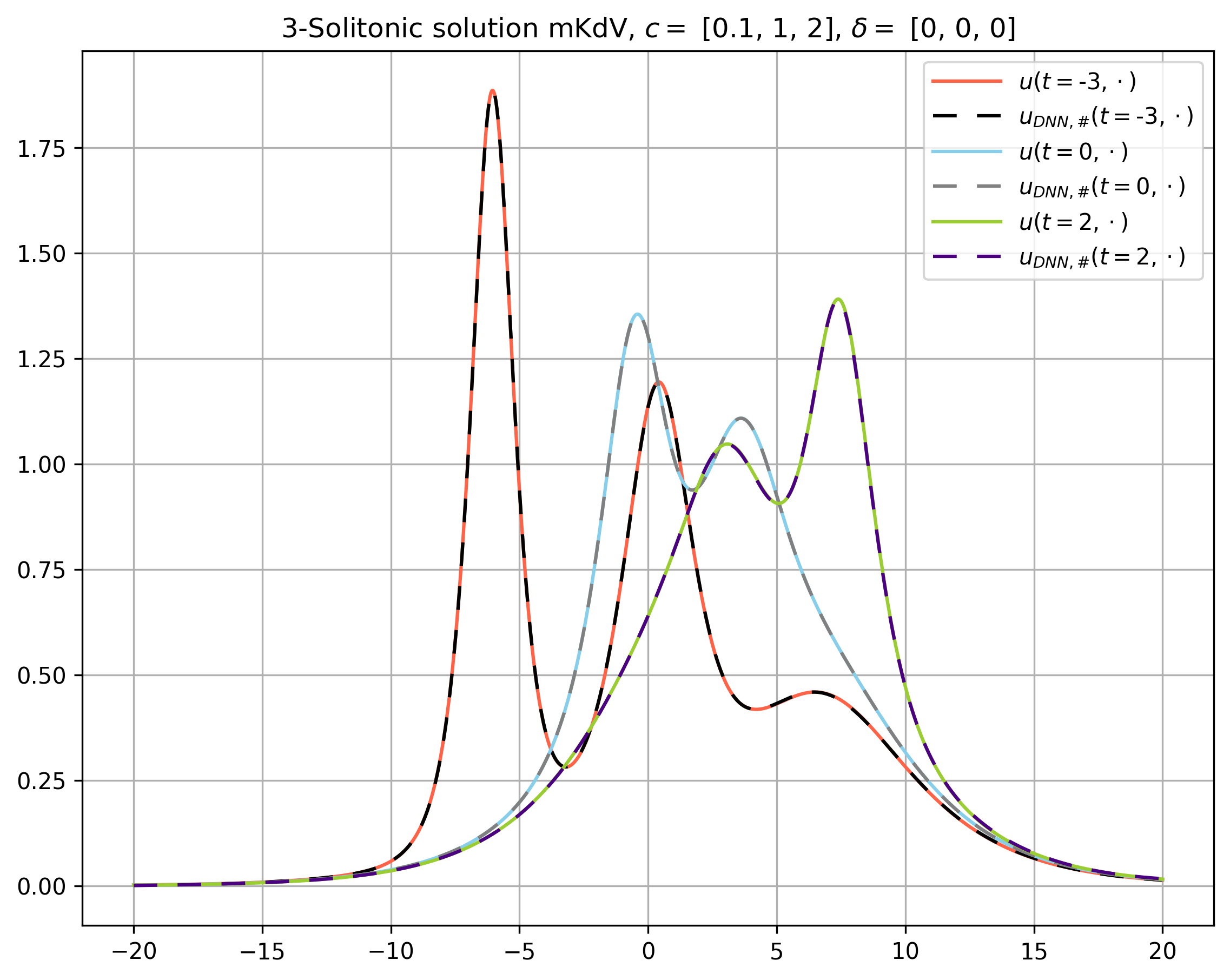}
\caption{Exact $u$ solution and predicted $u_{\text{DNN},\#}$ approximate solution in the 3-solitonic mKdV case at times $t=-3,0$ and $2$.}
\label{pastel}
\end{figure}

Theorem \ref{MT}, and in particular \eqref{eq:MT2}, provides an accurate way to find DNNs in the gKdV setting, by using the norms associated to the long time behavior theory. In this sense, it is of direct use in many applications, specially when one has a large domain, and one does not use the artificial periodic, zero or near zero boundary conditions. In computational work, we only use Strichartz based norms, and no boundary assumptions are needed. This is reflected in the case where one has several soliton solutions, as it is reflected in Fig. \ref{pastel}, that reflects how the proposed algorithm, strictly based on Theorem \ref{MT}, works, in the particular case of the complicated 3-soliton mKdV structure, and many other examples worked in more detail in Section \ref{sec:numerical}.  

\begin{figure}[ht]
\centering
\includegraphics[width=0.7\textwidth]{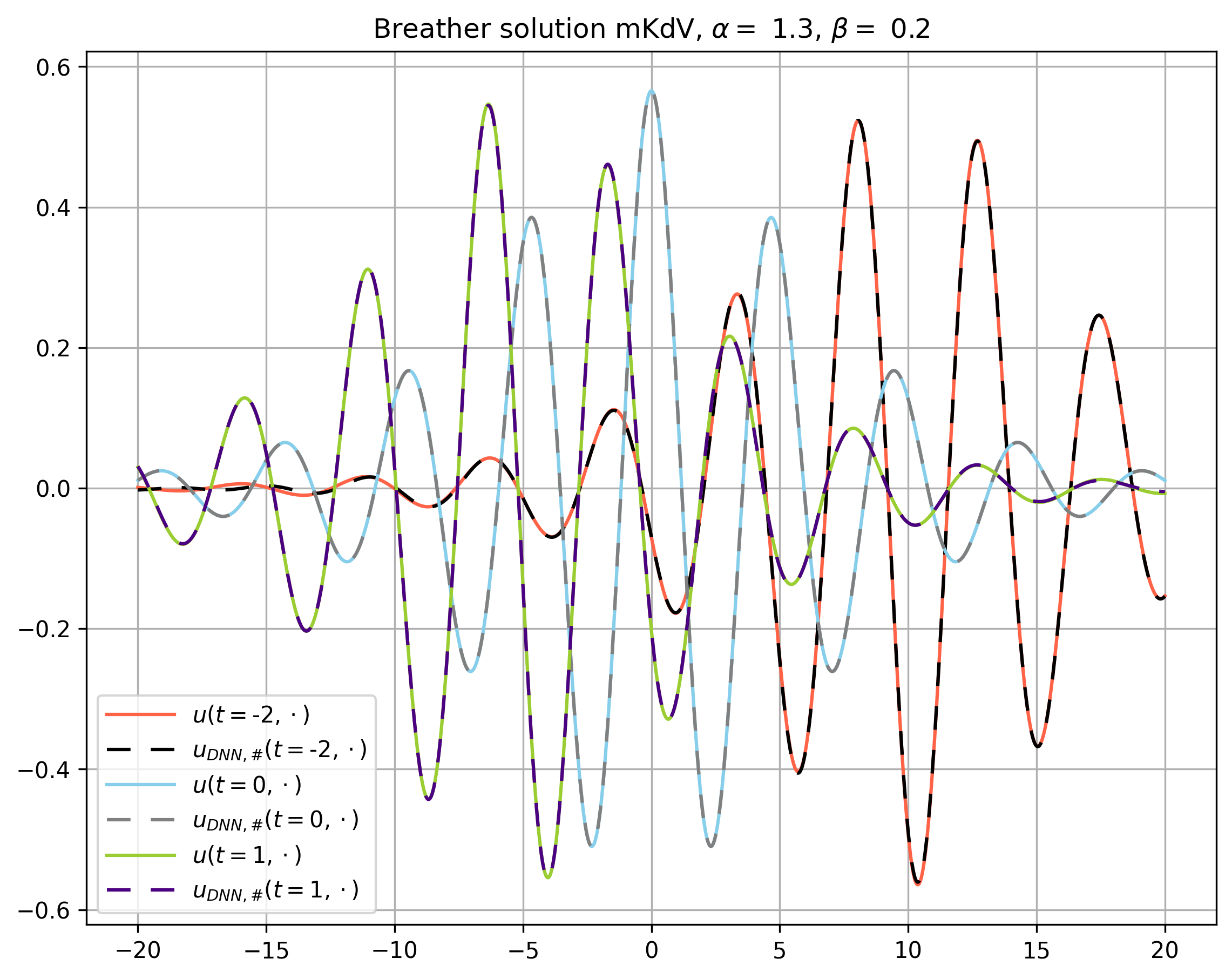}
\caption{Numerical simulation of breathers. Exact solution $u$ and predicted $u_{\text{DNN},\#}$ in the mKdV breather case with parameters $\alpha=1.3, \beta = 0.2$, and at times $t=-2, 0$ and $1$.}\label{pastel2}
\end{figure}

A key element to decide the accuracy of Theorem \ref{MT} when translated to numerical terms is to provide a good, accurate description of the highly oscillatory breather solution, in the mKdV case. This is a case that is not trivial since DNNs usually find hard to approximate high oscillations. In the case proposed in this paper, and treated in full detail in Section \ref{sec:numerical}, we provide a more than satisfactory numerical description of the breather solution, as shown in Fig. \ref{pastel2}.

\begin{remark}
Notice that the linear evolution is quickly solved using current numerical methods and the fact that $e^{-t\partial_x^3} f$ is given by
\begin{equation}\label{eq:FT-evo}
e^{-t\partial_x^3} f=\mathcal F_{\xi \to x}^{-1} \left( e^{8\pi^3it\xi^3} \mathcal F_{x\to \xi}(f)(\xi)\right),
\end{equation}
where $\mathcal F$ denotes the standard (unitary) Fourier transform with the normalization $\mathcal F_{x\to\xi}(f)(\xi)=\int_\R e^{-2\pi i x\xi}f(x)\,dx$ and $\mathcal F^{-1}_{\xi\to x}(g)(x)=\int_\R e^{2\pi i x\xi}g(\xi)\,d\xi$. In addition, for any $s \in \R$, $D_x^{s}f$ is given by
\begin{equation}\label{eq:FT-der}
D^s_xf = \mathcal F_{\xi \to x}^{-1} \left( |\xi|^s \mathcal F_{x\to \xi}(f)(\xi)\right),
\end{equation}
Previous definitions implies that for any function $f$ and scalar $s$, both $e^{-t\partial_x^3} f$ and $D^s_x f$ can be easily estimated and computed by using Fast Fourier Transform (FFT).
\end{remark} 

\begin{remark}
Theorem \ref{MT} provides the first setting for future work of the long time behavior of solitary waves, usually reflected with the notions of orbital and asymptotic stability, in the case of gKdV models. The idea is to improve and codify a combination of Theorem \ref{MT} analysis mixed with a deep use of stability techniques to provide a better description of the long time behavior, which in Theorem \ref{MT} is still influenced by the size of the time interval. This is a research direction under work that will be public in forthcoming months.
\end{remark}

\subsection{Method of proof} We prove Theorem \ref{MT} using a mixture between Kenig--Ponce--Vega stability type estimates, and the numerical properties ensured by the hypotheses (H1)-(H2)-(H3).  More precisely, the proof follows the method introduced in \cite{ACFMV24}, related to the nonlinear Schr\"odinger equation in 1D. Although the spirit is the same, we have had several complications when dealing with the gKdV models. The first and most important one is to establish the correct functional setting to prove stability estimates. This problem is not obvious since the local well-posedness theory is not fully adapted to the problem that we are considering here. The structure of the stability theorem then relies on combining the nonlinear stability framework of Kenig, Ponce, and Vega (KPV) with a rigorous analysis of the residual minimization principle underlying physics-informed neural networks (PINNs), a task that is performed both numerically and rigorously. The starting point is a given PINNs type artificial neural network approximation $u_{\mathrm{DNN},\#}$, minimizing the residual functional
\[
\mathcal E_k[u_{\mathrm{DNN},\#}] = \partial_t u_{\mathrm{DNN},\#} + \partial_x^3 u_{\mathrm{DNN},\#} + \partial_x(u_{\mathrm{DNN},\#}^k),
\]
and other related functionals. Among the other functionals, the one minimizing the linear evolution is of key interest. Notice that $\mathcal E_k[u_{\mathrm{DNN},\#}]$ can be regarded as an \emph{approximate solution} of the generalized KdV equation, but also as a nonlinear operator in suitable functional spaces, and therefore it will induce an associated functional setting for the corresponding solution. Using the Duhamel formulation and dispersive linear estimates, the difference $w = u - u_{\mathrm{DNN},\#}$ between the exact and approximate solutions satisfies an integral inequality controlled by the norms introduced in KPV. Under assumptions of small initial discrepancy and small residual (H2)--(H3), one obtains a Gronwall-type estimate leading to explicit and quantitative error bounds, roughly described as
\[
 \|w\|_{Y_{k,s}}  \lesssim C(A,\tilde A,L)\,\varepsilon.
\]
Here, $Y_{s,k}$ is the functional setting adapted to each gKdV model, obtained from \eqref{eq:Y2s_norm}-\eqref{eq:Y5s_norm}, depending on a regularity $s>s_k$ (see \eqref{s_k}) and a power nonlinearity $k$. This functional setting later will reduce to \eqref{eq:MT1}-\eqref{eq:MT3} after a suitable localization. Hence, the neural approximation inherits the same stability properties as the exact solution, with constants that can be computed explicitly, given the hypotheses. This argument provides, as far as we understand, the first rigorous quantitative justification of PINNs-based approximations for dispersive equations such as gKdV on unbounded domains.

\subsection*{Organization of this work} This work is organized as follows. Section \ref{sec:numerical} contains the set of numerical simulations that validate the main result of this paper.  In Section \ref{sec:2} we describe the basic ingredients needed for the rest of this paper. Section \ref{sec:3} is devoted to the statement of the main tool needed for the proof of Theorem \ref{MT}, which are the nonlinear long time stability estimates. Section \ref{sec:4} is completely devoted to the proof the previous statements announced in Section \ref{sec:3}. In Section \ref{sec:PMR} we prove Theorem \ref{MT}. Finally, Section \ref{sec:discussion} is devoted to discuss the current findings and the relationship with previous results in the literature, and to propose new directions for the future.

\subsection*{Acknowledgments} R. F. would like to thank prof. Felipe Poblete and UACh members, where part of this work was done. Their hospitality and support are deeply acknowledged. N. V. thanks BCAM members for their hospitality and support during research visits in 2023 and 2025. C. M. acknowledges prof. Juan Soler at Granada and BIRS Banff Canada, where part of this work was done. Part of this work was done while the authors were visiting INRIA Lille during June 2025, as part of the PANDA research project. The support and hospitality of INRIA members is greatly acknowledged.


\section{Numerical testing of Theorem \ref{MT}}\label{sec:numerical}

\subsection{Preliminaries}
This section is devoted to the numerical simulations for the solutions of gKdV equations using PINNs to validate the main Theorems, specially focusing on well known KdV solitons and breathers. All the simulations were carried out using Python on a 64-bit MacBook Pro M2 (2022) with 8GB of RAM.

Let $u_\theta := u_\theta(t,x)$ be a deep neural network with parameters $\theta$. The objective is to fit the parameters of $u_\theta$ in order to have an approximate solution satisfying the hypotheses of theorems presented in previous sections.

\subsubsection{Loss function} For each power of the nonlinearity we will use a different loss function, to be consistent with the proposed theoretical framework. For fixed number of collocation points $M_{\text{evol}}, N_{\text{evol}}, M_{\text{PDE}}, N_{\text{PDE}} \in \N$, the loss functions used are the following
\[
\mathcal L_{\text{PDE},k,s}(\theta) = \begin{cases}
    \mathcal J_{2,2,N_{\text{PDE}},M_{\text{PDE}}}\left[\mathcal E_k[u_\theta]\right] + \mathcal J_{2,2,N_{\text{PDE}},M_{\text{PDE}}}\left[D^s_x\mathcal E_k[u_\theta]\right] & \text{ if } ~ k < 5, \\
    \mathcal J_{1,2,N_{\text{PDE}},M_{\text{PDE}}}\left[D^{-1}_x\mathcal E_k[u_\theta]\right] & \text{ if } ~ k=5.
\end{cases}
\]
and
\[
\mathcal L_{\text{evol},k,s}(\theta) = \mathcal Y_{k,s,N_{\text{evol}},M_{\text{evol}}}[e^{-t\partial_x^3}(u_0-u_\theta(0))].
\]
Finally, define
\begin{equation}\label{loss}
\mathcal L_{k,s}(\theta) := \gamma_1 \mathcal L_{\text{evol},k,s}(\theta) + \gamma_2 \mathcal L_{\text{PDE},k,s}(\theta),
\end{equation}
where $\gamma_1,\gamma_2$ are properly chosen. In particular, unless we say the opposite, our algorithm will consider $\gamma_1$ as the inverse of the number of terms in $\mathcal Y_{k,s,N_1,M_1}[e^{-t\partial_x^3}(u_0-u_\theta(0))]$ and $\gamma_2=1$. Furthermore, the next simulations will be performed with $s=s_k$ for $k \geq 3$, and $s=s_k+10^{-6}$ for $k=2$, that is to say, using the loss functions $\mathcal L_{k,s_k}(\theta)$ and $\mathcal L_{2,s_k+10^{-6}}(\theta)$. 

\subsubsection{Activation function}
The choice of the activation function plays a fundamental role in the deep learning approximations. In the literature, many activation functions were considered. As for example, the usual functions are given by the hyperbolic tangent and the sigmoid, while there are several non-conventional activation functions, such as the sine \cite{SIREN}, or several wavelet-like \cite{WIRE,Wavelet} activation functions. Following \cite{WIRE}, and given the internal oscillating behavior of the solutions of gKdV and their exponentially decay, we will choose as activation function
\[
\sigma(x)=\sin(w_0x+b_0) e^{-(s_0x)^2},
\]
which corresponds to the imaginary part of a Gabor complex wavelet. To simplify, we will refer to $\sigma$ as the Wavelet activation function. The triple $(w_0,b_0,s_0)$ will be a trainable parameter in our model (in other words, it will be added to the parameters $\theta$ of the neural network), with initialization $(w_0,s_0) = \left(1,\frac{1}{\sqrt{2}}\right)$ and $b_0 \sim \mathcal U\left(-\frac\pi2,\frac\pi2\right)$. Additionally, each hidden layer will consider a different Wavelet activation function. This implies that we will have as many triples $(w_0,b_0,s_0)$ as hidden layers in the neural network.

\subsection{Methodology}

Recall that the operators $e^{-t\partial_x^3}$ and $D_x^s$ can be efficiently approximated over a grid by using fast Fourier transform (FFT) and their respective definitions in \eqref{eq:FT-evo} and \eqref{eq:FT-der}. Therefore the computation of $\mathcal L_{k,s}(\theta)$ can be performed if $(x_{1,j})_{j=1}^{N_{\text{evol}}}$ and $(x_{2,j})_{j=1}^{N_{\text{PDE}}}$ are uniform space grids.

\medskip

Let $M_{\text{evol}},M_{\text{PDE}},N_{\text{evol}},N_{\text{PDE}} \in \N$ and $R,T>0$ to be chosen in each example.

\begin{enumerate}
    \item Choose $(t_{1,\ell})_{\ell=1}^{M_{\text{evol}}}$, $(t_{2,\ell})_{\ell=1}^{M_{\text{PDE}}} \subseteq [-T,T]$ and $(x_{1,j})_{j=1}^{N_{\text{evol}}}$, $(x_{2,j})_{j=1}^{N_{\text{PDE}}} \subseteq [-R,R]$ be uniform time space grids, respectively, with $w_{1,\ell,j}=w_{2,\ell,j}=1$ for all $j$ and $\ell$.
    \item Initialize a DNN $\phi_{\theta}$ with $H$ hidden layers, $n$ neurons per layer and randomly chosen parameters $\theta$.
    \item Compute $\mathcal L_{k,s}(\theta)$ as in \eqref{loss}.
    \item Use an optimization algorithm to minimize $\mathcal L_{k,s}(\theta)$.
    \item Let
    \[
    \theta^* = \underset{\theta}{\arg\min} ~\mathcal L_{k,s}(\theta),
    \]
    be the optimal parameters obtained from minimization of previous step. $u_{\theta}:=\phi_{\theta^*}$ will be an approximation of $u$ solution of \eqref{eq:KdV}.
\end{enumerate}
The PINNs algorithm will be carried out with the L-BFGS optimizer and $n_{iter}=3000$ number of maximum iterations. In order to evaluate the method's approximation accuracy, we will compare the exact solution with the solution given by the PINN algorithm over an uniform grid $(t_\ell,x_j)_{\ell,j=1}^{M_{\text{test}},N_{\text{test}}}\subseteq[-T,T]\times[-R,R]$ of size $M_{\text{test}} \times N_{\text{test}}$. The weights $w_{\ell,j}$ will be chosen in such a way in order to recover the classical Riemann sums. There will be three different metrics for the errors, summarized as follows
\begin{enumerate}
    \item The $Y_{k,s}$ error, corresponding to an approximation of $\|u-u_\theta\|_{Y_{k,s}(I\times\R)}$ by using the operators $\mathcal Y_{k,s,N_{\text{test}},M_{\text{test}}}[u-u_\theta]$ defined in \eqref{eq:Y2s_appr}-\eqref{eq:Y5s_appr}. This will be denoted as ${\bf error}_Y$.
    \[
{\bf error}_Y := \mathcal Y_{k,s,N_{\text{test}},M_{\text{test}}}[u-u_\theta].
\]
    \item The $L_t^{\infty}H_x^s$ error, given by
    \[
    {\bf error}_{L_t^{\infty}H_x^s} := \mathcal K_{\infty,2,N_{\text{test}},M_{\text{test}}}[u-u_\theta] + \mathcal K_{\infty,2,N_{\text{test}},M_{\text{test}}}[D_x^s(u-u_\theta)],
    \]
    \item The $L_{t,x}^2$ relative error, given by
    \[
    {\bf error}_{rel} := \frac{\mathcal J_{2,2,N_{\text{test}},M_{\text{test}}}[u-u_\theta]}{\mathcal J_{2,2,N_{\text{test}},M_{\text{test}}}[u]}.
    \]
\end{enumerate}
Notice that although ${\bf error}_{rel}$ does not appear in the main Theorems, it is a common error to compare numerical approximations. In addition, some approximative norms presented in the hypothesis of Theorem \ref{MT} will be also estimated. In particular, we will give approximations to $\|u_\theta\|_{L_t^\infty H_x^s}$, $\|u_0-u_\theta(0)\|_{H_x^s}$ and $\|u_\theta\|_{Y_{k,s}}$, which yield to estimations for the constants $A$, $\widetilde A$ and $L$, respectively. All the comparisons will be performed with $N_{\text{test}}=M_{\text{test}}=300$.

\subsection{Results} Now we proceed to run the previous algorithm in the case of several test solutions.  In particular, we will focus on soliton/solitary wave solutions for the subcritical and critical cases of gKdV, the $N$-Soliton solutions for the KdV and mKdV cases, with $N=2,3$, and finally the breather solution for mKdV.

\subsubsection{The soliton case}
It is well known that for each $k \in \Z$, $k > 1$, the gKdV model \eqref{eq:KdV} has the solitary wave solution
\[
u_{c,k}(t,x) = \left(\frac{k+1}{2} c ~\hbox{sech}^2\left(\frac {k-1}2 \sqrt c (x-ct)\right)\right)^{\frac 1{k-1}},
\]
where $c>0$ is the propagation speed (see, e.g. \cite{LP}). Notice that this solution is $H^1$ stable if $k=2,3,4$ \cite{BSS}, but it is unstable if $k=5$ \cite{MtMe1}.

\medskip
 For each power of nonlinearity $k \in \{2,3,4,5\}$ we will consider the solitonic solution at velocities $c \in \{1,3\}$ over the time-space region $[-3,3]\times[-20,20]$. Notice that each setting will have a different trained neural network. In all the settings we consider $H=2$ hidden layers and $n=20$ neurons per hidden layer. In addition we choose $M_{\text{evol}} = M_{\text{PDE}} = 32$ in almost all settings except when $k=5$ and $c=3$, where we increase $M_{\text{evol}}$ to $M_{\text{evol}}= 128$, while $N_{\text{evol}},N_{\text{PDE}}$ will be increased when the power $k$ is large. For a full detail of the values of $N_{\text{evol}}$ and $N_{\text{PDE}}$, see Table \ref{tab:points-soliton}.

\begin{table}[ht]
\begin{tabular}{c|cccccccc}
\hline
$(k,c)$ & $(2,1)$ & $(2,3)$ & $(3,1)$ & $(3,3)$ & $(4,1)$ & $(4,3)$ & $(5,1)$ & $(5,3)$ \\ \hline
$N_{\text{evol}}$    & 64    & 64    & 64    & 64    & 128   & 128   & 128   & 128   \\
$N_{\text{PDE}}$    & 64    & 64    & 128   & 128   & 128   & 128   & 128   & 256   \\ \hline
\end{tabular}
\caption{Number of points $N_{\text{evol}}$ and $N_{\text{PDE}}$ used in each training, in the solitonic case. There is a natural growth in these variables as the power of the nonlinearity $k$ grows, and the speed $c$ increases, to better reflect instabilities close to the $L^2$ critical case, and strong oscillations in the high speed case.}
\label{tab:points-soliton}
\end{table}

\begin{figure}[ht]
\centering
\begin{subfigure}[t]{0.45\textwidth}
\centering
\includegraphics[width=\textwidth]{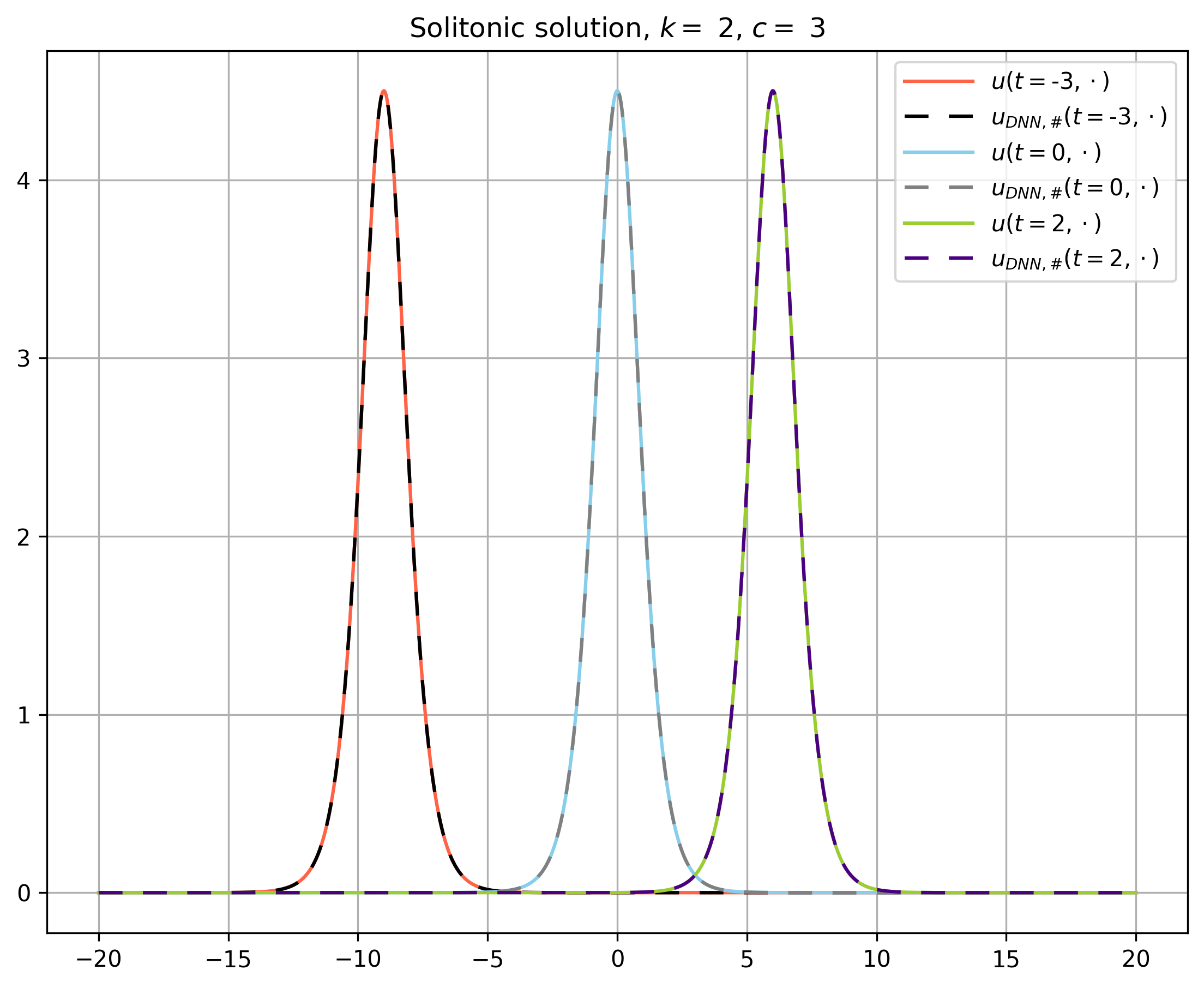}
\caption{$k=2$.}
\label{fig:Soliton-k2c3}
\end{subfigure}
\hspace{.2cm}
\begin{subfigure}[t]{0.45\textwidth}
\centering
\includegraphics[width=\textwidth]{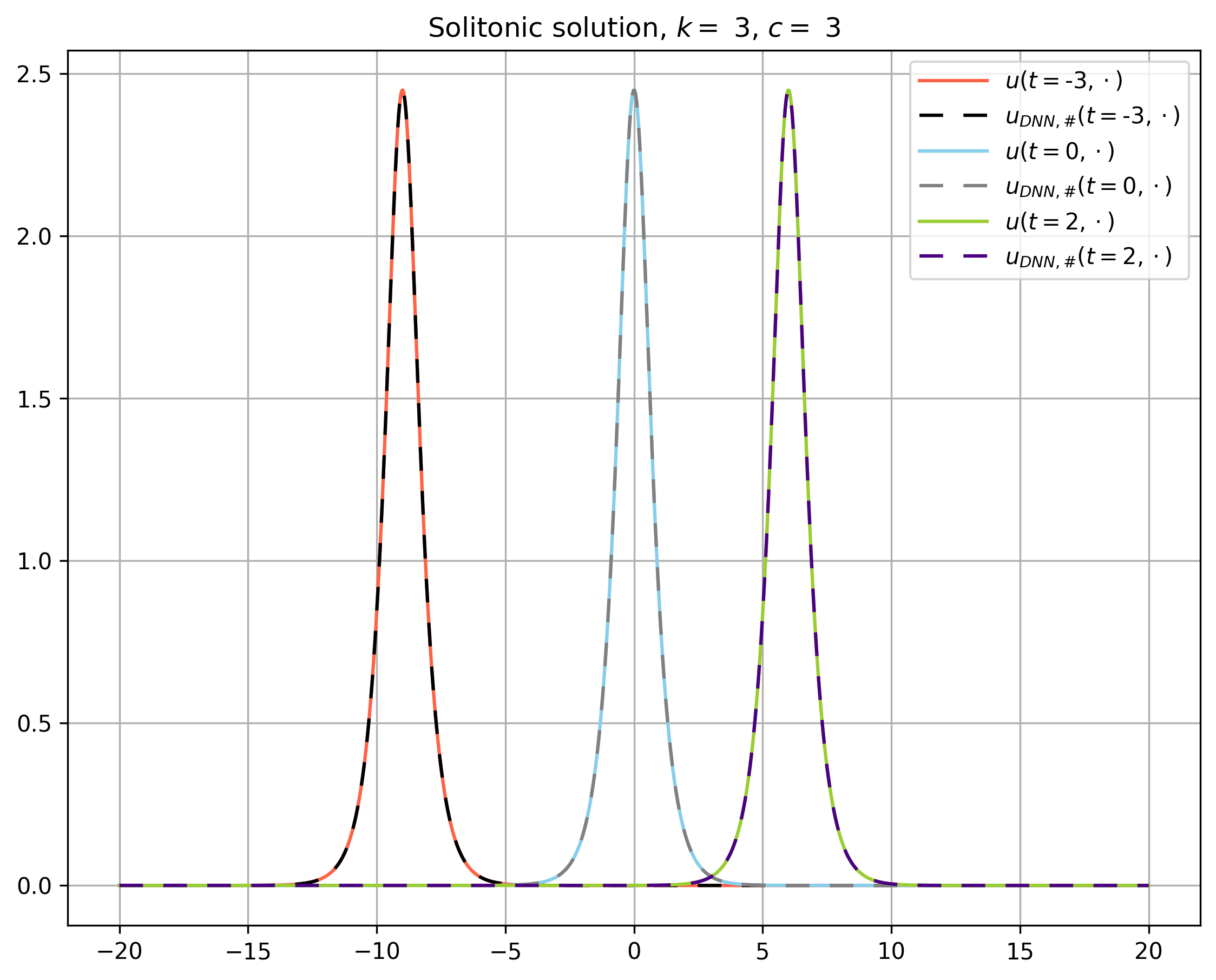}
\caption{$k=3$.}
\label{fig:Soliton-k3c3}
\end{subfigure}
\begin{subfigure}[t]{0.45\textwidth}
\centering
\includegraphics[width=\textwidth]{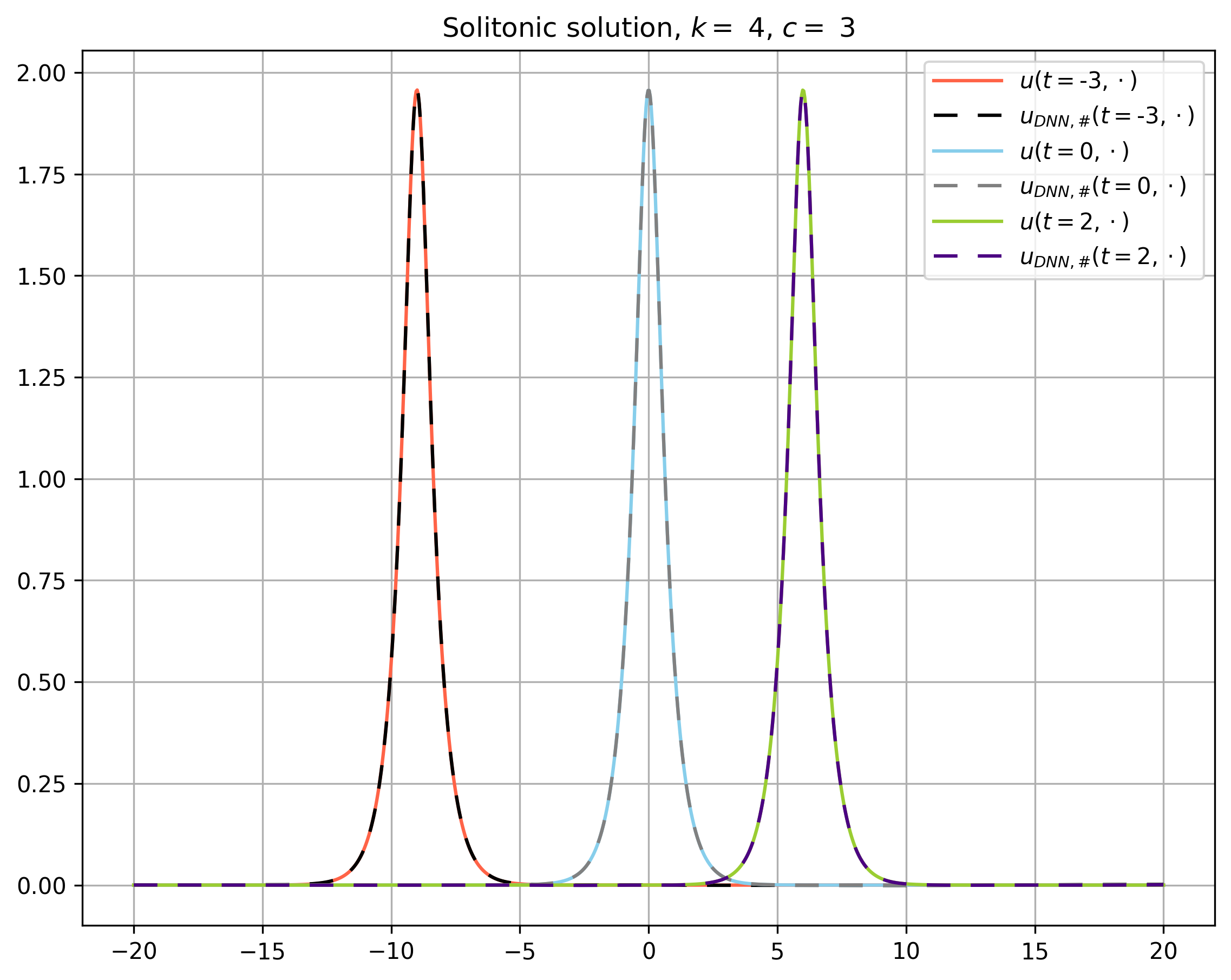}
\caption{$k=4$.}
\label{fig:Soliton-k4c3}
\end{subfigure}
\hspace{.2cm}
\begin{subfigure}[t]{0.45\textwidth}
\centering
\includegraphics[width=\textwidth]{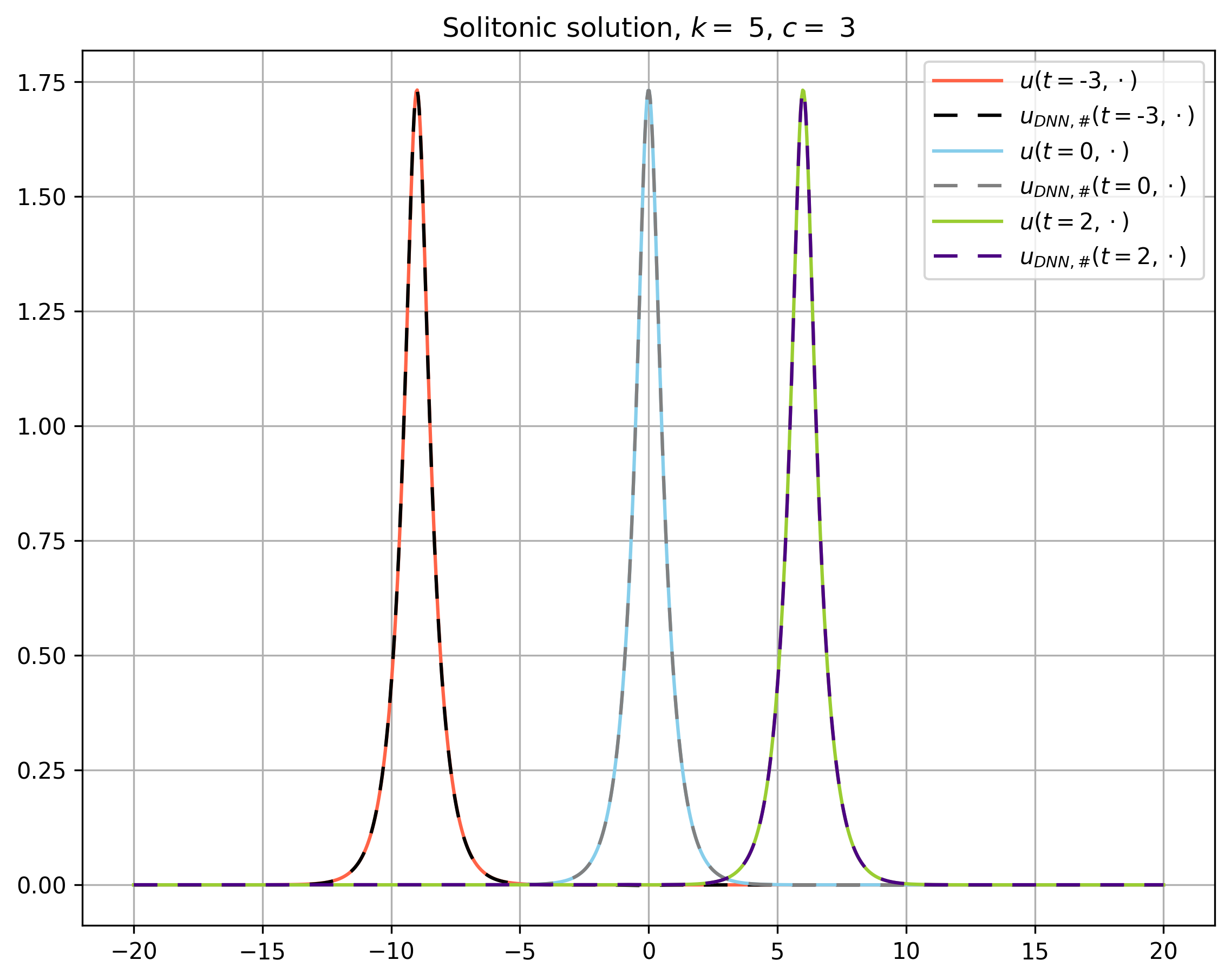}
\caption{$k=5$.}
\label{fig:Soliton-k5c3}
\end{subfigure}
\caption{Exact and predicted solution in the solitonic case, for different values of $k$.}
\label{fig:Soliton}
\end{figure}

First, Figure \ref{fig:Soliton} shows the exact and PINNs approximation in the soliton case with velocity $c=3$, for three different times and for each $k \in \{2,3,4,5\}$. Next, Tables \ref{tab:soliton_c1} and \ref{tab:soliton_c3} summarize the approximation metrics and constants presented in Theorem \ref{MT}, as well as the $L^2$ relative error, for the elected velocities and nonlinearities. The errors presented in Tables \ref{tab:soliton_c1} and \ref{tab:soliton_c3} were made by averaging 5 independent realizations of the algorithm in each setting. From Tables \ref{tab:soliton_c1} and \ref{tab:soliton_c3}, we can conclude that the approximation of solitons was greatly accurate: the $L^2$ relative error in the worst nonlinearity is $\sim0.06 \%$ for velocity $c=1$ and $\sim 0.3 \%$ for velocity $c=3$. For $c=1$ the worst values for the corresponding errors are {\bf error}$_Y\sim3.2 \times 10^{-3}$ and {\bf error}$_{L_t^{\infty}H_x^s} \sim  2.1 \times 10^{-3}$, achieved in the cases $k=4$ and $k=5$, respectively. In the second velocity, $c=3$, the worst values of the errors are {\bf error}$_Y\sim6.7 \times 10^{-2}$ and {\bf error}$_{L_t^{\infty}H_x^s} \sim  2.3 \times 10^{-2}$, both achieved in the $k=3$ scenario. Notice additionally that $\widetilde A$ and the loss $\mathcal L_{k,s}$ remain small and the constants $A,L$ are bounded between 1 and 10. Finally, the maximum time per simulation over all the experiments was less than 3 minutes.

\begin{table}[ht]
\centering
\setlength{\tabcolsep}{6pt}
\renewcommand{\arraystretch}{1.2}
\begin{tabular}{|c||c|c|c|c|}
\hline
$k$ & 2 & 3 & 4 & 5 \\ \hline \hline
$\widetilde A$ & $7.836\times10^{-5}$ & $1.222\times10^{-4}$ & $2.551\times10^{-4}$ & $1.479\times10^{-4}$ \\ \hline
$A$ & 3.539 & 3.529 & 3.346 & 1.647 \\ \hline
$L$ & 2.209 & 2.934 & 3.315 & 1.807 \\ \hline
{\bf error}$_Y$ & $2.715\times10^{-3}$ & $1.827\times10^{-3}$ & $3.242\times10^{-3}$ & $1.725\times10^{-3}$ \\ \hline
{\bf error}$_{L_t^\infty H_x^s}$ & $1.183\times10^{-3}$ & $1.144\times10^{-3}$ & $1.451\times10^{-3}$ & $2.080\times10^{-3}$ \\ \hline
{\bf error}$_{rel}$ & $1.542\times10^{-4}$  & $1.619 \times10^{-4}$ & $2.557\times10^{-4}$ & $6.329\times10^{-4}$  \\ \hline
$\mathcal L_{k,s}$ & $1.627\times10^{-4}$ & $1.614\times10^{-4}$ & $2.550\times10^{-4}$ & $3.821\times10^{-4}$ \\ \hline
Train time [s] & 23.99 & 84.89 & 97.51 & 101.44 \\ \hline
\end{tabular}
\caption{Constants and errors coming from Theorem \ref{MT}, in the case of a soliton with velocity $c=1$.}
\label{tab:soliton_c1}
\end{table}

\begin{table}[ht]
\centering
\setlength{\tabcolsep}{6pt}
\renewcommand{\arraystretch}{1.2}
\begin{tabular}{|c||c|c|c|c|}
\hline
 $k$& 2 & 3 & 4 & 5 \\ \hline \hline
$\widetilde A$ & $6.393\times10^{-4}$ & $1.668\times10^{-3}$ & $1.469\times10^{-3}$ & $4.961\times10^{-4}$ \\ \hline
$A$ & 9.893 & 4.963 & 3.781 & 1.647 \\ \hline
$L$ & 9.191 & 5.449 & 4.676 & 2.488 \\ \hline
{\bf error}$_Y$ & $4.482\times10^{-2}$ & $6.724\times10^{-2}$ & $5.541\times10^{-2}$ & $5.348\times10^{-3}$ \\ \hline
{\bf error}$_{L_t^\infty H_x^s}$ & $1.562\times10^{-2}$ & $2.262\times10^{-2}$ & $1.396\times10^{-2}$ & $4.924\times10^{-3}$ \\ \hline
{\bf error}$_{rel}$ & $5.431\times10^{-4}$ & $2.319\times10^{-3}$ & $2.085\times10^{-3}$ & $1.697\times10^{-3}$ \\ \hline
$\mathcal L_{k,s}$ & $7.464\times10^{-4}$ & $2.294\times10^{-3}$ & $1.767\times10^{-3}$ & $1.820\times10^{-3}$ \\ \hline
Train time [s] & 20.36 & 85.32 & 102.89 & 177.90 \\ \hline
\end{tabular}
\caption{Constants and errors coming from Theorem \ref{MT}, in the case of a soliton with velocity c=3.}
\label{tab:soliton_c3}
\end{table}

\subsubsection{The $N$-Soliton case in integrable models}
Now we consider the more demanding case of the $N$-soliton solution appearing in KdV and mKdV ($k=2,3$). The $N$-Solitons for both equations are given in \cite{Hir71} and \cite{Hir72}, respectively. For the KdV equation the $N$-Soliton has the form
\[
u_{c,\delta}(t,x) := 6 \frac{f\partial_x^2f -(\partial_xf)^2}{f^2},
\]
where $f:=f(t,x)$ is defined as
\[
f = 1 + \sum_{n=1}^N \sum_{{}_NC_n} \left(\prod_{k<\ell}^{(n)} \frac{(\sqrt{c_{i_k}}-\sqrt{c_{i_\ell}})^2}{(\sqrt{c_{i_k}}+\sqrt{c_{i_\ell}})^2}\right) \exp\left(s_{i_1}+\ldots+s_{i_n}\right),
\]
and for all $j \in \{1,\ldots,N\}$, $s_j = \sqrt{c_j} (x-c_jt) - \delta_j$. Here, ${}_NC_n$ indicates all possible combinations of $n$ elements taken from $N$, and $(n)$ indicates all possible combinations of the $n$ elements (with the specified condition $k<\ell$, as indicated). On the other hand, the $N$-Soliton solution for mKdV has the form
\[
u_{c,\delta}(t,x) := 2\sqrt{2} \frac{f\partial_xg-g\partial_xf}{f^2+g^2},
\]
where $f:=f(t,x)$ is defined as
\[
f = \sum_{n=0}^{[N/2]} \sum_{{}_N C_{2n}} \left(\prod_{k<\ell}^{(2n)} \frac{-(\sqrt{c_{i_k}}-\sqrt{c_{i_\ell}})^2}{(\sqrt{c_{i_k}}+\sqrt{c_{i_\ell}})^2}\right) \exp\left(s_{i_1}+s_{i_2}+\ldots+s_{i_{2n}}\right),
\]
and $g:=g(t,x)$ is defined as
\[
g = \sum_{n=0}^{[(N-1)/2]} \sum_{{}_N C_{2n+1}} \left(\prod_{k<\ell}^{(2n+1)} \frac{-(\sqrt{c_{i_k}}-\sqrt{c_{i_\ell}})^2}{(\sqrt{c_{i_k}}+\sqrt{c_{i_\ell}})^2}\right) \exp\left(s_{i_1}+s_{i_2}+\ldots+s_{i_{2n+1}}\right),
\]
and where $[N/2]$ denotes the maximum integer which does not exceed $N/2$. For simplicity in our simulations we only consider the cases of 2-Solitons and 3-Solitons for both KdV and mKdV equations, in the time-space region $[-3,3]\times[-20,20]$. We will consider four different 2-Solitons, with velocities 
\[
\vec{c} \in \{(0.1,0.4),(0.5,1),(0.3,1.8),(1,2)\},
\]
and the 3-Soliton with two velocities $\vec{c} \in \{(0.1,1,2),(0.5,1.5,2)\}$. For these examples, we choose $H=3$ hidden layers in all the experiments, with $n=32$ neurons per hidden layer for the 2-Soliton approximation and $n=40$ for the 3-Soliton. In addition, the number of collocation points were chosen as $M_{\text{evol}}=M_{\text{PDE}} = 32$, $N_{\text{evol}}=128$. The value of $N_{\text{PDE}}$ was settled as 128 for the 2-Soliton with velocities $\vec{c} = (0.1,0.4),(0.5,1)$, and 256 for the other simulations of 2 and 3-Solitons.

\begin{figure}[ht]
\centering
\begin{subfigure}[t]{0.45\textwidth}
\centering
\includegraphics[width=\textwidth]{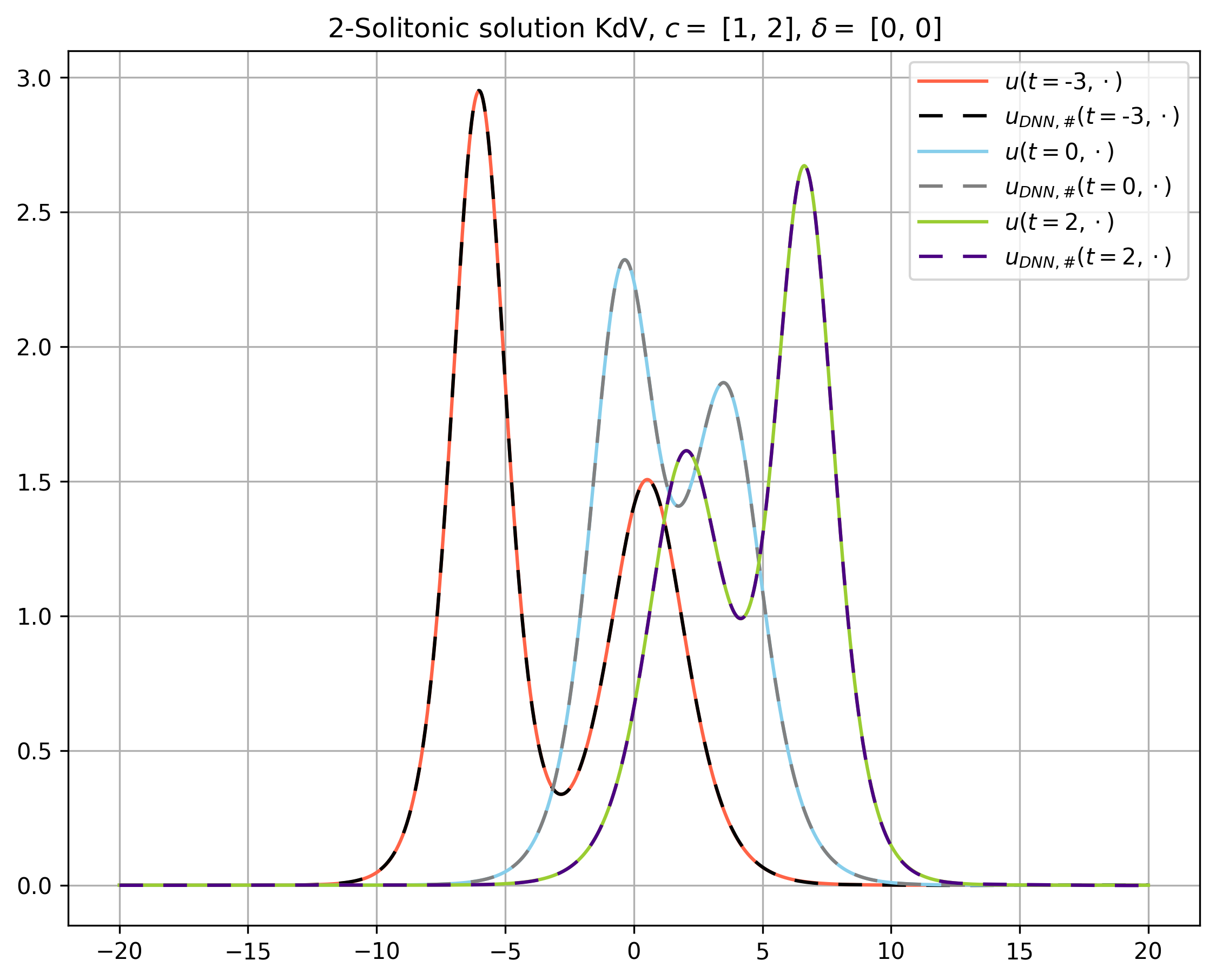}
\caption{$k=2$, $\vec c= (1,2)$.}
\label{fig:2Soliton-k2c4}
\end{subfigure}
\hspace{.2cm}
\begin{subfigure}[t]{0.45\textwidth}
\centering
\includegraphics[width=\textwidth]{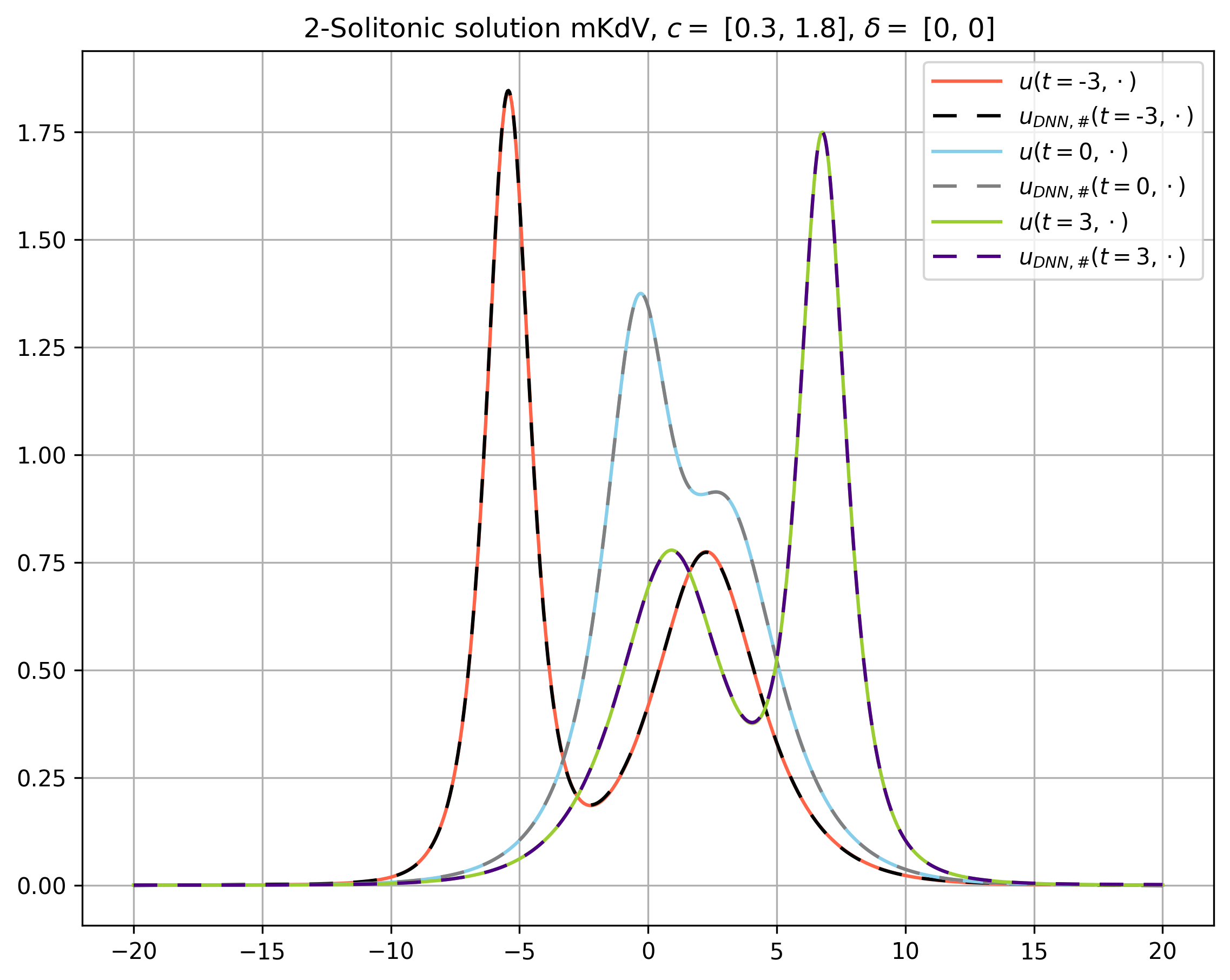}
\caption{$k=3$, $\vec c=(0.3,1.8)$.}
\label{fig:2Soliton-k3c3}
\end{subfigure}
\caption{Exact and predicted solution in the 2-solitonic case: left: $k=2$, right: $k=3$.}
\label{fig:2Soliton}
\end{figure}

\begin{figure}[ht]
\centering
\begin{subfigure}[t]{0.45\textwidth}
\centering
\includegraphics[width=\textwidth]{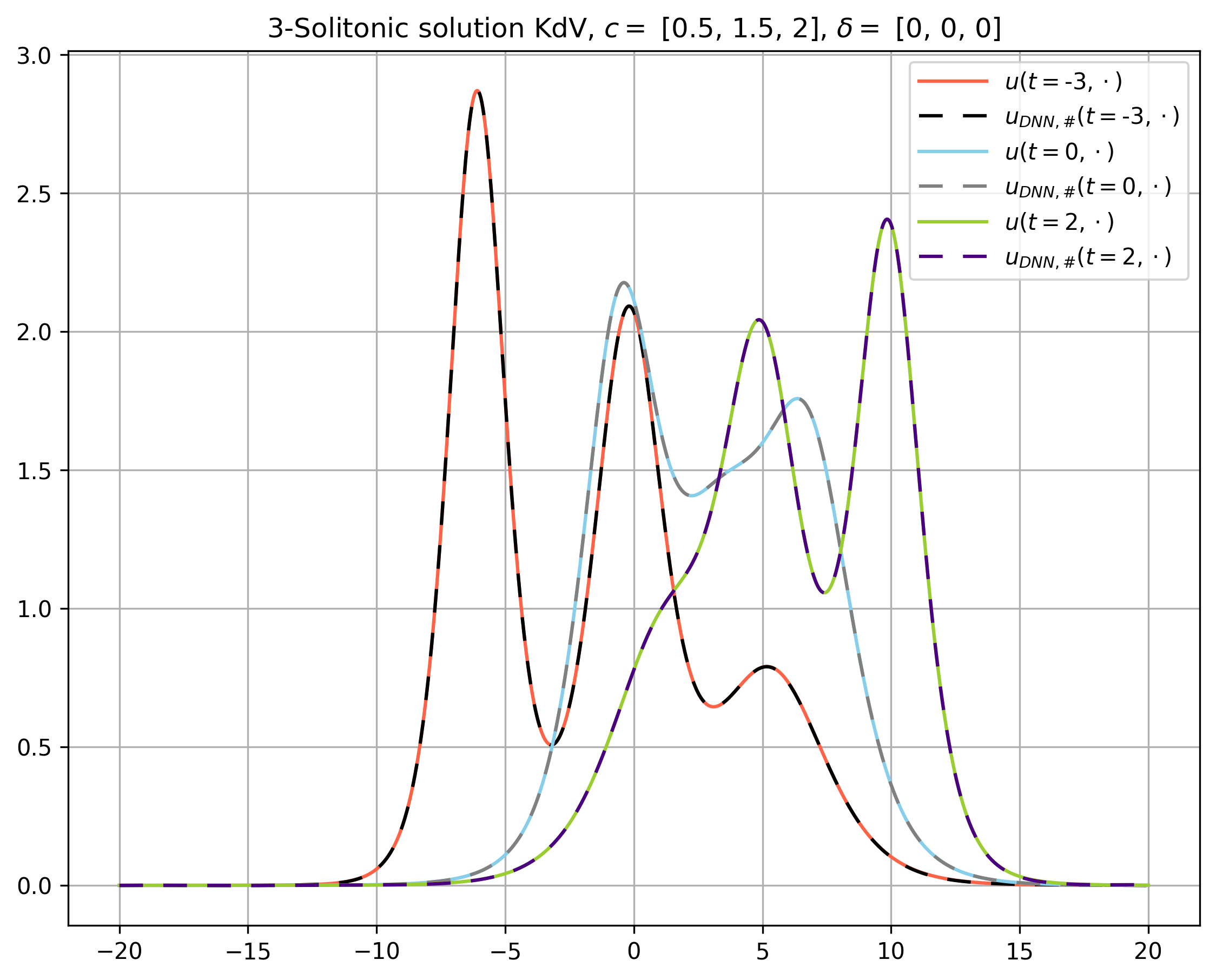}
\caption{$k=2$, $\vec c= (0.5,1.5,2)$.}
\label{fig:3Soliton-k2c2}
\end{subfigure}
\hspace{.2cm}
\begin{subfigure}[t]{0.45\textwidth}
\centering
\includegraphics[width=\textwidth]{images/3Soliton-k3c1.png}
\caption{$k=3$, $\vec c=(0.1,1,2)$.}
\label{fig:3Soliton-k3c1}
\end{subfigure}
\caption{Exact and predicted solution in the 3 solitonic case.}
\label{fig:3Soliton}
\end{figure}

Figure \ref{fig:2Soliton} shows the exact and PINNs approximation for two of the considered settings over three different times. In particular, Figure \ref{fig:2Soliton-k2c4} comprehends the KdV 2-Soliton with velocity $\vec c=(1,2)$, while Figure \ref{fig:2Soliton-k3c3} involves the mKdV 2-Soliton with velocity $\vec c=(0.3,1.8)$. Next, Tables \ref{tab:2soliton_k2} and \ref{tab:2soliton_k3} show the different errors and constants involved in Theorem \ref{MT} for the KdV and mKdV 2-Soliton, respectively. Additionally, Figure \ref{fig:3Soliton} and Table \ref{tab:3soliton} remain the 3-Soliton case. In particular, Figure \ref{fig:3Soliton-k2c2} shows the KdV 3-Soliton with velocity $\vec c = (0.5,1.5,2)$ and Figure \ref{fig:3Soliton-k3c1} encompasses the mKdV 3-Soliton with velocity $\vec c = (0.1,1,2)$. The results presented on Tables \ref{tab:2soliton_k2}, \ref{tab:2soliton_k3} and \ref{tab:3soliton} involve averages of 5 different simulations. 

From Tables \ref{tab:2soliton_k2} and \ref{tab:2soliton_k3} we can obtain the following conclusions: First, the worst $L^2$ relative error is $\sim 0.33\%$ for the KdV model, while it is $\sim 1.4 \%$ in the mKdV case. Additionally, for $k=2$ one has that the worst errors are around {\bf error}$_Y \sim 0.15$, and {\bf error}$_{L_t^\infty H_x^s} \sim 6.3 \times10^{-2}$, while for $k=3$ the worst errors correspond to {\bf error}$_Y \sim 0.5$, and {\bf error}$_{L_t^\infty H_x^s} \sim 0.24$. For both powers of nonlinearity, the loss error $\mathcal L_{k,s}$ and the constant $\widetilde A$ are small of order $O(10^{-2})$ and $O(10^{-3})$, respectively, while the values of the constants $A$ and $L$ are bounded between $0.6$ and $7.6$. The training time per simulation did not take more than $6$ minutes and 30 seconds. The training time for the KdV model is lower than the mKdV case, which involve a more complex nonlinearity.

Table \ref{tab:3soliton} presents a similar behavior: The worst $L^2$ relative error is less than $0.4\%$, achieving relative errors as low as $0.06\%$ in the KdV case, and the worst errors are {\bf error}$_{Y} \sim 0.17$ and {\bf error}$_{L_t^\infty H_x^s} \sim 8.5 \times 10^{-2}$, achieved in the mKdV case with velocity $(0.5,1.5,2)$. Meanwhile, $\mathcal L_{k,s}$ and $\widetilde A$ remain small of order $10^{-3}$ and $10^{-4}$, respectively, and the constants $A$ and $L$ are bounded as well, showing a decreasing trend for fixed velocity $\vec{c}$ when the nonlinearity increases. For this example, the training time per simulation takes around 5 minutes

\begin{table}[ht]
\centering
\begin{tabular}{|l||l|l|l|l|l|}
\hline
$\vec c$ & $\vec c_1$ & $\vec c_2$ & $\vec c_3$ & $\vec c_4$ \\ \hline \hline
$\widetilde A$ & $9.395\times 10^{-5}$ & $2.947\times10^{-4}$ & $5.142\times10^{-4}$  & $6.941\times10^{-4}$ \\ \hline
$A$ & 1.605 & 3.915 & 6.197 & 7.562 \\ \hline
$L$ & 0.612 & 1.890 & 4.136 & 4.902 \\ \hline
{\bf error}$_Y$ & $2.597\times10^{-3}$ & $1.036\times10^{-2}$ & $1.487\times10^{-1}$ & $1.166\times10^{-1}$ \\ \hline
{\bf error}$_{L_t^\infty H_x^s}$ & $1.150\times10^{-3}$ & $4.861\times10^{-3}$ & $6.281\times10^{-2}$ & $4.475\times10^{-2}$ \\ \hline
{\bf error}$_{rel}$ & $2.599\times10^{-4}$ & $5.631\times10^{-4}$ & $3.253\times10^{-3}$ & $2.261\times10^{-3}$ \\ \hline
$\mathcal L_{k,s}$ & $1.315\times10^{-4}$ & $8.442\times10^{-4}$ & $2.524\times10^{-3}$ & $4.944\times10^{-3}$ \\ \hline
Train time [s] & 141.99 & 125.68 & 246.70 & 235.62 \\ \hline
\end{tabular}
\caption{Constants and errors coming from Theorem \ref{MT}, in the case $k=2$, for a 2-Soliton with velocities $\vec c_1=(0.1,0.4)$, $\vec c_2=(0.5,1)$, $\vec c_3 = (0.3,1.8)$, $\vec c_4 = (1,2)$.}
\label{tab:2soliton_k2}
\end{table}

\begin{table}[ht]
\centering
\begin{tabular}{|l||l|l|l|l|l|}
\hline
 $\vec c$ & $\vec c_1$ & $\vec c_2$ & $\vec c_3$ & $\vec c_4$ \\ \hline \hline
$\widetilde A$ & $2.769\times10^{-4}$ & $6.038\times10^{-4}$ & $2.600\times10^{-3}$ & $1.179\times10^{-3}$ \\ \hline
$A$ & 2.859 & 4.302 & 4.701 & 5.462 \\ \hline
$L$ & 0.998 & 2.459 & 3.363 & 3.969 \\ \hline
{\bf error}$_Y$ & $2.279\times10^{-3}$ & $2.346\times10^{-2}$ & $1.025\times10^{-1}$ & $3.456\times10^{-1}$ \\ \hline
{\bf error}$_{L_t^\infty H_x^s}$ & $2.012\times10^{-3}$ & $1.538\times10^{-2}$ & $4.353\times10^{-2}$ & $1.752\times10^{-1}$ \\ \hline
{\bf error}$_{rel}$ & $3.063\times10^{-4}$ & $1.284\times10^{-3}$ & $4.067\times10^{-3}$ & $9.396\times10^{-3}$ \\ \hline
$\mathcal L_{k,s}$ & $5.940\times10^{-4}$ & $1.055\times10^{-3}$ & $6.267\times10^{-3}$ & $7.194\times10^{-3}$ \\ \hline
Train time [s] & 153.63 & 213.41 & 373.36 & 391.23 \\ \hline
\end{tabular}
\caption{Constants and errors coming from Theorem \ref{MT}, in the case $k=3$, for a 2-Soliton with velocities $\vec c_1=(0.1,0.4)$, $\vec c_2=(0.5,1)$, $\vec c_3 = (0.3,1.8)$, $\vec c_4 = (1,2)$.}
\label{tab:2soliton_k3}
\end{table}

\begin{table}[ht]
\centering
\begin{tabular}{|c||c|c||c|c|}
\hline
 $k$& \multicolumn{2}{c||}{2} & \multicolumn{2}{c|}{3} \\ \hline
 $\vec c$ & $\vec c_1$ & $\vec c_2$ & $\vec c_1$ & $\vec c_2$ \\ \hline \hline
$\widetilde A$ & $5.897\times 10^{-4}$ & $5.420\times 10^{-4}$ & $1.084\times 10^{-3}$ & $7.190\times 10^{-4}$ \\ \hline
$A$ & 7.428 & 8.295 & 5.373 & 6.124 \\ \hline
$L$ & 4.207 & 4.400 & 3.215 & 3.711 \\ \hline
{\bf error}$_Y$ & $1.599\times10^{-2}$ & $2.357\times 10^{-2}$ & $4.513\times10^{-2}$ & $1.696\times10^{-1}$ \\ \hline
{\bf error}$_{L_t^\infty H_x^s}$ & $8.015\times10^{-3}$ & $1.136\times10^{-2}$ & $2.501\times10^{-2}$ & $8.540\times10^{-2}$ \\ \hline
{\bf error}$_{rel}$ & $4.542\times10^{-4}$ & $5.900\times10^{-4}$ & $1.553\times10^{-3}$ & $3.758\times10^{-3}$  \\ \hline
$\mathcal L_{k,s}$ & $2.867\times10^{-3}$ & $3.825\times10^{-3}$ & $3.584\times10^{-3}$ & $5.954\times10^{-3}$ \\ \hline
Train time [s] & 295.59 & 306.93 & 308.90 & 303.28 \\ \hline
\end{tabular}
\caption{Constants and errors coming from Theorem \ref{MT}, in the case of a 3-Soliton with velocities $\vec c_1=(0.1,1,2)$, $\vec c_2=(0.5,1.5,2)$. Here we have used $n_\text{iter}=5000$.}
\label{tab:3soliton}
\end{table}

\subsubsection{The breather solution} For the last example we will focus only on the mKdV case, where breathers are admissible solutions. Recall \eqref{breather}. More generally, the mKdV breather \cite{Wad} takes the form
\[
\begin{aligned}
B_{\alpha,\beta}(t,x)&:=2\sqrt{2}\partial_x\left[ \arctan\left(\frac{\beta}{\alpha} \frac{\sin(\alpha(x+\delta t))}{\cosh(\beta(x+\gamma t))}\right) \right], \hspace{.5cm} \alpha,\beta >0. \\
&~= 2\sqrt{2}\beta ~\hbox{sech}(\beta(x+\gamma t)) \\
& \quad \times \left[\frac{\cos(\alpha(x+\delta t))-(\beta/\alpha)\sin(\alpha(x+\delta t))\tanh(\beta(x+\gamma t))}{1+(\beta/\alpha)^2 \sin^2(\alpha(x+\delta t)) ~\hbox{sech}^2(\beta(x+\gamma t))}\right].
\end{aligned}
\]
with $\delta:= \alpha^2-3\beta^2$, and $\gamma:=3\alpha^2-\beta^2$. Here, the velocity of the breather is $-\gamma$.  We considered four different pairs $(\alpha,\beta)$: $(\alpha,\beta) \in \{(0.5,0.5), (0.9,0.3), (1,0.5), (1.3,0.2)\}$. For each case, we train a neural network with $H=3$ hidden layers and $n=40$ neurons per layer. The number of iterations was increased to $n_{iter}=5000$. All the pairs will be approximated in the time space region $[-2,2]\times[-20,20]$, where we have considered the number of collocation points as $N_{\text{evol}}=M_{\text{evol}}=N_{\text{PDE}}=128$, and $M_{\text{PDE}} =64$. 

\begin{figure}[ht]
\centering
\begin{subfigure}[t]{0.45\textwidth}
\centering
\includegraphics[width=\textwidth]{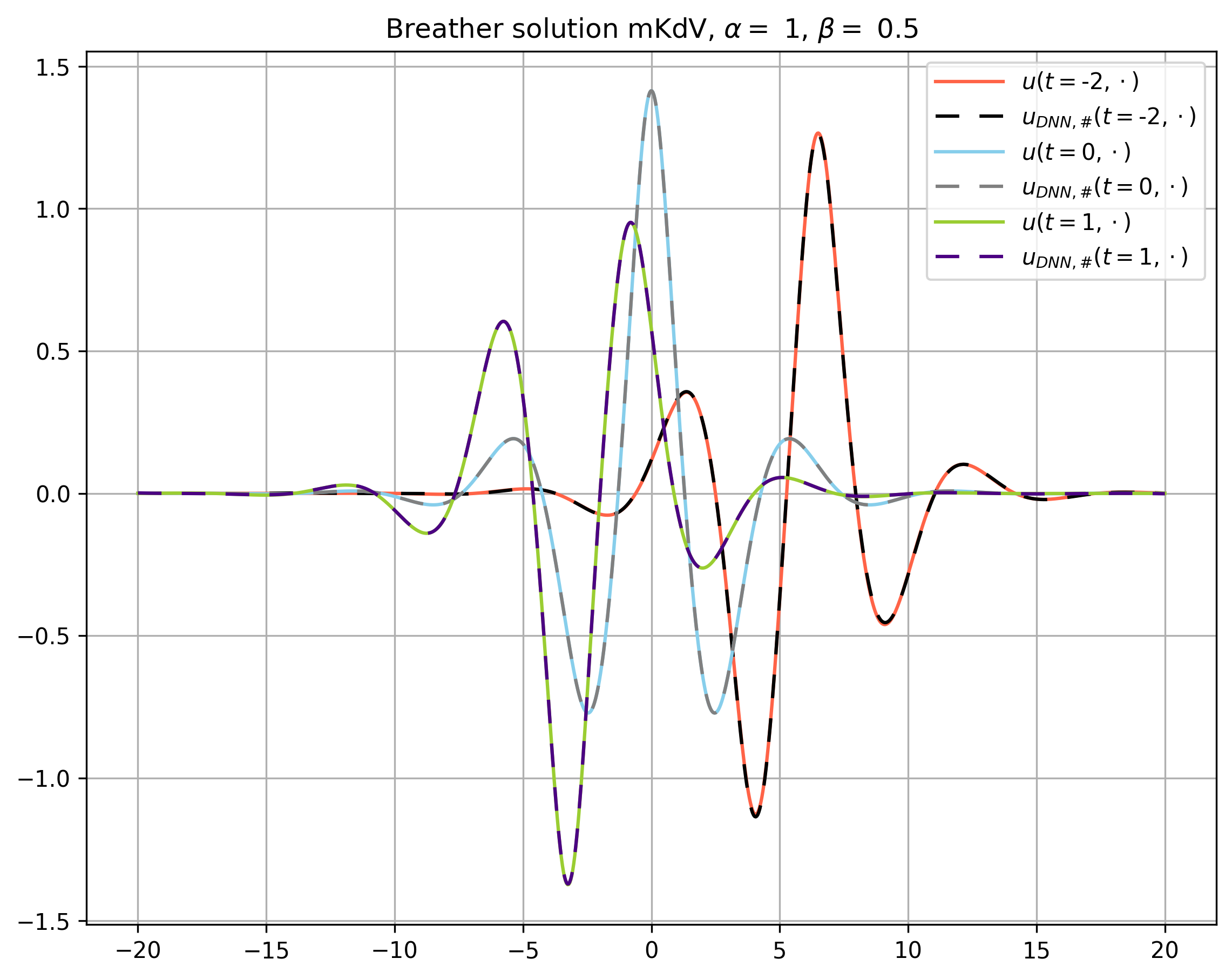}
\caption{$\alpha=1, \beta = 0.5$.}
\label{fig:Breather1}
\end{subfigure}
\hspace{.2cm}
\begin{subfigure}[t]{0.45\textwidth}
\centering
\includegraphics[width=\textwidth]{images/Breather2.png}
\caption{$\alpha=1.3, \beta = 0.2$.}
\label{fig:Breather2}
\end{subfigure}
\caption{Exact and Predicted solution in the breather case.}
\label{fig:Breather}
\end{figure}

Figure \ref{fig:Breather} shows two breather approximation for three different times. In particular, Figure \ref{fig:Breather1} involves the case $\alpha = 1$ and $\beta= 0.5$ and Figure \ref{fig:Breather2} the case $\alpha = 1.3$ and $\beta = 0.2$. Next, Table \ref{tab:breather} comprehends the metrics and constants in Theorem \ref{MT} for the considered pairs $(\alpha,\beta)$. From Table \ref{tab:breather} the following conclusions can be made: First, the worst $L^2$ relative error is $\sim 2.1\%$, achieved in the case $(\alpha,\beta)=(1.3,0.2)$. Meanwhile, the same case of $(\alpha,\beta)$ has the worst {\bf error}$_Y$ and {\bf error}$_{L_t^\infty H_x^s}$ as well, and they are around $0.22$ and $0.13$. The loss error $\mathcal L_{k,s}$ and the constant $\widetilde A$ are small ($\sim 5.3 \times 10^{-3}$ and $\sim 1.5 \times 10^{-3}$, respectively), but higher than previous experiments. In addition, the constants $A,L$ remain bounded between $1.9$ and $4.1$. The training time per simulation is approximate 12 minutes and 30 seconds, longer than the other examples due the augmented number of steps in the LFBGS optimizer and neurons per layer.

\begin{table}[ht]
\centering
\begin{tabular}{|l||l|l|l|l|l|}
\hline
 $(\alpha,\beta)$ & $(0.5,0.5)$ & $(0.9,0.3)$ & $(1,0.5)$ & $(1.3,0.2)$ \\ \hline \hline
$\widetilde A$ & $1.146 \times 10^{-4}$ & $2.750 \times 10^{-4}$ & $6.846\times 10^{-4}$ & $1.530\times 10^{-3}$ \\ \hline
$A$ & 3.859 & 3.078 & 4.061 & 2.617 \\ \hline
$L$ & 3.864 & 2.439 & 3.892 & 1.994 \\ \hline
{\bf error}$_Y$ & $1.355\times 10^{-2}$ & $1.709\times 10^{-2}$ & $9.212\times 10^{-2}$ & $2.160\times10^{-1}$ \\ \hline
{\bf error}$_{L_t^\infty H_x^s}$ & $6.515 \times 10^{-3}$ & $1.234\times 10^{-2}$ & $5.926\times 10^{-2}$ & $1.268\times10^{-1}$ \\ \hline
{\bf error}$_{rel}$ & $8.490 \times 10^{-4}$ & $2.048 \times 10^{-3}$ & $7.182\times 10^{-3}$ & $2.107\times10^{-2}$ \\ \hline
$\mathcal L_{k,s}$ & $7.450 \times 10^{-4}$ & $2.299 \times 10^{-3}$ & $5.346\times10^{-3}$ & $4.260\times10^{-3}$ \\ \hline
Train time [s] & 752.14 & 722.95 & 755.23 & 734.67 \\ \hline
\end{tabular}
\caption{Constants and errors coming from Theorem \ref{MT}, in the case $k=3$, for a breather solution with scaling parameters $\alpha$ and $\beta$ given on each column.}
\label{tab:breather}
\end{table}

\subsubsection{The defocusing mKdV kink case}

\begin{figure}[ht]
\centering
\begin{subfigure}[t]{0.45\textwidth}
\centering
\includegraphics[width=\textwidth]{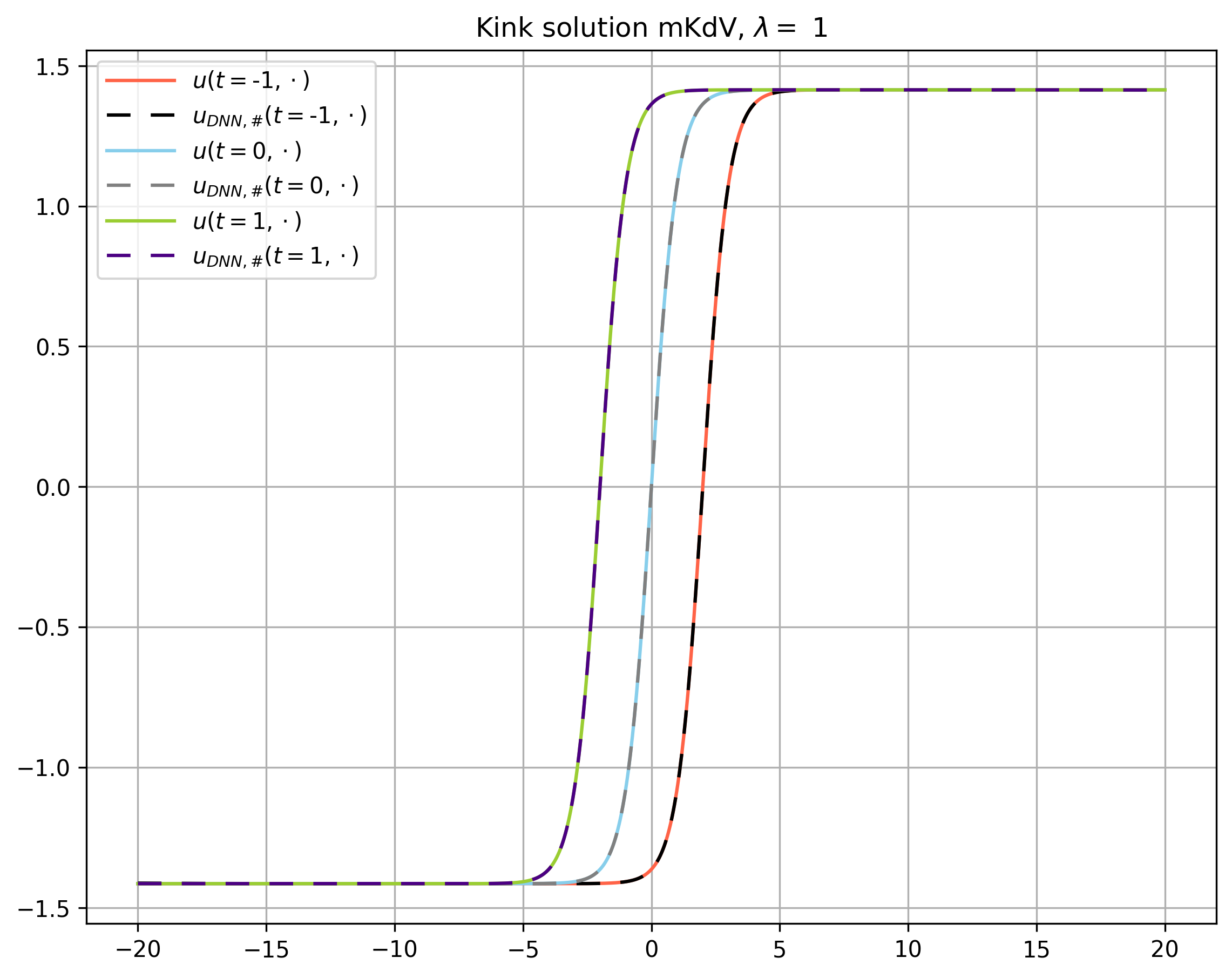}
\caption{Case $\lambda=1$.}
\label{fig:Kink1}
\end{subfigure}
\hspace{.2cm}
\begin{subfigure}[t]{0.45\textwidth}
\centering
\includegraphics[width=\textwidth]{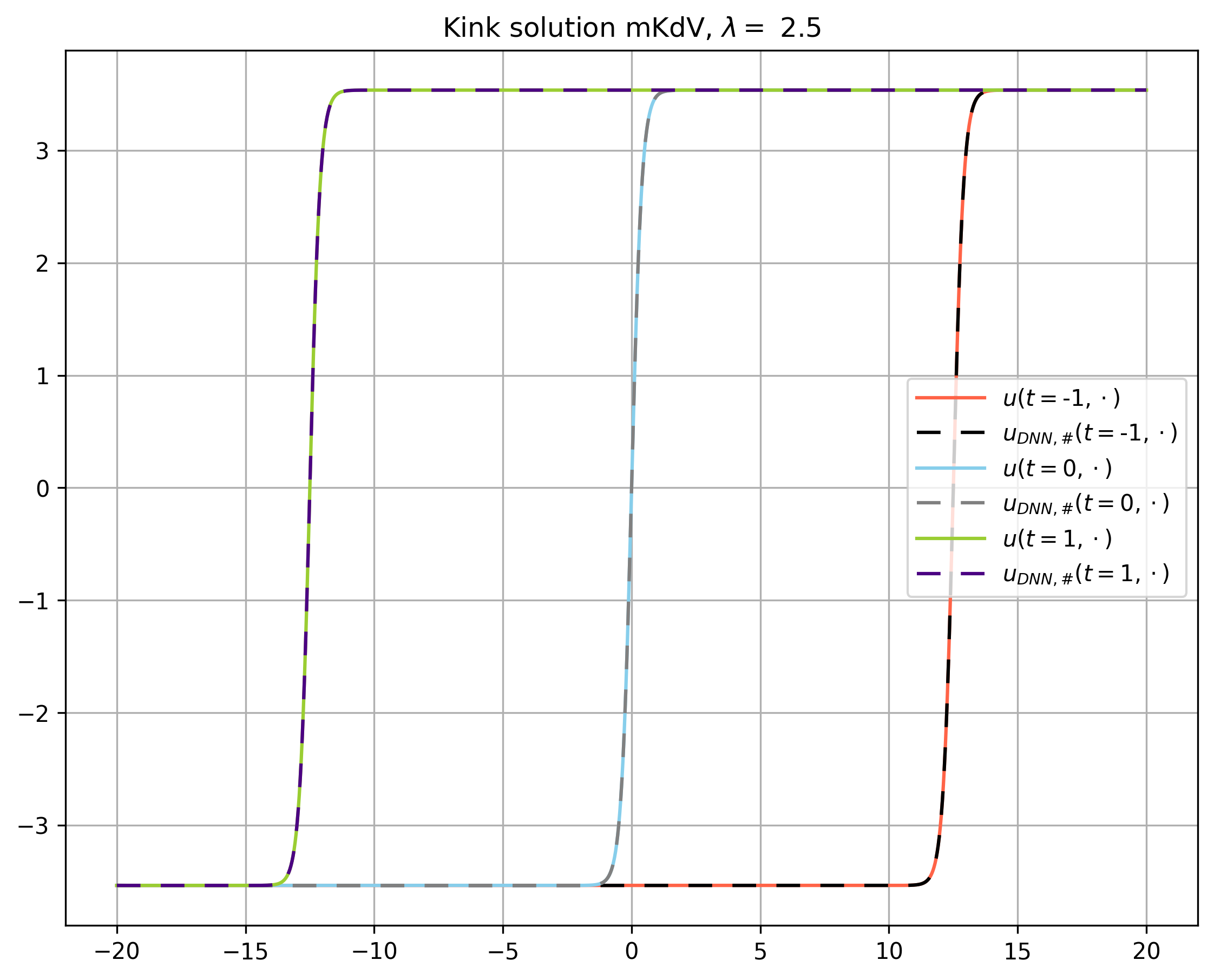}
\caption{Case $\lambda=2.5$.}
\label{fig:Kink2}
\end{subfigure}
\caption{Exact and predicted kink solutions in the kink case \eqref{kink_mKdV}.}
\label{fig:Kink}
\end{figure}

We now consider the case of the defocusing mKdV model
\[
\partial_t u +\partial_{xxx} u = \partial_x(u^3), \quad (t,x)\in \mathbb R\times \mathbb R. 
\]
In this case, we shall consider the case of nonlinear objects that escape from the hypotheses stated in Theorem \ref{MT}, in the sense that solutions are no longer in $H^s(\mathbb R)$. Indeed, for any $\lambda \in\mathbb R$, we consider the kink solution, given by
\be\label{kink_mKdV}
u(t,x) = \sqrt{2} \lambda\tanh(\lambda (x+2\lambda^2t)),
\ee
which is a renowned solution of the defocusing mKdV equation. Although the theoretical results presented in this paper consider only the focusing case, the proposed numerical implementation described in this work allows us to treat the defocusing case, since at the end of the day our main result mainly considers differences among functions. We stress, however, that the kink \eqref{kink_mKdV} does not belong to $H^s(\R)$ (it has non-zero boundary values at $\pm\infty$), so it lies strictly outside the hypotheses of Theorem \ref{MT}; the quantities reported below are therefore not constants and errors furnished by Theorem \ref{MT}, but rather the \emph{same approximation metrics} $A,\widetilde A,L$, ${\bf error}_Y$, ${\bf error}_{L^\infty_tH^s_x}$ and ${\bf error}_{rel}$, evaluated for this reference solution. They are meaningful because the residual-minimization framework and the local-in-$x$ control of the difference $u-u_{\text{DNN},\#}$ underlying our scheme remain well-defined for bounded, smooth data such as \eqref{kink_mKdV}, even when the global $H^s$-theory does not formally apply. First, Fig. \ref{fig:Kink} presents the approximation of the kink solution for three different times between $[-1,1]$ and two different values of $\lambda$: $1$ and $2.5$ (Figs. \ref{fig:Kink1} and \ref{fig:Kink2}, respectively). Then, Table \ref{tab:kink} shows the metrics of Theorem \ref{MT} for the case of the kink solution. In particular, as one can see in Table \ref{tab:kink}, the relative errors are near the error $O(10^{-5})$, while the error associated to the operator $Y_{k,s}$ does not surpass the value $2\times10^{-3}$, and then leading to very good results.

\begin{table}[ht]
\centering
\begin{tabular}{|l|l|l|}
\hline
 $\lambda$ & $1$ & $2.5$ \\ \hline \hline
$\widetilde A$ & $2.303 \times 10^{-4}$ & $1.779 \times 10^{-3}$  \\ \hline
$A$ & 14.790 & 37.760  \\ \hline
$L$ & 4.788 & 11.290  \\ \hline
{\bf error}$_Y$ & $2.739\times 10^{-3}$ & $1.405\times 10^{-2}$  \\ \hline
{\bf error}$_{L_t^\infty H_x^s}$ & $2.499 \times 10^{-3}$ & $1.124\times 10^{-2}$  \\ \hline
{\bf error}$_{rel}$ & $5.731 \times 10^{-5}$ & $1.251 \times 10^{-4}$  \\ \hline
$\mathcal L_{k,s}$ & $4.757 \times 10^{-4}$ & $1.495 \times 10^{-2}$  \\ \hline
Train time [s] & 197.05 & 267.99 \\ \hline
\end{tabular}
\caption{Approximation metrics (defined as in Theorem \ref{MT}, but here computed for a solution outside its hypotheses) in the case $k=3$, for a kink solution with scaling parameter $\lambda$ given at each column. Notice the reasonable quality in the errors found despite the small computation time.}
\label{tab:kink}
\end{table}

\subsection{Evolution over the number of iterations}
In this section we will present how the estimated quantities obtained by the PINN algorithm evolve through the number of inner iterations of the LBFGS optimizer. In particular we consider only the mKdV model with initial conditions given by the soliton and 2-soliton solutions. The other presented examples behave in a similar way. 

\begin{figure}[ht]
\centering
\begin{subfigure}[t]{0.45\textwidth}
\centering
\includegraphics[width=\textwidth]{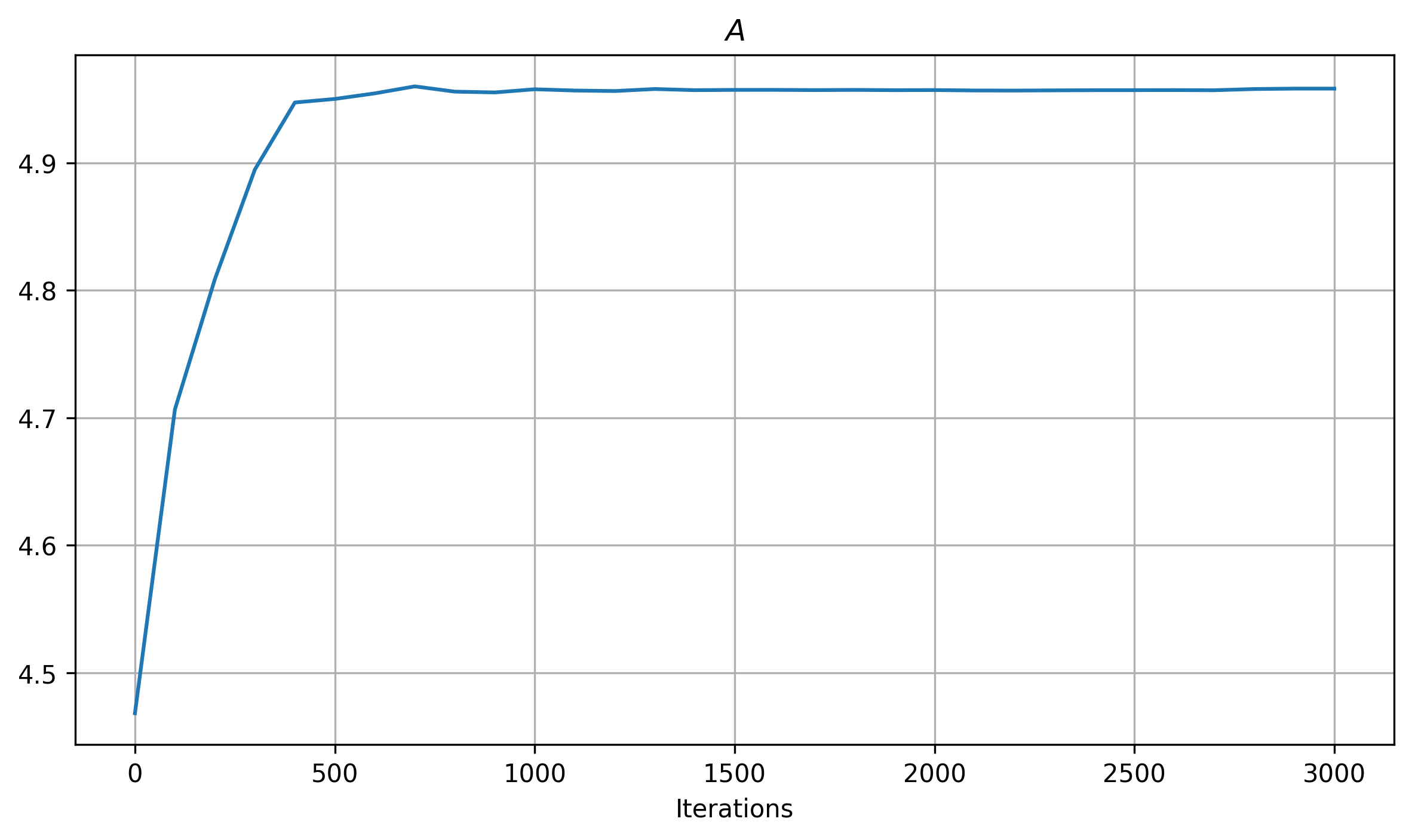}
\caption{$A$.}
\label{fig:A-Soliton}
\end{subfigure}
\hspace{.2cm}
\begin{subfigure}[t]{0.45\textwidth}
\centering
\includegraphics[width=\textwidth]{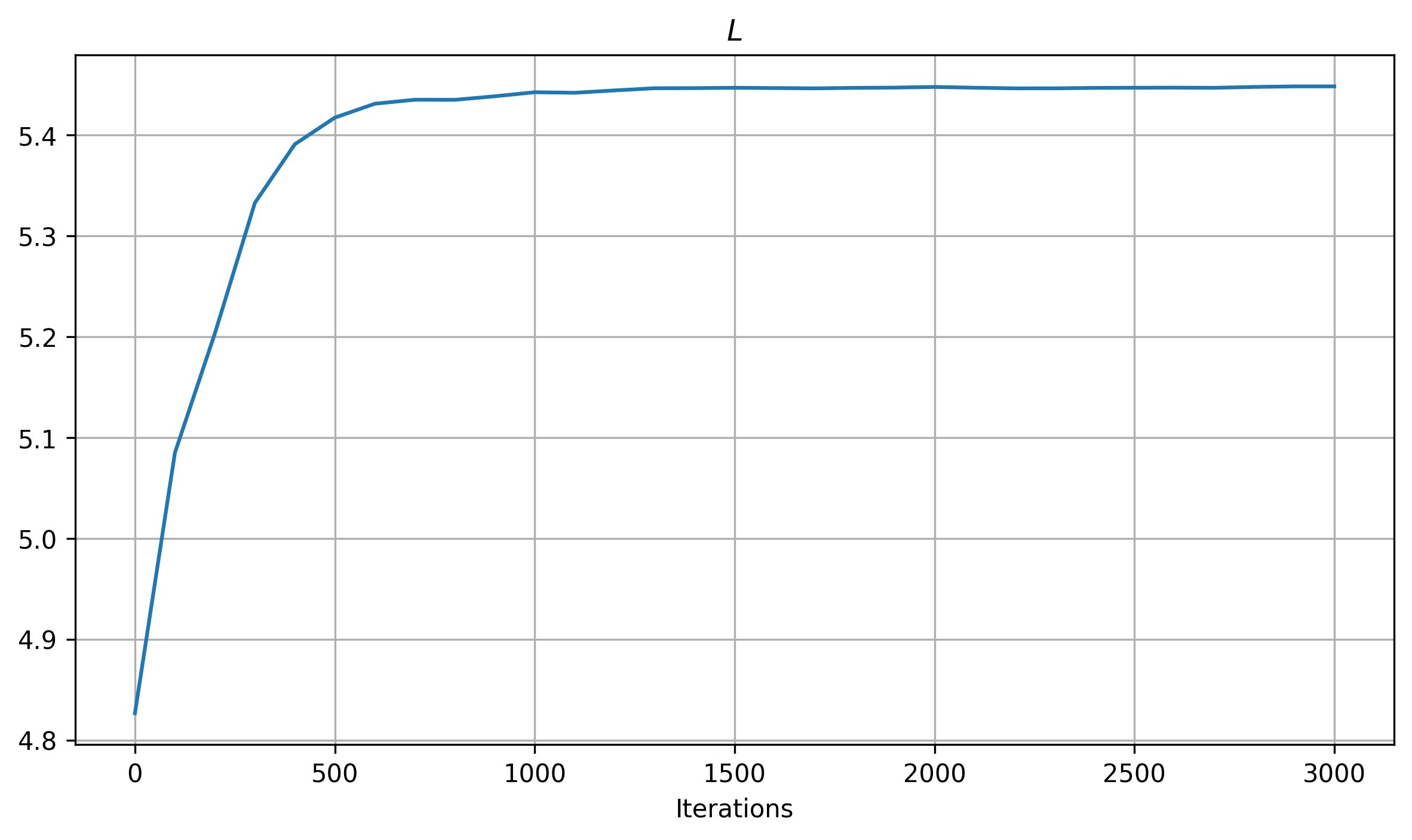}
\caption{$L$.}
\label{fig:L-Soliton}
\end{subfigure}
\begin{subfigure}[t]{0.45\textwidth}
\centering
\includegraphics[width=\textwidth]{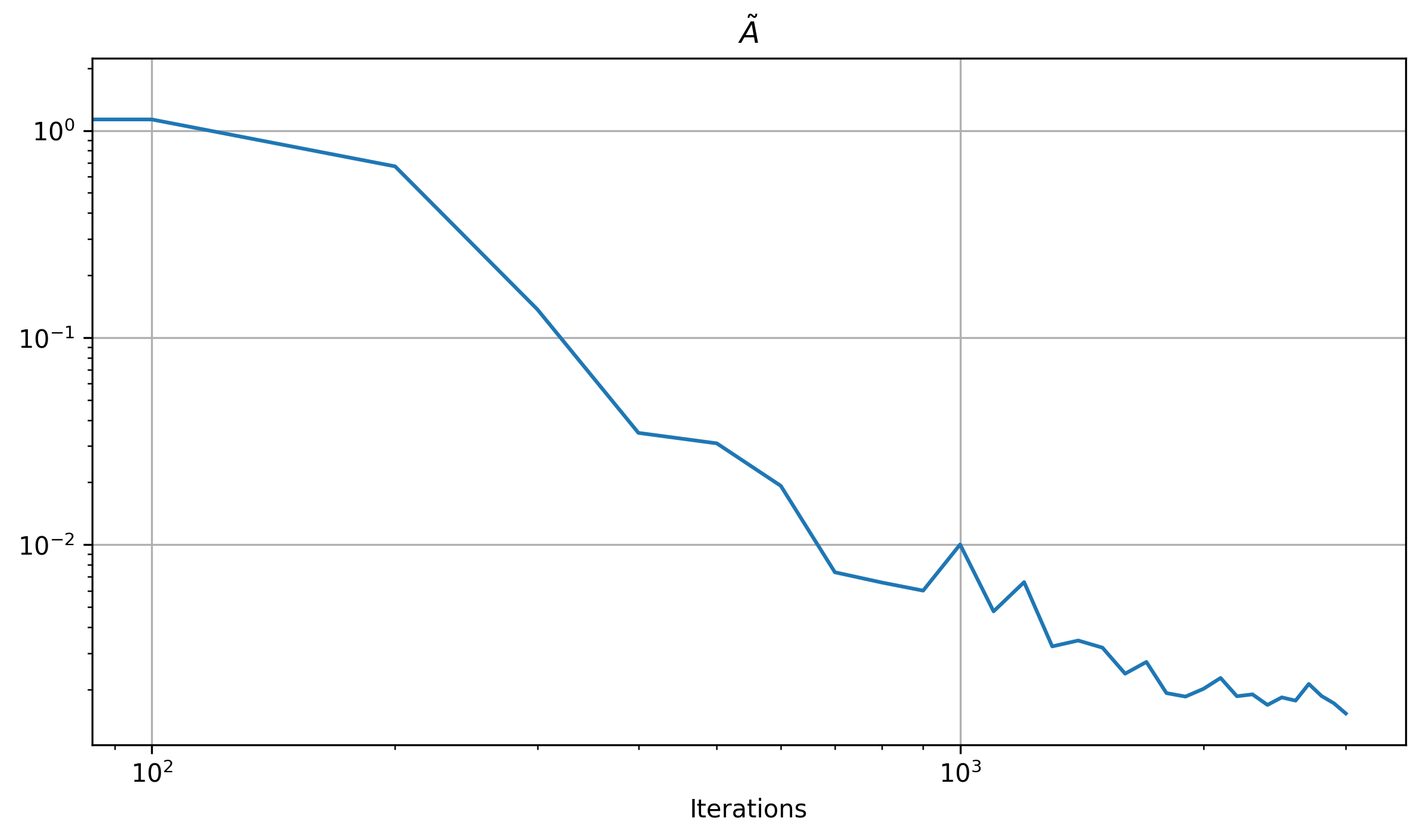}
\caption{$\widetilde A$.}
\label{fig:tildeA-Soliton}
\end{subfigure}
\hspace{.2cm}
\begin{subfigure}[t]{0.45\textwidth}
\centering
\includegraphics[width=\textwidth]{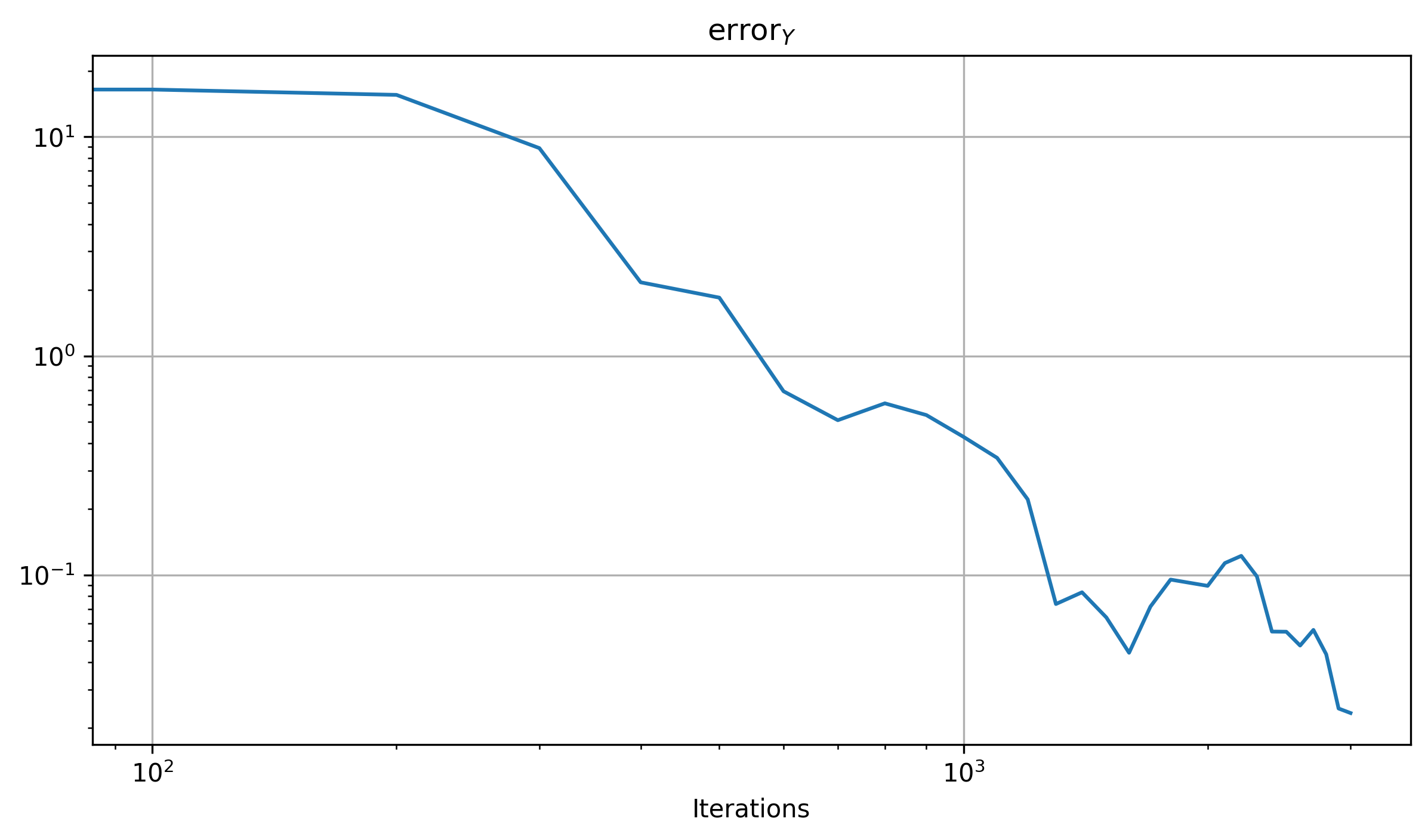}
\caption{{\bf error}$_Y$.}
\label{fig:Y-Soliton}
\end{subfigure}
\begin{subfigure}[t]{0.45\textwidth}
\centering
\includegraphics[width=\textwidth]{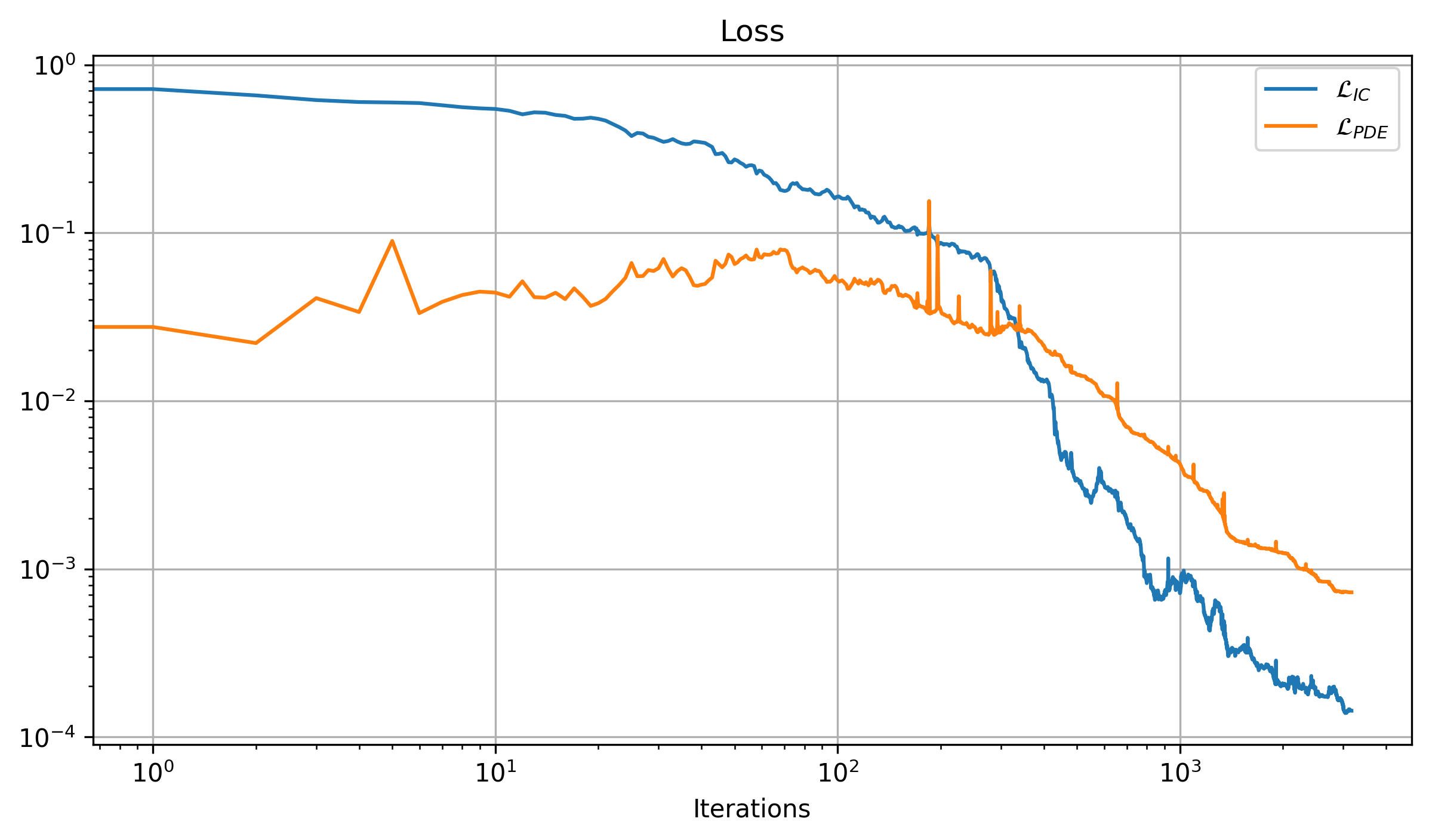}
\caption{$\mathcal L_{\text{evol}}$ (blue) and $\mathcal L_{\text{PDE}}$ (orange).}
\label{fig:Loss-Soliton}
\end{subfigure}
\caption{Estimation of quantities involved in Theorem \ref{MT} through the number of iterations, in the soliton regime with $c=3$. In \eqref{fig:A-Soliton}, \eqref{fig:L-Soliton}, the constants $A$ and $L$ are computed in terms of the number of iterations of the algorithm. Same procedure in panels \eqref{fig:tildeA-Soliton} and \eqref{fig:Y-Soliton}. Finally, in \eqref{fig:Loss-Soliton} the loss functions are depicted.}
\label{fig:Soliton-iters}
\end{figure}

\begin{figure}[ht]
\centering
\begin{subfigure}[t]{0.45\textwidth}
\centering
\includegraphics[width=\textwidth]{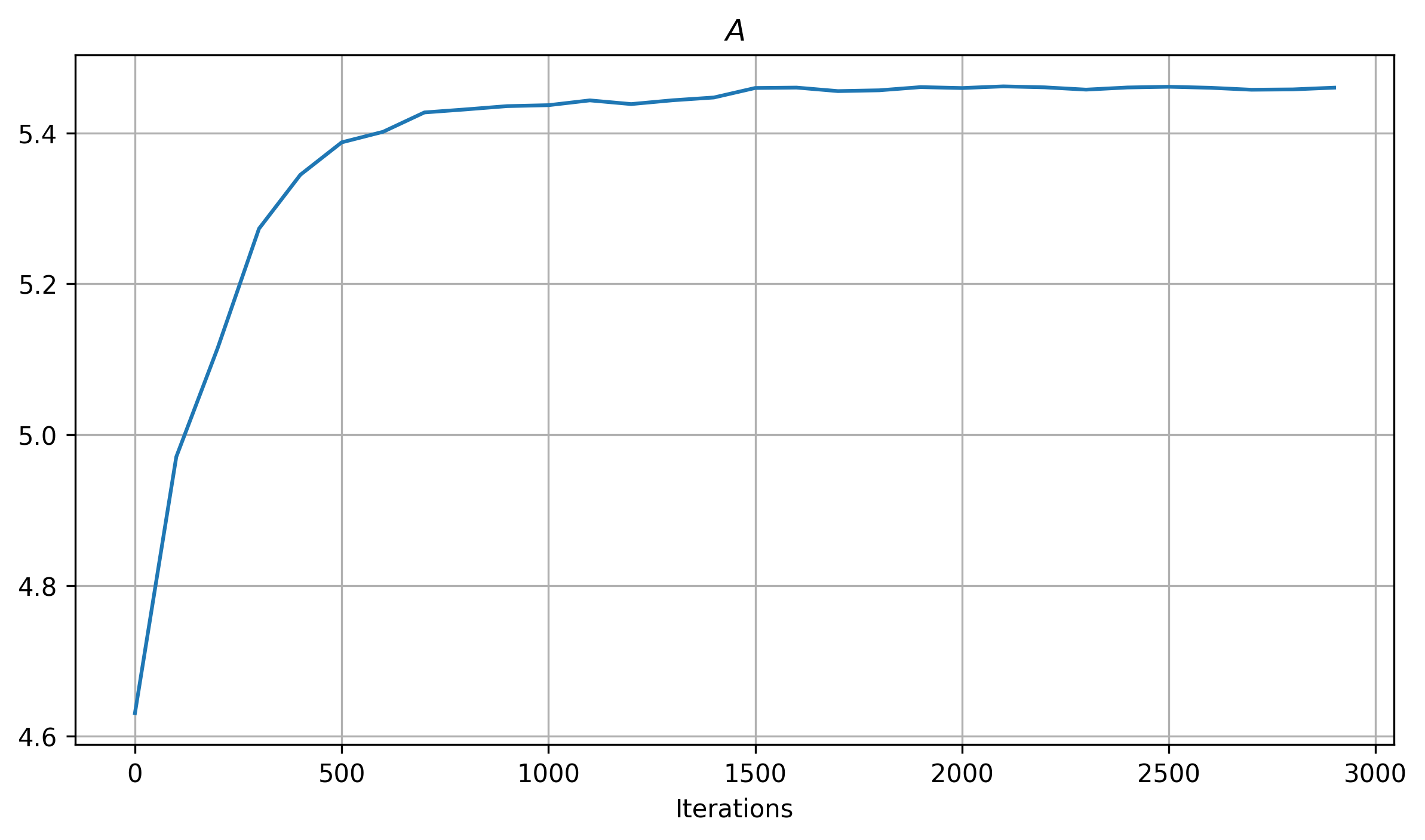}
\caption{$A$.}
\label{fig:A-2Soliton}
\end{subfigure}
\hspace{.2cm}
\begin{subfigure}[t]{0.45\textwidth}
\centering
\includegraphics[width=\textwidth]{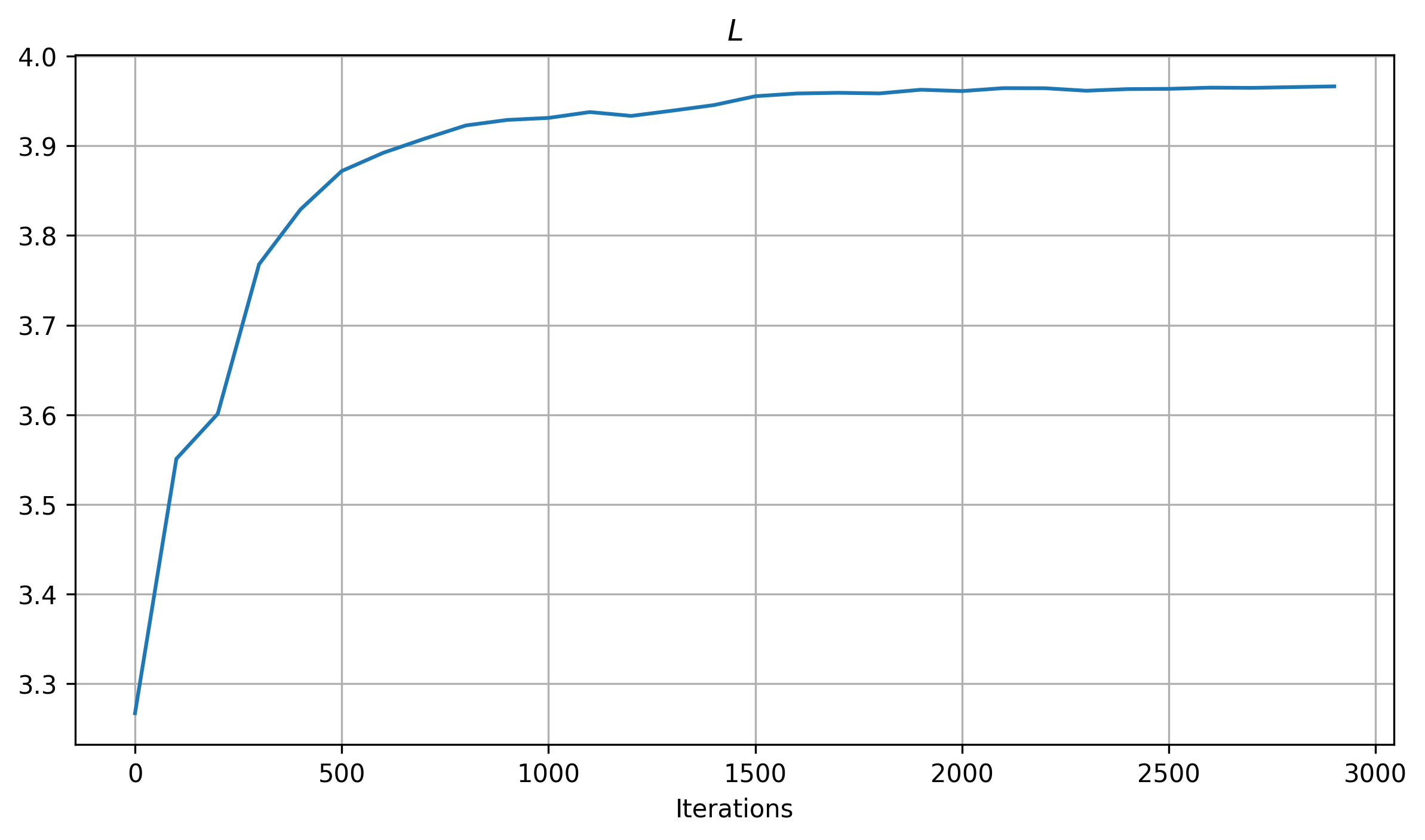}
\caption{$L$.}
\label{fig:L-2Soliton}
\end{subfigure}
\begin{subfigure}[t]{0.45\textwidth}
\centering
\includegraphics[width=\textwidth]{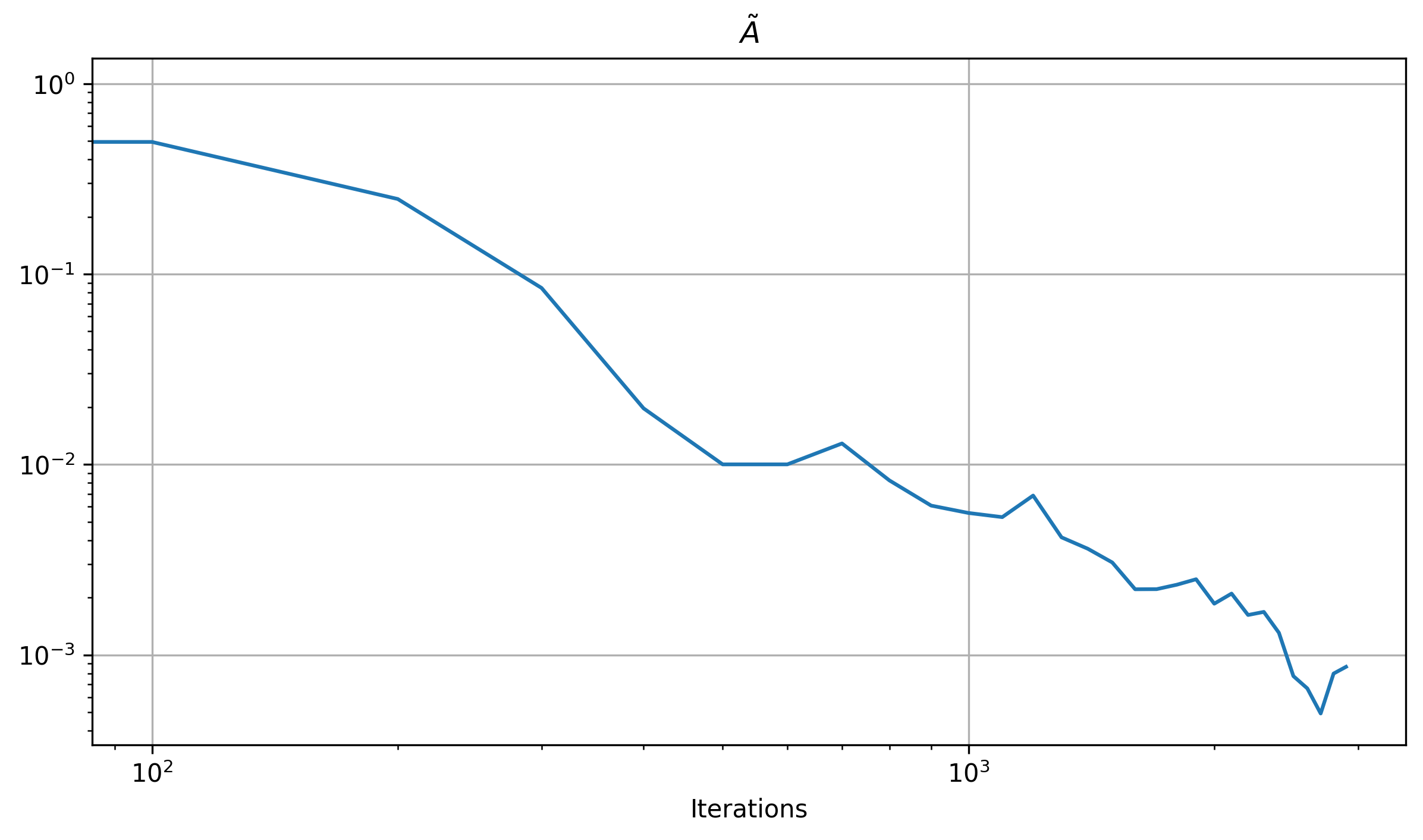}
\caption{$\widetilde A$.}
\label{fig:tildeA-2Soliton}
\end{subfigure}
\hspace{.2cm}
\begin{subfigure}[t]{0.45\textwidth}
\centering
\includegraphics[width=\textwidth]{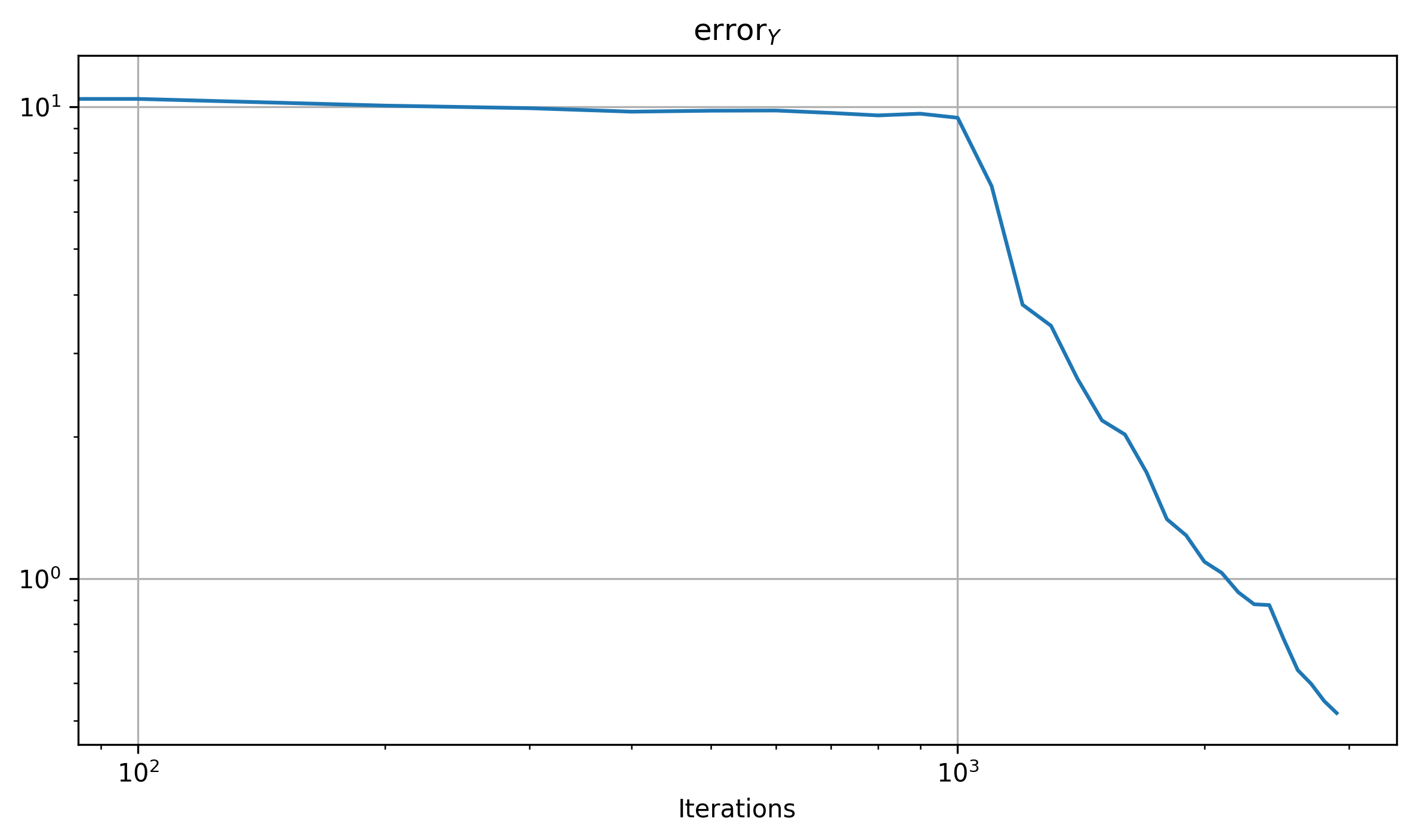}
\caption{{\bf error}$_Y$.}
\label{fig:Y-2Soliton}
\end{subfigure}
\begin{subfigure}[t]{0.45\textwidth}
\centering
\includegraphics[width=\textwidth]{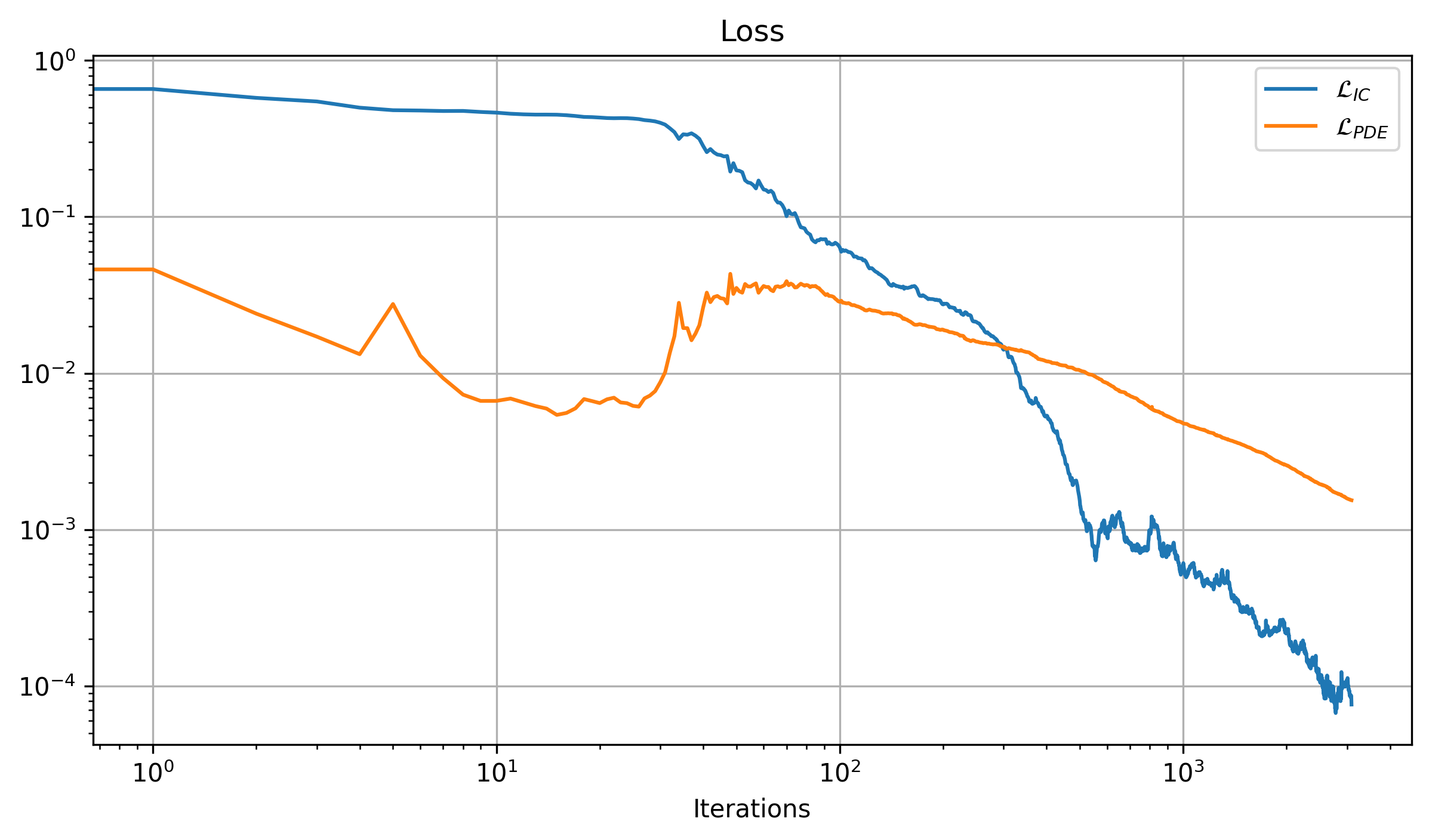}
\caption{$\mathcal L_{\text{evol}}$ (blue) and $\mathcal L_{\text{PDE}}$ (orange).}
\label{fig:Loss-2Soliton}
\end{subfigure}
\caption{Estimation of quantities involved in Theorem \ref{MT} through the number of iterations, in the 2-soliton regime with $\vec c=(1,2)$.}
\label{fig:2Soliton-iters}
\end{figure}

First, Figure \ref{fig:Soliton-iters} shows the evolution of the constants $A,\widetilde A,L$ and the errors $\mathcal L_{k,s}$ and {\bf error}$_Y$ through the number of iterations for the soliton solution with velocity $c=3$, and Figure \ref{fig:2Soliton-iters} presents the evolution of constants and errors obtained with the PINNs algorithm in the 2-soliton case with velocity $\vec c=(1,2)$. In particular, Figures \ref{fig:A-Soliton}, \ref{fig:L-Soliton}, \ref{fig:A-2Soliton} and \ref{fig:L-2Soliton}, show that both constants $A$ and $L$ have an abrupt variation on the first iterations (500 in Figures \ref{fig:A-Soliton}, \ref{fig:L-Soliton} and 1500 in Figures \ref{fig:A-2Soliton}, \ref{fig:L-2Soliton}), but after this threshold they get stabilized and the variation is minimal. Contrarily, Figures \ref{fig:Y-Soliton} and \ref{fig:Y-2Soliton} show that the change on {\bf error}$_Y$ is very small in the first 300 and 1000 iterations, respectively, but then they drop. This is not the case of the evolution of $\widetilde A$, as shown in Figures \ref{fig:tildeA-Soliton} and \ref{fig:tildeA-2Soliton}, where their values constantly go down. Lastly, Figures \ref{fig:Loss-Soliton} and \ref{fig:Loss-2Soliton} show the evolution of both losses appearing in $\mathcal L_{k,s}$. For both cases, the algorithm struggles to lower the evolution part of the loss function $\mathcal L_{\text{evol},k,s}$, but after a certain threshold, the decay is notable. Although the value of $\mathcal L_{\text{PDE},k,s}$ starts lower than $\mathcal L_{\text{evol},k,s}$, its descent begins almost at the same time than $\mathcal L_{\text{evol},k,s}$.

%

\section{Proof of Theorem \ref{MT}: preliminaries}\label{sec:2}

\subsection{Linear estimates} To present the results of the paper, we start by presenting some preliminary results that will be important in the development of the work. We recall the definition of the unitary group $ \{e^{-t\partial_x^3}\}$ associated to the solution of the linear problem
\begin{equation}
\begin{cases}
    \partial_t v + \partial_x^3 v= 0,\\
    v(x,0) = v_0(x),
\end{cases}    
\end{equation}
The function \( v(x,t) \) is given by
\begin{equation}
v(x, t) = e^{-t\partial_x^3}v_0 (x) = S_t \ast v_0 (x),
\end{equation}
where
\begin{equation}
S_t (x) = \int_{-\infty}^{\infty} e^{2\pi i x \xi} e^{8\pi^3 i t \xi^3} d\xi.
\end{equation}
We also consider the associated inhomogeneous problem
\begin{equation}
\partial_t v + \partial_x^3 v = f, \quad (x, t) \in \mathbb{R} \times \mathbb{R},
\end{equation}
with initial condition $v(x, 0) = 0$, whose solution is given by the Duhamel formula $v(x,t)=\int_0^t e^{-(t-t')\partial_x^3}f(\cdot,t')(x)\,dt'$. Recall that $D_x^s := \mathcal F^{-1} (|\xi|^s \mathcal F(\cdot))$. We now collect several space-time estimates for the Airy group $e^{-t\partial_x^3}$ acting on suitable functions, which will be used throughout the remainder of this work.

\begin{lemma}\label{Linear estimatess}
Let $v_0=v_0(x)\in L^2(\R)$ and let $f=f(t,x)$, $g=g(t,x)$ be functions on $\R\times\R$ for which the right-hand sides below are finite. The group $ \{e^{-t\partial_x^3}\}$ satisfies
  \begin{equation}
 \|e^{-t\partial_x^3}v_0\|_{L^5_x L^{10}_t} + \|\partial_x e^{-t\partial_x^3}v_0\|_{L_x^\infty L_t^2} \leq C \|v_0\|_{L^2_x},
\end{equation}
\begin{equation}
\|D_xe^{-t\partial_x^3}v_0\|_{L_x^{20}L_t^{5/2}}+ \|e^{-t\partial_x^3}v_0\|_{L^4_x L^\infty_t} \leq C \|D_x^{1/4} v_0\|_{L^2_x}.
\end{equation}
\begin{equation}\label{linear 4-kdv}
\left\| \int_{-\infty}^{\infty} \partial_x e^{-(t - t')\partial_x^3} f (t',x) dt' \right\|_{L_x^\infty L_t^2} \leq C \|f\|_{L_x^1 L_t^2}.
\end{equation}
Furthermore, for any \( ( \theta,\alpha) \in [0,1] \times [0, 1/2] \), we have
\begin{equation}
\|D_x^{\theta \alpha/2}e^{-t\partial_x^3}v_0\|_{L_{t}^{q}L_{x}^{p}} \leq C \| v_0\|_{L^2_x},
\end{equation}
\begin{equation}
    \left\| \int_{-\infty}^{\infty}D_x^{\theta \alpha} e^{-(t - t')\partial_x^3} g(t',x) \, dt' \right\|_{L^q_{t} L^{p}_{x}} \leq C \|g\|_{L^{q'}_{t} L^{p'}_{x}}.
\end{equation}
and
\begin{equation}
    \left\| \int_{0}^{t}D_x^{\theta \alpha} e^{-(t - t')\partial_x^3} g(t',x) \, dt' \right\|_{L^q_t([-T,T]) L^{p}_x} \leq C \|g\|_{L^{q'}_t([-T,T]) L^{p'}_x},
\end{equation}
where \( (q, p) \) satisfies the relation
\[
(q, p) = \left( \frac{6}{\theta(\alpha + 1)}, \frac{2}{1 - \theta} \right),
\]
and the exponents satisfy the duality condition $\frac{1}{q} + \frac{1}{q'} = \frac{1}{p} + \frac{1}{p'}  = 1$. In particular,  if \( g \in L^{5/4}_x L^{10/9}_t \), then
\begin{equation}
\left\| \int_{0}^{t} e^{-(t - t')\partial_x^3} g(t',x) \, dt' \right\|_{L^5_x L^{10}_t} \leq C \|g\|_{L^{5/4}_x L^{10/9}_t}.
\end{equation}
\end{lemma}

The previous estimates are classical in the literature; we refer e.g. to Kenig, Ponce and Vega \cite{KPV1} and to Linares and Ponce \cite{LP}. 

\subsection{Product and chain rules} We also require below the use of classical product and chain rules.

\begin{lemma}[Product rule] \label{fraLR}
Given  $s>0.$ Let $f=f(x)$ and $g=g(x)$ be functions on $\R^d$, let $1<r<\infty$, and let $1<p_1,p_2,q_1,q_2  \leq \infty$
with $$\frac{1}{p_j}+\frac{1}{q_j}=\frac{1}{r}.$$
Here $J^s:=\mathcal F^{-1}\big((1+|\xi|^2)^{s/2}\mathcal F(\cdot)\big)$ denotes the inhomogeneous Bessel-potential operator. Then,
\begin{equation}
\begin{split}
  \left\|D^s(fg)\right\|_{L^r(\mathbb{R}^d)} &\lesssim  \left\|D^s f\right\|_{L^{p_1}(\mathbb{R}^d)} \left\|g\right\|_{L^{q_1}(\mathbb{R}^d)}+ \left\|f\right\|_{L^{p_2}(\mathbb{R}^d)} \left\|D^sg\right\|_{L^{q_2}(\mathbb{R}^d)} \label{libRule} \\
\end{split}
    \end{equation}
    and 
    \begin{equation}
    \begin{split}
     \left\|J^s(fg)\right\|_{L^r(\mathbb{R}^d)} & \lesssim  \left\|J^s f\right\|_{L^{p_1}(\mathbb{R}^d)} \left\|g\right\|_{L^{q_1}(\mathbb{R}^d)}+ \left\|f\right\|_{L^{p_2}(\mathbb{R}^d)} \left\|J^s g\right\|_{L^{q_2}(\mathbb{R}^d)}. \label{eqfraLR}
    \end{split}
    \end{equation}    
\end{lemma}

\begin{lemma}[Improved product rule]\label{G.N.S.0}
    Let $s\in (0,1)$ and  $ s_1, s_2 \in [0,\alpha]$ with $s = s_1 +s_2. $ Let $p,q,p_1,p_2,q_1,q_2 \in(1,\infty)$. Let $f=f(t,x)$ and $g=g(t,x)$ be functions on $\R\times\R$, and recall the notation $L^q_T:=L^q_t([-T,T])$, so that $\|f\|_{L^p_xL^q_T}=\big(\int_\R(\int_{-T}^T|f(t,x)|^q\,dt)^{p/q}dx\big)^{1/p}$. Then
    \begin{equation}\label{Est. space-time}
        \|D_x^s(fg)-fD^s_xg-gD^s_xf \|_{L^p_xL^q_T}~{} \leq~{} C\|D^{s_1}_xf\|_{L^{p_1}_xL^{q_1}_T}\|D^{s_2}_xg\|_{L^{p_2}_xL^{q_2}_T}
    \end{equation}
    where $\frac{1}{p}= \frac{1}{p_1} +  \frac{1}{p_2} $ and $\frac{1}{q}=\frac{1}{q_1} + \frac{1}{q_2}$.
\end{lemma}

Let us also recall the fractional Gagliardo--Nirenberg inequality.

\begin{lemma}\label{G.N.S.}
    Let $s > 0$ and $\rho \in (0,s)$ and  $p,q,r \in(1,\infty)$, then 
    \begin{equation}\label{Embbeding}
        \|D^{\rho}(f)\|_{L^p}~{} \leq~{} C\|f\|_{L^{r}}^{1-\rho/s}\|D^{s}f\|_{L^{q}}^{\rho/s}
    \end{equation}
    with 
    $$
    \frac{1}{p} = \left(1 - \frac{\rho}{s}\right)\frac{1}{r} + \left( \frac{\rho}{s}\right)\frac{1}{q}. 
    $$
\end{lemma}



\subsection{Local well-posedness theory}
Recall the norms $\|\cdot\|_{Y_{k,s}}$ defined in equations \eqref{eq:Y2s_norm}-\eqref{eq:Y5s_norm} for $k\in\{2,3,4,5\}$ and $s \geq 0$ . The results in next sections are followed by the ideas of Kenig, Ponce and Vega \cite{KPV1}. For the sake of completeness we will enunciate their results, adapted to our notation. They are summarized in the following Theorem
\begin{theorem}[LWP gKdV] Let $s \geq s_k,$ if $k\in\{3,4,5\},$ and $s > s_k,$ if $k = 2$. Then for any $u_0 \in H^s(\R)$ there exists $T=T(\|u_0\|_{H^s(\R)})>0$ (with $T(\rho)\to \infty$ as $\rho \to 0$) and a unique strong solution $u(t)$ of the IVP \eqref{eq:KdV} satisfying
\begin{equation}
    u \in C([-T,T];H^s(\R)) \quad \mbox{and} \quad \|u\|_{Y_{k,s}([-T,T]\times\R)} < \infty.
\end{equation}
The only exception is the $L^2$-critical case $k=5$ at the endpoint $s=s_5=0$: by \cite[Theorem 2.10]{KPV1} the existence time then depends on the profile of the initial datum itself, $T=T(u_0)$, and not merely on $\|u_0\|_{L^2_x}$; for $s>0$ in the case $k=5$ the dependence $T=T(\|u_0\|_{H^s})$ is recovered through \cite[Corollary 2.11]{KPV1}. For $k\in\{2,3,4\}$ the dependence $T=T(\|u_0\|_{H^s})$ holds as stated, in accordance with \cite{KPV1}.
\end{theorem}
\begin{remark}
This general version of the theorem above for the cases $k\in\{3,5\}$ can be found on \cite{LP}.
\end{remark}

\section{Long time stability results in g-KdV models}\label{sec:3}

In the following Lemmas, we will consider $I$ be a fixed compact interval containing zero, for example $I = [-T,T]$. We start this section by recalling the definition introduced in \eqref{def_E}.

\begin{definition}[Approximate solution]\label{def:appr}
	We say that $u_{\#}$ is an approximate solution to the gKdV model~\eqref{eq:KdV} with nonlinearity $k \in \{2,3,\ldots\}$ if $u_{\#}$ satisfies the perturbed equation
	\begin{equation}\label{eq:perturbed}
	\partial_t u_{\#} + \partial_{xxx} u_{\#} + \partial_x(u^k_\#)= \mathcal{E}_k[u_{\#}].
\end{equation} 
for some error function $\mathcal{E}_k.$
\end{definition}
Let us first describe a long-time stability result in the cases $k=2,3,4,5$.
\begin{theorem}[Long time stability, cases $k=2,3,4,5$]\label{thm:long}
    Let \( I \) be a time interval containing zero, $T:= \sup_{t\in I}|t|$, $s \geq s_k,$ if $k\in\{3,4,5\},$ and $s > s_k,$ if $k = 2$. Let \( u_{\#} \) be an approximate solution to \eqref{eq:KdV}, on \( I \times \mathbb{R} \) in the sense of Definition \ref{def:appr}. Assume that
\begin{align}
       \|u_{\#}\|_{L^\infty_tH^s_x  (I\times \R)} & \leq A, \label{eq:long1}\\
       \|u_{\#}\|_{Y_{k,s}(I\times\R)}&\leq L, \label{eq:long2}
     \end{align}
for some constants $A, L>0.$ Let $u_0\in H^{s}_x(\R)$ such that
\begin{equation}\label{eq:long3}
    \|u_0-u_{\#}(0)\|_{H^{s}_x(\R)}\leq \widetilde{A},
\end{equation}
for some $\widetilde{A}>0.$  Assume the smallness conditions
\begin{align}
   \|e^{-t\partial_x^3}(u_0-u_{\#}(0))\|_{Y_{k,s}(I\times\R)} & \leq \varepsilon, \label{eq:long4}\\
   \left\|\mathcal{E}_k[u_{\#}]\right\|_{L^2_t H^{s}_x (I\times \R)}&\leq \varepsilon, \quad \mbox{if} \quad k \leq 4,\label{eq:long5}\\
    \left\|\partial^{-1}_x\mathcal{E}_k[u_{\#}]\right\|_{L^1_x L^2_t (I\times \R)}&\leq \varepsilon, \quad \mbox{if} \quad k = 5. \label{eq:long6}
\end{align}
for some $0<\varepsilon\leq \varepsilon_1,$ with $\varepsilon_1=\varepsilon_1(A, \widetilde{A},T)$ a small constant and $\mathcal E_k$ defined in \eqref{eq:perturbed}. Then there exists a solution $u \in C(I;H^s(\R))$ to~\eqref{eq:KdV} on $I\times \R$ with initial datum $u_0\in H^{s}_x(\R)$ such that
\begin{align}
\label{L2L2}
\|\partial_x(u^k-u^k_{\#})\|_{L_x^2L_t^2(I\times\R)}\, & \leq\, C(A, \widetilde{A},L)\varepsilon,\\
\label{Y general}
    \|u-u_\#\|_{Y_{k,s}(I\times\R)}\,& \leq\, C(A, \widetilde{A},L)\varepsilon\\
 \label{final gkdv Hs}
 \|u-u_\#\|_{L_t^{\infty}H_x^s(I\times\R)}\, & \leq\, C(A, \widetilde{A},L)(1 +  \varepsilon)\\
 \label{final u norm gkdv}
  \|u\|_{Y_{k,s}(I\times\R)}\,& \leq\, C(A, \widetilde{A},L)
\end{align}
\end{theorem}
The previous result encompasses the $L^2$ critical case $k=5$, which has been already proved  in \cite{KKSV}, that for sake of completeness we will enunciate here
\begin{theorem}[Long-time stability $L^2$ critical gKdV] Let $I$ be a time interval containing zero and let $u_\#$ be an approximate solution \eqref{eq:KdV} on $I \times \R$ in the sense of Definition \eqref{def:appr}. Assume that 
\[
\|u_\#\|_{L_t^\infty L_x^2(I\times \R)} \leq A,
\]
for some constant $A>0$. Let $u_0 \in L_x^2(\R)$ such that
\[
\|u_0-u_\#(0)\|_{L_x^2(\R)} \leq \widetilde A,
\]
for some $\widetilde A >0$. Assume
\[
\begin{aligned}
    \|u_\#\|_{L_x^5 L_t^{10}(I\times \R)} &\leq B, \\
    \|e^{-t\partial_x^3}(u_0 - u_\#(0))\|_{L_x^5 L_t^{10}(I\times \R)} &\leq \varepsilon, \\
    \|D_x^{-1} \mathcal E_5[u_\#]\|_{L_x^1 L_t^2(I\times\R)} &\leq \varepsilon,
\end{aligned}
\]
for some constant $B>0$ and some small $0 < \varepsilon < \varepsilon_1 = \varepsilon_1(A,\widetilde A,B)$. Then there exists a solution $u$ to \eqref{eq:KdV} on $I \times \R$ with initial data $u_0$ at time $t=0$ satisfying
\[
\begin{aligned}
    \|u-u_\#\|_{L_x^5L_t^{10}(I\times\R)} &\leq C(A,\widetilde A,B) \varepsilon, \\
    \|u^5-u_\#^5\|_{L_x^1L_t^2(I\times\R)} &\leq C(A,\widetilde A,B)\varepsilon, \\
    \|u-u_\#\|_{L_t^{\infty}L_x^2(I\times\R)} + \|D_x^{1/6}(u-u_\#)\|_{L_{t,x}^6(I\times\R)} &\leq C(A,\widetilde A,B), \\
    \|\partial_x u\|_{L_x^\infty L_t^2(I\times\R)} + \|D_x^{-1/4}u\|_{L_x^4L_t^\infty(I\times\R)} + \|u\|_{L_t^\infty L_x^2(I\times\R)} & \notag \\
    +\|D^{1/6}_x u\|_{L_{t,x}^6(I\times \R)} + \|u\|_{L_x^5L_t^{10}(I\times\R)} &\leq C(A,\widetilde A,B).
\end{aligned}
\]
\end{theorem}

\subsection{Improved estimates} We next state the lemmas on stability over small intervals of the real line, which will be used in the proof of the main results. In order to proceed, we shall work with the auxiliary norms defined through the following notation. We assume $s>s_k$ defined in \eqref{s_k}. If $k=2$,
\begin{equation}\label{eq:tildeY2s_norm}
\begin{aligned}
\|u\|_{\widetilde Y_{2,s}(I\times\R)}&:= \|\partial_x u\|_{L_t^4L_x^{\infty}(I\times\R)} 
+ \|D_x^s\partial_x u\|_{L_x^{\infty}L_t^2(I\times\R)}  \\
& +  |I|^{1/4}\|u\|_{L_x^2L_t^\infty(I\times\R)} 
+ |I|^{1/4}\|u\|_{L_t^{\infty}H^s(I\times\R)},
\end{aligned}
\end{equation}
while if $k=3$,
\begin{equation}\label{eq:tildeY3s_norm}
\begin{aligned}
\|u\|_{\widetilde Y_{3,s}(I\times\R)} &:= \|\partial_x u\|_{L_x^\infty L_t^2(I\times\R)} 
 + \|D_x^s u\|_{L_x^5L_t^{10}(I\times\R)} 
+ \|D^s_x \partial_x u\|_{L_x^\infty L_t^2(I\times\R)} \\
&~{} \quad + |I|^{1/4}\|u\|_{L_x^4L_t^\infty(I\times\R)} 
+ \|\partial_x u\|_{L_x^{20}L_t^{5/2}(I\times\R)}; 
\end{aligned}
\end{equation}
and if $k=4,$
\begin{equation}\label{eq:tildeY4s_norm}
\begin{aligned}
\|u\|_{\widetilde Y_{4,s}(I\times\R)}&:=\|  D_x^{s}\partial_x u\|_{L_x^\infty L_T^2} +  \|u \|_{ L_x^{\frac{42}{13}}L_T^{\frac{21}4} } +  \|u \|_{L_x^{\frac{60}{13}}L_T^{15} }\\
 &~{}\quad + \| u\|_{ L_x^{\frac{10}{3}}L_T^{\frac{30}7} } +  \| D_x^{s} u\|_{ L_x^{\frac{10}{3}}L_T^{\frac{30}7} } +  \| \partial_x u\|_{L_x^\infty L_T^2}.
\end{aligned}
\end{equation}

 where $\rho > 3/4$ is a fixed constant. In the case of $k=5$, the $\widetilde Y_{5,s}$ norm will coincide with the $Y_{k,5}$ norm. Now we are able to enunciate a version of Theorem \ref{thm:long} over a short time regime.


\begin{lemma}[Short-time stability gKdV]\label{short time general}
Let I be a fixed compact interval containing zero with $T \in (0,1]$. Let $u_{\#}$ be an approximate solution to ~\eqref{eq:KdV} on $I\times \R$ in the sense of Definition~\ref{def:appr}. Assume that 
\begin{equation*}
    \|u_{\#}\|_{L^\infty_t H^s_x (I\times \R)}\leq A,
\end{equation*}
for some constant $A>0$ and $s \geq s_k,$ if $k\in\{3,4,5\},$ and $s > s_k,$ if $k = 2.$
Let $u_0\in H^{s}_x(\R)$,  such that
\[
    \|u_0-u_{\#}(0)\|_{H^{s}_x(\R)}\leq \widetilde{A},
\]
for some $\widetilde{A}>0.$
 Assume the smallness conditions
\begin{align}
    \label{short gkdvug}
   \|u_{\#}\|_{\widetilde Y_{k,s}(I\times\R)}&\leq \varepsilon_0, \\
   \label{short gkdvflux}
   \|e^{-t\partial_x^3}(u_0-u_{\#}(0))\|_{Y_{k,s}(I\times\R)} & \leq \varepsilon,\\
   \label{short gkdvE}
    \left\|\mathcal{E}_k[u_{\#}]\right\|_{L^2_t H^{s}_x (I\times \R)}&\leq \varepsilon, \quad \mbox{if} \quad k \leq 4,\\
    \left\|\partial^{-1}_x\mathcal{E}_k[u_{\#}]\right\|_{L^1_x L^2_t (I\times \R)}&\leq \varepsilon, \quad \mbox{if} \quad k = 5.
\end{align}
for some $0<\varepsilon\leq \varepsilon_0,$ with $\varepsilon_0=\varepsilon_0(A, \widetilde{A})$ a small constant. 

Then there exists a solution $u$ to~\eqref{eq:KdV} on $I\times \R$ with initial datum $u_0\in H^{s}_x(\R)$ such that 
\begin{align}
\|\partial_x(u^k-u^k_{\#})\|_{L_x^2L_t^2(I\times\R)} & \lesssim \varepsilon,\\
  \label{Yks final}
  \|u -u_{\#}\|_{Y_{k,s}(I\times\R)}& \lesssim \varepsilon,\\
  \label{short gkdv Hs}
 \|u-u_\#\|_{L_t^{\infty}H_x^s(I\times\R)} & \lesssim \widetilde{A} +  \varepsilon,\\
  \|u\|_{Y_{k,s}(I\times\R)} &  \lesssim A + \widetilde{A} +    \varepsilon_0. \label{u norm mkdv}
\end{align}

\end{lemma}

\begin{remark}
    The condition imposed $\|u_{\#}\|_{L^\infty_t H^s_x (I\times \R)}\leq A$ on Theorem \ref{thm:long} and Lemma \ref{short time general} is intended to ensure that the norms appearing in the inequalities can be bounded by $\varepsilon_0$, given that they can be controlled by a constant $CA$ through the well-posedness theory.
\end{remark}

\section{Short time stability: proofs}\label{sec:4}

In this section we prove the short time stability estimates needed for the proof of Theorem \ref{MT}, and the basis for the numerical approximations run in Section \ref{sec:numerical}. 
Let $w := u-u_\#$. Then $w$ satisfies the following initial-value problem
\[
\begin{cases}
\partial_t w + \partial_x^3 w = - \partial_x ((w + u_\#)^k - u_\#^k) - \mathcal E_k[u_\#], \\
w(0) = u_0 - u_\#(0).
\end{cases}
\]
Passing to the integral representation of the problem via Duhamel’s principle
\begin{equation}\label{Duhamel}
 \begin{aligned}
    w(t) = &~{} e^{-t\partial_x^3}(u_0 - u_\#(0)) \\
    &~{} - \int_0^t e^{-(t-s)\partial_x^3}\partial_x ((w + u_\#)^k - u_\#^k)ds - \int_0^t e^{-(t-s)\partial_x^3}\mathcal E_k[u_\#]ds.
    \end{aligned}
\end{equation}
The idea here is for sake of completeness, to present the proof for each case of gKdV. For $t \in I$ and each $k\in\{2,3,4\}$, we write
\begin{equation}\label{Bk estimates}
    B_k(t) := \|\partial_x ((w+u_\#)^k - u_\#^k)\|_{L^2_t H^{s}_x([0,t]\times \R)}.
\end{equation}
and 
\begin{equation}\label{BkI estimates}
    B_k(I) := \|\partial_x ((w+u_\#)^k - u_\#^k)\|_{L^2_t H^{s}_x(I\times \R)}.
\end{equation}
Since $(w+u_\#)^k-u_\#^k=\sum_{j=1}^{k}\binom{k}{j}w^j u_\#^{k-j}$, every monomial in this expansion (and hence in its $\partial_x$ derivative) contains at least one factor $w$; in particular the pure term $u_\#^{k-1}\partial_x u_\#$ coming from $\partial_x(u_\#^k)$ has already been cancelled algebraically by the subtraction of $u_\#^k$, before taking any norm. Therefore
\begin{equation}\label{Bk cota}
\begin{aligned}
B_k(t) \lesssim &~{}\sum_{j=0}^{k-1} \|w^ju_\#^{k-1-j} \partial_x w\|_{L^2_t H^{s}_x} +\sum_{j=1}^{k-1}\|w^ju_\#^{k-1-j} \partial_x u_\#\|_{L^2_t H^{s}_x},
\end{aligned}
 \end{equation}
 where both sums now run only over terms that carry at least one factor $w$ or $\partial_x w$.
We may assume without loss of generality that $\inf I = 0$. To simplify the notation, we will use $Y_{2,s}:=Y_{2,s}(I\times \R)$ and  $\widetilde Y_{2,s}:= \widetilde Y_{2,s}(I\times \R).$
\begin{proof}[Short-time stability in the KdV case]
Let us consider $k=2.$ Firstly, we will bound each term involved in $\|w\|_{Y_{2,s}}$, using Duhamel formula, the linear estimates in Lemma \ref{Linear estimatess} and Lemma \ref{G.N.S.},
\medskip 

\noindent 
$\bullet$ \underline{$\|\partial_xw\|_{L_t^4L_x^{\infty}}$}:
\begin{align*}
\|\partial_xw\|_{L_x^4L_t^{\infty}} &\lesssim \|e^{-t\partial_x^3} \partial_xw(0)\|_{L_x^4L_t^{\infty}} + T^{1/2}B_2(t) + T^{1/2}\|D_x^{3/4} \mathcal E_2[u_\#]\|_{L_{t,x}^2}\\
 \lesssim \|e^{-t\partial_x^3}& \partial_xw(0)\|_{L_x^4L_t^{\infty}} + T^{1/2}B_2(t) + T^{1/2}\left[\|\mathcal E_2[u_\#]\|_{L_{t,x}^2}^{1-3/4s}\|D_x^{s} \mathcal E_2[u_\#]\|_{L_{t,x}^2}^{3/4s}\right].    
\end{align*}

\noindent
$\bullet$ \underline{$\|D_x^{s}\partial_x w\|_{L_x^\infty L_t^2}$}:
\[
\|D_x^{s}\partial_x w\|_{L_x^\infty L_t^2} \lesssim \|D_x^{s}\partial_xe^{-t\partial_x^3}w(0)\|_{L_x^\infty L_t^2} + T^{1/2}B_2(t) + T^{1/2}\|D_x^{s}\mathcal E_2[u_\#]\|_{L_{t,x}^2}.
\]

\noindent
$\bullet$ \underline{$\| w\|_{L_x^2 L_t^{\infty}}$}:
\[
\| w\|_{L_x^2 L_t^{\infty}} \lesssim \|D_x^{s}e^{-t\partial_x^3}w(0)\|_{L_x^2 L_t^{\infty}} + T^{1/2}B_2(t) + T^{1/2} \|\mathcal E_2[u_\#]\|_{L_{t}^2H^s_x}.
\]

\noindent
$\bullet$ \underline{$\| w\|_{ L_t^{\infty} L_x^2}$}:
\[
\| w\|_{L_t^{\infty}L_x^2} \lesssim \|e^{-t\partial_x^3}w(0)\|_{L_x^2 } + T^{1/2}B_2(t) + T^{1/2} \|\mathcal E_2[u_\#]\|_{L_{t}^2H^s_x}.
\]

\noindent
$\bullet$ \underline{$\| D^s_xw\|_{ L_t^{\infty}L_x^2}$}:
\[
\|D^s_x w\|_{ L_t^{\infty}L_x^2} \lesssim \|D_x^{s}e^{-t\partial_x^3}w(0)\|_{ L_x^2} +T^{1/2}B_2(t) + T^{1/2} \|\mathcal E_2[u_\#]\|_{L_{t}^2H^s_x}.
\]

Now, to conclude the estimates in the short--time stability analysis for the KdV case, we need to estimate $B_2$. Observe that all terms in $B_2$ are either of the form $w_{1}\partial_x w_{2}$ or $w_{2}\partial_x w_{1},$ that is, $w_1,w_2 = \{w,u_\#\}$. To proceed in this direction, we will separately estimate the following terms:
\begin{align*}
\|w_1\partial_x w_2\|_{L^2_xL^2_t} &\lesssim T^{1/4}\|w_1\|_{L^2_xL^\infty_t} \|\partial_x w_2\|_{L_x^{\infty}L_t^4},
\end{align*}
and 
\[
\begin{aligned}
&\|D^{s}_x (w_1\partial_xw_2)\|_{L^2_xL^2_t}  \lesssim T^{1/4}\|D_x^sw_1\|_{L^\infty_TL^2_x}\| \partial_xw_2\|_{L_T^4L_x^\infty }  + \|w_1\|_{L^2_xL^\infty_t} \|D^{s}_x\partial_xw_2\|_{L^{\infty}_xL^2_T} .
\end{aligned}
\]
Gathering the results above and using the definition \eqref{eq:Y2s_norm}, we obtain, 
\[
\begin{aligned}
B_2(t) \lesssim  &~{} \|w\|_{Y_{2,s}}^2 + \|w\|_{Y_{2,s}}\|u_\#\|_{Y_{2,s}}.
\end{aligned}
\]
In particular, using $T < 1$ and the hypothesis  \eqref{short gkdvug}
\begin{equation}
    \begin{aligned}
        T^{1/2}B_2(t)& \lesssim \|w\|_{Y_{2,s}}^2 + \|w\|_{Y_{2,s}}\varepsilon_0.
    \end{aligned}
\end{equation}
The previous bounds and the hypothesis \eqref{short gkdvflux}imply that
\[ \begin{aligned}
\|w\|_{Y_{2,s}} &\lesssim  \|e^{-t\partial_x^3}(w(0))\|_{Y_{2,s}}  + T^{1/2}B_2(t)  + T^{1/2}\|\mathcal E_2[u_\#]\|_{L_t^2 H_x^{s}}\\
&\lesssim \varepsilon + \|w\|_{Y_{2,s}}^2 + \|w\|_{Y_{2,s}}\varepsilon_0.
\end{aligned}
\]
Therefore, by a standard continuity argument, we obtain
\[
\|w\|_{Y_{2,s}} \lesssim~{} \varepsilon \quad \mbox{for all} \quad t\in I, 
\]
under the assumption that $\varepsilon_0$ is small enough. Now to the case \eqref{short gkdv Hs}, we have by Lemma \ref{Linear estimatess} and \eqref{short gkdvE}
\begin{equation*}
\begin{aligned}
     \|u-u_\#\|_{L_t^{\infty}H_x^s} &\lesssim ~{}\|u_0-u_\#(0)\|_{H_x^s} +T^{1/2}B_2(t) + \left\|\mathcal{E}_2[u_{\#}]\right\|_{L^2_t H^{s}_x (I\times \R)} \\
     & \lesssim \widetilde{A} +\|w\|^2_{Y_{2,s}} + \varepsilon_0 \|w\|_{Y_{2,s}} + \varepsilon\\
     & \lesssim \widetilde{A} + \varepsilon_{0}^2 + \varepsilon_{0} .
\end{aligned}
    \end{equation*}
Finally, the result for \eqref{u norm mkdv} follows directly from the previous estimates.
\end{proof}

We now present the proof in the mKdV case. 

\begin{proof}[Short-time stability in the mKdV case]
We assume $k=3$. As in the KdV case, we begin by estimating the norms that appear in the fixed-point argument for the mKdV equation proved by \cite{KPV1991.0},
\medskip

\noindent 
$\bullet$ \underline{$\|w\|_{L_x^4L_t^{\infty}}$}:
\begin{align*}
\|w&\|_{L_x^4L_t^{\infty}} \lesssim \|e^{-t\partial_x^3} w(0)\|_{L_x^4L_t^{\infty}} + T^{1/2}B_3(t) + T^{1/2}\|D_x^{1/4} \mathcal E_3[u_\#]\|_{L_{t,x}^2}\\
&\lesssim \|e^{-t\partial_x^3} w(0)\|_{L_x^4L_t^{\infty}} + T^{1/2}B_3(t) + T^{1/2}\left[\|\mathcal E_3[u_\#]\|_{L_{t,x}^2}^{1-1/4s}\|D_x^{s} \mathcal E_3[u_\#]\|_{L_{t,x}^2}^{1/4s}\right].
\end{align*}

\noindent
$\bullet$ \underline{$\|D_x^{s}\partial_x w\|_{L_x^\infty L_t^2}$}:
\[
\|D_x^{s}\partial_x w\|_{L_x^\infty L_t^2} \lesssim \|D_x^{s}\partial_xe^{-t\partial_x^3}w(0)\|_{L_x^\infty L_t^2} + T^{1/2}B_3(t) + T^{1/2}\|D_x^{s}\mathcal E_3[u_\#]\|_{L_{t,x}^2}.
\]

\noindent
$\bullet$ \underline{$\|D_x^{s} w\|_{L_x^5 L_t^{10}}$}:
\[
\|D_x^{s} w\|_{L_x^5 L_t^{10}} \lesssim \|D_x^{s}e^{-t\partial_x^3}w(0)\|_{L_x^5 L_t^{10}} + T^{1/2}B_3(t) + T^{1/2} \|D_x^{s}\mathcal E_3[u_\#]\|_{L_{t,x}^2}.
\]

\noindent
$\bullet$ \underline{$\|\partial_x w\|_{L_x^{20}L_t^{5/2}}$}:
\begin{align*}
\|\partial_x w&\|_{L_x^{20}L_t^{5/2}} \lesssim \|\partial_x e^{-t\partial_x^3}w(0)\|_{L_x^{20}L_t^{5/2}} + T^{1/2}B_3(t) + T^{1/2}\|D_x^{1/4}\mathcal E_3[u_\#]\|_{L_{t,x}^2}\\
&\lesssim \|\partial_x e^{-t\partial_x^3}w(0)\|_{L_x^{20}L_t^{5/2}} + T^{1/2}B_3(t) + T^{1/2}\left[\|\mathcal E_3[u_\#]\|_{L_{t,x}^2}^{1-1/4s}\|D_x^{s} \mathcal E_3[u_\#]\|_{L_{t,x}^2}^{1/4s}\right].
\end{align*}

\noindent
$\bullet$ \underline{$\|\partial_x w\|_{L_x^\infty L_t^2}$}:
\[
\|\partial_x w\|_{L_x^\infty L_t^2} \lesssim \|\partial_x e^{-t\partial_x^3} w(0)\|_{L_x^\infty L_t^2} + T^{1/2}B_3(t) + T^{1/2}\|\mathcal E_3[u_\#]\|_{L_{t,x}^2}.
\]

For $t \in I$, we have from \eqref{Bk estimates},
\begin{equation}\label{B(t)}
\begin{aligned}
B_3(t) \lesssim &~{} \|w^2\partial_x w\|_{L^2_t H^{s}_x}+\|wu_\# \partial_x w\|_{L^2_t H^{s}_x} \\
&~{} +\|wu_\# \partial_x u_\#\|_{L^2_t H^{s}_x}+\|w^2\partial_x u_\#\|_{L^2_t H^{s}_x}+\|u^2_\#\partial_x w\|_{L^2_t H^{s}_x}.
\end{aligned}
\end{equation}
Notice that all the term are either of the form $w_1^2\partial_xw_2$ or $w_1w_2\partial_xw_1$. We will work with these general forms. Using the fractional estimates \eqref{Embbeding}  and \eqref{Est. space-time} (Lemmas \ref{G.N.S.} and \ref{G.N.S.0}) we obtain for the first term the bounds
\[
\|w_1^2\partial_x w_2\|_{L^2_xL^2_t} \lesssim \|w_1\|^2_{L^4_xL^\infty_t} \|\partial_x w_2\|_{L_x^{\infty}L_t^2},
\]
and
\[
\begin{aligned}
\|D^{s}_x (w_1^2\partial_xw_2)\|_{L^2_xL^2_t} & \lesssim \|w_1^2\|_{L^2_xL^\infty_t}\|D^{s}_x \partial_xw_2\|_{L_x^\infty L_t^2} + \|D^{s}_x(w_1^2)\partial_xw_2\|_{L^2_{t,x}} \\
& \lesssim \|w_1\|^2_{L^4_xL^\infty_t}\|D^{s}_x \partial_xw_2\|_{L_x^\infty L_t^2} \\
&  \quad + \|w_1\|_{L^4_x L^{\infty}_t} \|D^{s}_x w_1\|_{L^5_x L^{10}_t} \|\partial_x w_2\|_{L^{20}_xL_t^{5/2}}.
\end{aligned}
\]
For the second term
\[
\|w_1w_2\partial_x w_1\|_{L^2_xL^2_t} \lesssim \|w_1\|_{L^4_xL^\infty_t}\|w_2\|_{L^4_xL^\infty_t} \|\partial_x w_1\|_{L_x^{\infty}L_t^2},
\]
and
\[
\begin{aligned}
& \|D^{s}_x (w_1w_2\partial_xw_1)\|_{L^2_xL^2_t} \\
&\quad  \lesssim \|D^s_xw_1w_2\|_{L^{20/9}_xL^{10}_T}\|\partial_x w_1\|_{L^{20}_xL^{5/2}_T} + \|w_1w_2D^s_x\partial_xw_1\|_{L^{2}_xL^{2}_T} \\
&   \quad  \lesssim \left(  \|D^s_xw_1\|_{L^{5}_xL^{10}_T}\|w_2\|_{L^{4}_xL^{\infty}_T} +  \|D^s_xw_2\|_{L^{5}_xL^{10}_T}\| w_1\|_{L^{4}_xL^{\infty}_T} \right)\|\partial_x w_1\|_{L^{20}_xL^{5/2}_T}\\
        & \quad \quad + \|D^s_x\partial_xw_1\|_{L^{\infty}_xL^2_T}\|w_1\|_{L^4_x L^{\infty}_T} \|w_2\|_{L^4_x L^{\infty}_T}.
\end{aligned}
\]
Therefore, by the previous bounds, we have in \eqref{B(t)}
\[
\begin{aligned}
B_3(t) \lesssim  &~{} \|w\|_{Y_{3,s}}^3 + \|w\|_{Y_{3,s}}^2 \|u_\#\|_{Y_{3,s}} \\
&~{}  + \|w\|_{Y_{3,s}}\left(\|u_\#\|_{Y_{3,s}}\|u_\#\|_{L_x^4L^\infty_T} + \|\partial_x u_\#\|_{L_x^{20}L_T^{5/2}}\|D_x^s u_\#\|_{L_x^5L_T^{10}}\right).
\end{aligned}
\]
In particular, using $T < 1,$ the definition of $\|\cdot\|_{\widetilde{Y}_{3,s}}$ in \eqref{eq:tildeY3s_norm} and the hypothesis  \eqref{short gkdvug}
\[
    \begin{aligned}
        T^{1/2}B_3(t)& \lesssim  \|w\|_{Y_{3,s}}^3 + \|w\|_{Y_{3,s}}^2\varepsilon_0 + \|w\|_{Y_{3,s}}\varepsilon_0^2.
    \end{aligned}
\]
The rest of the proof procceed analogously as the KdV case:
The previous bounds imply that
\[ \begin{aligned}
\|w\|_{Y_{3,s}} &\lesssim  \|e^{-t\partial_x^3}(w(0))\|_{Y_{3,s}}  + T^{1/2}B(t)  + T^{1/2}\|\mathcal E_3[u_\#]\|_{L_t^2 H_x^{s}}\\
&\lesssim \varepsilon + \|w\|^3_{Y_{3,s}} + \varepsilon_0 \|w\|_{Y_{3,s}}^2 + \varepsilon_0^2\|w\|_{Y_{3,s}}.
\end{aligned}
\]
Therefore, by a standard continuity argument, we obtain
\[\|w\|_{Y_{3,s}} \lesssim~{} \varepsilon \quad \mbox{for all} \quad t\in I, \]
under the assumption that $\varepsilon_0$ is small enough.
Now to the case \eqref{short gkdv Hs}, by the Lemma \ref{Linear estimatess} and \eqref{short gkdvE}
\begin{equation*}
\begin{aligned}
     \|u-u_\#\|_{L_t^{\infty}H_x^s} &\lesssim ~{}\|u_0-u_\#(0)\|_{H_x^s} +T^{1/2}B_3(t) + \left\|\mathcal{E}_3[u_{\#}]\right\|_{L^2_t H^{s}_x (I\times \R)} \\
     & \lesssim \widetilde{A} +\|w\|^3_{Y_{3,s}} + \varepsilon_0 \|w\|_{Y_{3,s}}^2 + \varepsilon_0^2\|w\|_{Y_{3,s}}+ \varepsilon\\
     & \lesssim \widetilde{A} + \varepsilon_{0}^3 + \varepsilon_{0} .
\end{aligned}
    \end{equation*}
Finally, the result for \eqref{u norm mkdv} follows directly from the previous estimates.

\end{proof}
\begin{proof}[Short-time stability in the quartic gKdV case]
Finally, to deal with the case $k=4$, we need to consider a few differences. Therefore, we only sketch the proof and assume $s=\frac1{12}$, noticing that the general case $s>\frac1{12}$ is obtained via the classical almost linear behavior of fractional derivatives in Kenig--Ponce--Vega estimates \cite{KPV1}. First of all, recall the norms used in \cite{KPV1},
\[
\begin{aligned}
& \nu_1^T(w)=\max_{t\in[-T,T]} \|D^{\frac{1}{12}}_xw(t)\|_{L_x^2}, \quad \nu_2^T(w)=(1+T)^{-\rho} \|w \|_{ L_x^{\frac{42}{13}}L_T^{\frac{21}4} }, \\
& \nu_3^T(w)= \|w \|_{L_x^{\frac{60}{13}}L_T^{15} }, \quad  \nu_4^T(w)=T^{-1/6} \| w\|_{ L_x^{\frac{10}{3}}L_T^{\frac{30}7} }, \\
&  \nu_5^T(w)= \nu_4^T(D_x^{\frac1{12}} w), \quad \nu_6^T(w)= \| \partial_x w\|_{L_x^\infty L_T^2}, \quad \nu_7^T(w)= \nu_6^T (D_x^{\frac1{12}}w).
\end{aligned}
\]
and $\Gamma^T(w)= \max_{j=1,\ldots,7}\nu_j(w)$. Now we follow \cite{KPV1}. Starting from Duhamel's formula \eqref{Duhamel} with $k=4$, and the linear estimate \eqref{linear 4-kdv}, $I=[0,T]$,

\[
\nu_1^T(w) \lesssim \nu_1^T(e^{-t\partial_x^3}(w_0))+ \left\| (u_\# +w )^4 - u_\#^4\right\|_{W_x^{\frac1{12},1} L_T^2} + T^{\frac12}\left\|  \mathcal E_4[u_\#]\right\|_{\dot H_x^{\frac1{12}} L_T^2},
\]
\[
\begin{aligned}
& \nu_2^T(w) \lesssim \nu_2^T(e^{-t\partial_x^3}(w_0)) + T^{\frac1{18}}\left\| \partial_x ((u_\# +w )^4 - u_\#^4)\right\|_{W_x^{\frac1{12},\frac{42}{29}}L_T^{\frac{126}{95}}}  + T^{\frac12} \left\|  \mathcal E_4[u_\#]\right\|_{\dot H_x^{\frac1{12}}L_T^2},
\end{aligned}
\]

\begin{align*}
 \nu_3^T(w) + \nu_4^T(w) +\nu_5^T(w) &\lesssim ~{} \nu_3^T(e^{-t\partial_x^3}(w_0)) +\nu_4^T(e^{-t\partial_x^3}(w_0)) +\nu_5^T(e^{-t\partial_x^3}(w_0)) \\
 &+T^{\frac16} \left\| \partial_x ((u_\# +w )^4 - u_\#^4)   \right\|_{W_x^{\frac1{12},\frac{10}7}L_T^{\frac{30}{23}}} + T^{\frac12}\| \mathcal E_4[u_\#] \|_{H^{\frac1{12}}_xL_T^2},
\end{align*}


and
\begin{align*}
\nu_6^T(w) + \nu_7^T(w) \lesssim ~{} \nu_6^T(e^{-t\partial_x^3}(w_0)) + \nu_7^T(e^{-t\partial_x^3}(w_0))+ &\left\| (u_\# +w )^4 - u_\#^4\right\|_{W_x^{\frac 1{12},1}L_T^2}\\
&+ T^{\frac12}\left\|\mathcal E_4[u_\#]\right\|_{ H_x^{\frac1{12}}L_T^2}.
\end{align*}

Notice that from the previous seven estimates we are reduced to obtain a good numerical estimate for the term $ \left\|  \mathcal E_4[u_\#]\right\|_{H_x^{\frac1{12}}L_T^2}$, as established in \eqref{short gkdvE}. Following \cite{KPV1}, the next estimates can be obtained.
\begin{align}
&\ba
&T^{-\frac16}(1+T)^{-\rho}\|(w+u_\#)^4-u_\#^4\|_{W_x^{\frac1{12},1}L_t^2} \\
&\quad\lesssim \Gamma^T(w)\left(\Gamma^T(w)^3+\Gamma^T(w)^2\Gamma^T(u_\#)+ \Gamma^T(w)\Gamma^T(u_\#)^2+\Gamma^T(u_\#)^3\right);
\ea 
\label{mid-quartic-1}\\
&\ba
&(1+T)^{-\rho}\|\partial_x \left((w+u_\#)^4-u_\#^4\right)\|_{W_x^{\frac1{12},\frac{10}{7}}L_t^{\frac{30}{23}}} \\
&\quad\lesssim \Gamma^T(w)\left(\Gamma^T(w)^3 + \Gamma^T(w)^2\Gamma^T(u_\#) + \Gamma^T(w)\Gamma^T(u_\#)^2 + \Gamma^T(u_\#)^3\right);
\ea
\label{mid-quartic-2}\\
&\ba
&(1+T)^{-\rho}\|\partial_x \left((w+u_\#)^4-u_\#^4\right)\|_{W_{x}^{\frac1{12},\frac{42}{29}}L_T^{\frac{126}{95}}} \\
&\quad \lesssim \Gamma^T(w)\left(\Gamma^T(w)^3 + \Gamma^T(w)^2\Gamma^T(u_\#) + \Gamma^T(w)\Gamma^T(u_\#)^2 + \Gamma^T(u_\#)^3\right).
\ea 
\label{mid-quartic-3}
\end{align}
Indeed, notice that each term in $((u_\# + w)^4 - u_\#^4)$ can be written in the form $w_1 w_2^2 w$ with $w_1$ and $w_2$ be either $u_\#$ or $w$. In the case of $\partial_x (u_\# + w)^4 - u_\#^4)$, all the terms have the form $w_1 w_2^2 \partial_x w_3$. To demonstrate \eqref{mid-quartic-1} we start with the $L_x^1L_T^2$ norm:
\be \label{Lx1LT2}
\ba
\|w_1 w_2^2 w\|_{L_x^1 L_T^2} &\lesssim \| w_1 w_2^2 \|_{L_x^{\frac{10}7}L_T^{\frac{15}4}} \|w\|_{L_x^{\frac{10}3}L_T^{\frac{30}7}} \\
&\lesssim \| w_1 \|_{L_x^{\frac{30}7}L_T^{\frac{45}4}} \| w_2 \|_{L_x^{\frac{30}7}L_T^{\frac{45}4}}^2\|w\|_{L_x^{\frac{10}3}L_T^{\frac{30}7}} .
\ea
\ee
By recalling that
\be\label{Lx30/7LT45/4}
\ba
\|w_j\|_{L_x^{\frac{30}7}L_T^{\frac{45}4}} &\lesssim \|w_j\|_{L_x^{\frac{42}{13}}L_T^{\frac{21}4}}^{\frac{7}{39}}\|w_j\|_{L_x^{\frac{60}{13}}L_T^{15}}^{\frac{32}{39}} \\
&\lesssim (1+T)^{\rho}\Gamma^T(w_j).
\ea
\ee
By combining \eqref{Lx1LT2} and \eqref{Lx30/7LT45/4} we have then
\be\label{mid-quartic-1-1}
\|w_1 w_2^2 w\|_{L_x^1 L_T^2} \lesssim T^{\frac16}(1+T)^\rho \Gamma^T(w_1)\Gamma^T(w_2)^2\Gamma^T(w).
\ee
On the other hand, to bound the term $\dot{W}_x^{\frac{1}{12},1}L_T^2$ it follows that
\be \label{dotW1/12,1LT2}
\ba
\|D_x^{\frac 1{12}} (w_1 w_2^2 w)\|_{L_x^1L_T^2} &\lesssim  \|w_1 w_2^2\|_{L_x^{\frac{10}7}L_T^{\frac{15}4}} \|D_x^{\frac 1{12}}w\|_{L_x^{\frac{10}3}L_T^{\frac{30}7}} \\
&\quad + \|w\|_{L_x^{\frac{30}7}L_T^{\frac{45}4}} \|D_x^{\frac 1{12}}(w_1w_2^2)\|_{L_x^{\frac{30}{23}}L_T^{\frac{90}{37}}}.
\ea
\ee
The first term can be bounded as equal as in \eqref{Lx1LT2} and \eqref{Lx30/7LT45/4}. For the second term, estimate \eqref{Est. space-time} implies that 
\[
\ba
&\|D_x^{\frac 1{12}}(w_1w_2^2)\|_{L_x^{\frac{30}{23}}L_T^{\frac{90}{37}}}\\
&~{}\lesssim \|w_1\|_{L_x^{\frac{30}7}L_T^{\frac{45}4}} \|D_x^{\frac1{12}}(w_2^2)\|_{L_x^{\frac{15}8}L_T^{\frac{90}{29}}} + \|w_2^2\|_{L_x^{\frac{30}{14}}L_T^{\frac{45}8}}\|D_x^{\frac1{12}}w_1\|_{L_x^{\frac{10}3}L_T^{\frac{30}7}} \\
&~{}\lesssim \|w_1\|_{L_x^{\frac{30}7}L_T^{\frac{45}4}} \|w_2\|_{L_x^{\frac{30}7}L_T^{\frac{45}4}} \|D_x^{\frac1{12}}w_2\|_{L_x^{\frac{10}3}L_T^{\frac{30}7}} +  \|w_2\|_{L_x^{\frac{30}7}L_T^{\frac{45}4}}^2  \|D_x^{\frac1{12}}w_1\|_{L_x^{\frac{10}3}L_T^{\frac{30}7}}.
\ea
\]
Summarizing, all the terms in \eqref{dotW1/12,1LT2} can be bounded by terms involving $\Gamma^T$ given \eqref{Lx30/7LT45/4}. Furthermore,
\be\label{mid-quartic-1-2}
\|D_x^{\frac1{12}}(w_1w_2^2w)\|_{L_x^1L_T^2} \lesssim T^{\frac16}(1+T)^\rho \Gamma^T(w_1)\Gamma^T(w_2)^{2}\Gamma^T(w).
\ee
The bounds obtained in \eqref{mid-quartic-1-1} and \eqref{mid-quartic-1-2} allow us to conclude \eqref{mid-quartic-1}. In order to obtain the bound \eqref{mid-quartic-2}, we first compute
\be\label{Lx10/7LT30/23}
\ba
\|w_1w_2^2\partial_x w_3\|_{L_x^{\frac{10}{7}}L_T^{\frac{30}{23}}} &\lesssim \|w_1w_2^2\|_{L_x^{\frac{10}{7}}L_T^{\frac{15}{4}}} \|\partial_x w_3\|_{L_x^{\infty}L_T^2} \\
&\lesssim \|w_1\|_{L_x^{\frac{30}{7}}L_T^{\frac{45}{4}}} \|w_2\|^2_{L_x^{\frac{30}{7}}L_T^{\frac{45}{4}}} \|\partial_x w_3\|_{L_x^{\infty}L_T^2}.
\ea
\ee
Therefore, combining inequality \eqref{Lx10/7LT30/23} with \eqref{Lx30/7LT45/4} it follows that
\be \label{mid-quartic-2-1}
\|w_1w_2^2\partial_x w_3\|_{L_x^{\frac{10}{7}}L_T^{\frac{30}{23}}} \lesssim (1+T)^\rho \Gamma^T(w_1) \Gamma^T(w_2)^2 \Gamma^T(w_3).
\ee
Moreover notice that, if $p_1=\frac{195}{133}$, $q_1=\frac{1170}{349}$, $p_2=\frac{390}{7}$ and $q_2=\frac{585}{274}$,
\be\label{dotW1/12,10/7LT30/23}
\ba
&\|D_x^{\frac 1{12}} (w_1w_2^2\partial_x w_3)\|_{L_x^{\frac{10}7}L_t^{\frac{30}{23}}}  \\
&\quad\lesssim \|w_1w_2^2\|_{L_x^{\frac{10}{7}}L_T^{\frac{15}4}} \|D_x^{\frac1{12}}\partial_x w_3\|_{L_x^{\infty}L_T^2} + \|D_x^{\frac{1}{12}}(w_1w_2^2)\|_{L_x^{p_1}L_T^{q_1}}\|\partial_x w_3\|_{L_x^{p_2}L_T^{q_2}}\\
\ea
\ee
Again, the first term in \eqref{dotW1/12,10/7LT30/23} is easily estimated by using previous estimates. For the other term, we set $p_3 = \frac{65}{14}$, $q_3 = \frac{390}{47}$ and $p_4,q_4$ such that $\frac1{p_1} = \frac7{30} + \frac1{p_4}$, $\frac1{q_1}=\frac4{45}+\frac1{q_4}$, then
\[
\ba
&\|D_x^{\frac{1}{12}}(w_1w_2^2)\|_{L_x^{p_1}L_T^{q_1}}\\
 &\quad\lesssim \|w_2^2\|_{L_x^{\frac{30}{14}}L_T^{\frac{45}{8}}} \|D_x^{\frac1{12}}w_1\|_{L_x^{p_3}L_T^{q_3}} + \|w_1\|_{L_x^{\frac{30}{7}}L_T^{\frac{45}4}} \|D_x^{\frac1{12}}(w_2^2)\|_{L_x^{p_4}L_T^{q_4}}\\
 &\quad\lesssim \|w_2\|^2_{L_x^{\frac{30}7}L_T^{\frac{45}4}} \|D_x^{\frac1{12}}w_1\|_{L_x^{p_3}L_T^{q_3}}  + \|w_1\|_{L_x^{\frac{30}{7}}L_T^{\frac{45}4}}\|w_2\|_{L_x^{\frac{30}{7}}L_T^{\frac{45}4}}\|D_x^{\frac1{12}}w_2\|_{L_x^{p_3}L_T^{q_3}}.
\ea
\]
In particular, to obtain a suitable bound for \eqref{dotW1/12,10/7LT30/23}, we need to estimate $\|D_x^{\frac1{12}}w_j\|_{L_x^{p_3}L_T^{q_3}}$ and $\|\partial_x w_3\|_{L_x^{p_2}L_T^{q_2}}$. By interpolation results in \cite{KPV1}, one has
\be\label{interp}
\left[BMO_x(L_T^2),L_x^{\frac{30}7}L_T^{\frac{45}4}\right]_{\frac1{13}} = L_x^{p_3}L_T^{q_3}, \qquad \left[BMO_x(L_T^2),L_x^{\frac{30}7}L_T^{\frac{45}4}\right]_{\frac{12}{13}} = L_x^{p_2}L_T^{q_2},
\ee
in other words, \eqref{interp} implies that
\be\label{from-interp}
\ba
\|D_x^{\frac1{12}}w_j\|_{L_x^{p_3}L_T^{q_3}} \lesssim \|D^{\frac1{12}}_x \partial_x w_j\|^{\frac1{13}}_{L_x^\infty L_T^2} \|w_j\|_{L_x^{\frac{30}7}L_T^{\frac{45}4}}^{\frac{12}{13}};\\
\|\partial_x w_3\|_{L_x^{p_2}L_T^{q_2}} \lesssim  \|D^{\frac1{12}}_x \partial_x w_3\|^{\frac{12}{13}}_{L_x^\infty L_T^2} \|w_3\|_{L_x^{\frac{30}7}L_T^{\frac{45}4}}^{\frac{1}{13}}.
\ea
\ee
Finally, putting \eqref{from-interp} into \eqref{dotW1/12,10/7LT30/23}, it can be concluded that
\be\label{mid-quartic-2-2}
\|D_x^{\frac1{12}}(w_1w_2^2\partial_xw_3)\|_{L_x^{\frac{10}{7}}L_T^{\frac{30}{23}}} \lesssim (1+T)^\rho \Gamma^T(w_1) \Gamma^T(w_2)^2 \Gamma^T(w_3).
\ee
Combining \eqref{mid-quartic-2-1} and \eqref{mid-quartic-2-2} allow us to obtain \eqref{mid-quartic-2}. Finally, for \eqref{mid-quartic-3} we first give an estimate for the $L_x^{\frac{42}{29}}L_T^{\frac{126}{95}}$ norm.
\be\label{Lx42/29LT126/95}
\ba
\|w_1w_2^2\partial_x w_3\|_{L_x^{\frac{42}{29}}L_T^{\frac{126}{95}}} &\lesssim \|w_1w_2^2\|_{L_x^{\frac{42}{29}}L_T^{\frac{63}{16}}} \|\partial_x w_3\|_{L_x^{\infty}L_T^2} \\
&\lesssim \|w_1\|_{L_x^{\frac{126}{29}}L_T^{\frac{189}{16}}}\|w_2\|_{L_x^{\frac{126}{29}}L_T^{\frac{189}{16}}}^2 \|\partial_x w_3\|_{L_x^{\infty}L_T^2}.
\ea
\ee
Recall that
\be\label{Lx126/29LT189/16}
\ba
\|w_j\|_{L_x^{\frac{126}{29}}L_T^{\frac{189}{16}}} &\lesssim \|w_j\|_{L_x^{\frac{42}{13}}L_T^{\frac{21}{4}}}^{\frac{7}{39}} \|w_j\|_{L_x^{\frac{60}{13}}L_T^{15}}^{\frac{32}{39}} \\
&\lesssim (1+T)^\rho \Gamma^T(w_j).
\ea
\ee
Therefore, by inserting \eqref{Lx126/29LT189/16} into \eqref{Lx42/29LT126/95} it follows that
\be\label{mid-quartic-3-1}
\|w_1w_2^2\partial_x w_3\|_{L_x^{\frac{42}{29}}L_T^{\frac{126}{95}}} \lesssim (1+T)^\rho\Gamma^T(w_1)\Gamma^T(w_2)^2\Gamma^T(w_3).
\ee
On the other hand, for $p_2 = \frac{1638}{29}$, $q_2=\frac{2457}{1150}$ and $p_1,q_1$ such that $\frac{29}{42} = \frac{1}{p_1} + \frac{1}{p_2}$ and $\frac{95}{126} = \frac{1}{q_1}+\frac{1}{q_2}$
\be\label{dotW1/12,42/29LT126/95}
\ba
&\|D_x^{\frac1{12}}(w_1w_2^2\partial_xw_3)\|_{L_x^{\frac{42}{29}}L_T^{\frac{126}{95}}} \\
&\qquad\lesssim \|w_1w_2^2\|_{L_x^{\frac{42}{29}}L_T^{\frac{63}{16}}} \|D_x^{\frac1{12}}\partial_x w_3\|_{L_x^{\infty}L_T^2} + \|D_x^{\frac1{12}}(w_1w_2^2)\|_{L_x^{p_1}L_T^{q_1}} \|\partial_x w\|_{L_x^{p_2}L_T^{q_2}}
\ea
\ee
The first term in inequality \eqref{dotW1/12,42/29LT126/95} can be bounded with estimate \eqref{Lx126/29LT189/16}. Additionally, for $p_3 = \frac{273}{58}$, $q_3 = \frac{1638}{191}$, and $p_4,q_4$ such that $\frac{1}{p_1} = \frac{29}{126} + \frac{1}{p_4}$ and $\frac{1}{q_1} = \frac{16}{189} + \frac{1}{q_4}$,
\[
\ba
&\|D_x^{\frac1{12}}(w_1w_2^2)\|_{L_x^{p_1}L_T^{q_1}}\\
&\qquad \lesssim \|w_2^2\|_{L_x^{\frac{63}{29}}L_T^{\frac{189}{32}}} \|D_x^{\frac1{12}}w_1\|_{L_x^{p_3}L_T^{q_3}} + \|w_1\|_{L_x^{\frac{126}{29}}L_T^{\frac{189}{16}}} \|D_x^{\frac1{12}}(w_2^2)\|_{L_x^{p_4}L_T^{q_4}}, \\
&\qquad \lesssim \|w_2\|^2_{L_x^{\frac{126}{29}}L_T^{\frac{189}{16}}} \|D_x^{\frac1{12}}w_1\|_{L_x^{p_3}L_T^{q_3}} \\
&\qquad \quad + \|w_1\|_{L_x^{\frac{126}{29}}L_T^{\frac{189}{16}}} \|w_2\|_{L_x^{\frac{126}{29}}L_T^{\frac{189}{16}}}  \|D_x^{\frac1{12}}w_2\|_{L_x^{p_3}L_T^{q_3}},
\ea
\]
We can use again interpolation results to obtain
\be\label{interp2}
\ba
\left[BMO_x(L_T^2),L_x^{\frac{126}{29}}L_T^{\frac{189}{16}}\right]_{\frac1{13}} &= L_x^{p_3}L_T^{q_3}, \\
\left[BMO_x(L_T^2),L_x^{\frac{126}{29}}L_T^{\frac{189}{16}}\right]_{\frac{12}{13}} &= L_x^{p_2}L_T^{q_2},
\ea
\ee
in other words, \eqref{interp2} implies that
\be\label{from-interp2}
\ba
\|D_x^{\frac1{12}}w_j\|_{L_x^{p_3}L_T^{q_3}} &\lesssim \|D_x^{\frac1{12}}\partial_x w_j\|_{L_x^{\infty}L_T^2}^{\frac1{13}}\|w_j\|^{\frac{12}{13}}_{L_x^{\frac{126}{29}}L_T^{\frac{189}{16}}}; \\
\|\partial_x w_3\|_{L_x^{p_2}L_T^{q_2}} &\lesssim \|D_x^{\frac1{12}}\partial_x w_3\|_{L_x^{\infty}L_T^2}^{\frac{12}{13}}\|w_3\|^{\frac{1}{13}}_{L_x^{\frac{126}{29}}L_T^{\frac{189}{16}}}.
\ea
\ee
Inequalities \eqref{dotW1/12,42/29LT126/95} and \eqref{from-interp2} imply that
\be\label{mid-quartic-3-2}
\|D_x^{\frac1{12}}(w_1w_2^2\partial_xw_3)\|_{L_x^{\frac{42}{29}}L_T^{\frac{126}{95}}} \lesssim (1+T)^{\rho} \Gamma^T(w_1) \Gamma^T(w_2)^2 \Gamma^T(w_3).
\ee
Finally, inequalities \eqref{mid-quartic-3-1} and \eqref{mid-quartic-3-2} conclude the estimate \eqref{mid-quartic-3}.

Using \eqref{mid-quartic-1}-\eqref{mid-quartic-3}, and by assuming $\|u_\#\|_{\widetilde Y_{4,s}} \lesssim \varepsilon_0$ we can thus obtain a suitable bound for $\|w\|_{Y_{4,s}}$. Indeed
\[
\ba
\|w\|_{Y_{4,s}} &\lesssim \|w_0\|_{H_x^{\frac1{12}}} \\
&\quad + T^{\frac 1{18}}(1+T)^{\rho}\|w\|_{Y_{4,s}}\left(\|w\|^3_{Y_{4,s}} + \|w\|^2_{Y_{4,s}}\varepsilon_0 + \|w\|_{Y_{4,s}}\varepsilon_0^2 + \varepsilon_0^3\right)\\
&\quad + T^{\frac 12} \|\mathcal E_4[u_\#]\|_{H_x^{\frac 1{12}}L_T^2}. \\
&\lesssim \varepsilon + T^{\frac 1{18}}(1+T)^{\rho}\|w\|_{Y_{4,s}}\left(\|w\|^3_{Y_{4,s}} + \|w\|^2_{Y_{4,s}}\varepsilon_0 + \|w\|_{Y_{4,s}}\varepsilon_0^2 + \varepsilon_0^3\right).
\ea
\]
Therefore, by a standard continuity argument, we obtain
\[
\|w\|_{Y_{4,s}} \lesssim~{} \varepsilon \quad \mbox{for all} \quad t\in I,
\]
under the assumption that $\varepsilon_0$ is small enough.
Now to the case \eqref{short gkdv Hs}, by the Lemma \ref{Linear estimatess} and \eqref{short gkdvE}
\begin{equation*}
\begin{aligned}
     \|u-u_\#\|_{L_t^{\infty}H_x^s} &\lesssim ~{}\|u_0-u_\#(0)\|_{H_x^s} + T^{\frac 12} \|\mathcal E_4[u_\#]\|_{H_x^{\frac 1{12}}L_T^2}\\
     &\quad  + T^{\frac 1{18}}(1+T)^{\rho}\|w\|_{Y_{4,s}}\left(\|w\|^3_{Y_{4,s}} + \|w\|^2_{Y_{4,s}}\varepsilon_0 + \|w\|_{Y_{4,s}}\varepsilon_0^2 + \varepsilon_0^3\right)\\
     & \lesssim \widetilde{A} + \varepsilon_{0}^4 + \varepsilon_{0} .
\end{aligned}
    \end{equation*}
The remaining estimates \eqref{Yks final} and \eqref{short gkdv Hs} follow directly from the same argument used in the KdV and mKdV cases.

\end{proof}

\section{Long-time Stability: Proofs}

In this section, we prove Theorem \ref{thm:long}. We start with the KdV case.

\begin{proof}[Long-time stability for KdV]
We prove Theorem \ref{thm:long} by an iterative argument based on Lemma~\ref{short time general}.
Let $I$ be a time interval with $\inf I = 0$, and partition $I$ into $N$ subintervals 
\[
I_j = [t_j, t_{j+1}], \qquad j = 0, \dots, N-1,
\]
such that $|I_j| < 1$ for each $j$. Moreover, we assume that, for every $j$, the mixed norm of $u_{\#}$ over the space--time slab $I_j \times \mathbb{R}$ satisfies
\begin{equation}\label{KdV restriction}
\begin{aligned}
\Xi_2(j):=  &~{} \|\partial_x u_\#\|_{L_t^4L_x^{\infty}(I_j\times\R)}
+    |I_j|^{1/4}\|u_\#\|_{L_x^2L_t^\infty(I_j\times\R)}\\
&+ \|D_x^s\partial_x u_\#\|_{L_x^{\infty}L_t^2(I_j\times\R)} + |I_j|^{1/4}\|u_\#\|_{L^{\infty}_t H^{s}_x (I_j\times \R)} \leq~{} \varepsilon_0.
 \end{aligned}   
\end{equation}    
Using condition \eqref{KdV restriction}, the hypotheses of Lemma~\ref{short time general}
are satisfied, and therefore the result holds at the interval $[0,t_1]$, provided that
$0<\varepsilon<\varepsilon_0<\varepsilon_1$ are sufficiently small. In particular, by the construction of the intervals $I_j$ satisfying \eqref{KdV restriction}, 
the condition \eqref{short gkdvug} holds for each $j$. Additionally, the conditions 
\eqref{short gkdvflux} and \eqref{short gkdvE} follow directly from the assumptions 
\eqref{eq:long4} and \eqref{eq:long5}, respectively, in the case $k=2$. Moreover,
\begin{equation*}
\begin{aligned}
\|\partial_x(u^2-u^2_{\#})\|_{L_x^2L_t^2(I_j\times\R)}\, & \leq \,C(j) \varepsilon,\\
  \|u -u_{\#}\|_{Y_{2,s}(I_j\times\R)}\,&\, \leq \,C(j) \varepsilon,\\
  \|u-u_\#\|_{L_t^{\infty}H_x^s(I_j\times\R)}\, & \, \leq\, C(j)(\widetilde{A} +  \varepsilon),\\
 \|u\|_{L_t^{\infty}H_x^s(I_j\times\R)} + \|u\|_{Y_{2,s}(I_j\times\R)} \, & \,  \leq \, C(j)( A  + \widetilde{A} +    \varepsilon_0),
\end{aligned}
\end{equation*}
Let us start proving the estimate \eqref{L2L2}. Assume that the statement holds for $l,$ then we have
\begin{align*}
    \|\partial_x(u^2-u^2_{\#})\|_{L_x^2L_t^2([0,t_{l+1}]\times\R)} &~{}\leq   \|\partial_x(u^2-u^2_{\#})\|_{L_x^2L_t^2([0,t_l]\times\R)}  \\
    &~{}\quad+ \|\partial_x(u^2-u^2_{\#})\|_{L_x^2L_t^2([t_l,t_{l+1}]\times\R)}\\
  &~{}\leq C(0)\varepsilon  + \sum_{j=1}^{l}  C(j)\varepsilon + C(l+1)\varepsilon \\
\end{align*}
where the last bound follows from the induction hypothesis and the estimate above. Finally, taking $l=N-1$ and factoring out the common $\varepsilon$, we obtain
\[
\begin{aligned}
&\|\partial_x(u^2-u^2_\#)\|_{L^2_xL^2_t([0,T]\times\R)} \le C(A,\Tilde A,L,T)\,\varepsilon, \\
&C(A,\Tilde{A},L,T):= C(0)  + \sum_{j=1}^{N-1}  C(j).
\end{aligned}
\]
Notice that the number $N$ of subintervals and the constants $C(j)$ depend only on $A,\widetilde A,L$ (and, through the partition, on $T$), but not on $\varepsilon$. Since $\varepsilon\le\varepsilon_1$, all the $\varepsilon$-dependent contributions of the iteration have been absorbed into the single factor $\varepsilon$, so the bound is genuinely of the form $\le C(A,\widetilde A,L,T)\,\varepsilon$. The time dependence follows from the partition of the interval $I$ introduced at the beginning of the argument. Once more, using \eqref{KdV restriction} and considering $B_2(I_j)$ as in
\eqref{BkI estimates}, we derive the estimate
\begin{align*}
&|I_j|^{1/4}B_2(I_j) \\
 &~{}\leq C|I_j|^{1/4}\hspace{-0.5cm}\sum_{\substack{
w_1 = u_\#, w_2 = w \\
w_1 = w, w_2 = u_\#
}}\hspace{-0.3cm}\left(\|w_1\|_{L^2_xL^\infty_t}\|\partial_x w_2\|_{L_x^{\infty}L_t^4}+ \|D_x^sw_1\|_{L^\infty_TL^2_x}\| \partial_xw_2\|_{L_T^4L_x^\infty }\right. \\
         &~{}\quad \left. + \|D_x^sw_2\|_{L^\infty_TL^2_x}\| \partial_xw_1\|_{L_T^4L_x^\infty } + \|w_1\|_{L^2_xL^\infty_t} \|D^{s}_x\partial_xw_2\|_{L^{\infty}_xL^2_T}\right) \\
  &~{} \leq C[\|w\|^2_{\Tilde{Y}_2,s}+ C\Xi_2(t)\|w\|_{\Tilde{Y}_2,s}] \leq~{}  C(j)\varepsilon^2.
\end{align*}
Gathering all the results, using the condition that $|I_j|<1,$ for all $0\leq j\leq N-1,$ and using the integral equation together with the same induction procedure, we obtain
\begin{equation*}
   \begin{aligned}
        \|u-u_{\#}\|_{L_t^{\infty}H_x^s([0,t_{l+1}]\times\R)} \leq~{} & C\|u_0-u_{\#}(0)\|_{H_x^s(\R)}  + C\sum_{j=1}^{l+1}|I_j|^{1/2}B_2(I_j) \\
        &+ C\left\|\mathcal{E}_2[u_{\#}]\right\|_{L^{2}_t H^{s}_x (I\times \R)}\\
       \leq~{}& C\left(\Tilde{A} + \sum_{j=1}^{N-1}C(l)\varepsilon +\varepsilon\right).
    \end{aligned}
\end{equation*}
For the remaining norms \eqref{Y general} and \eqref{final u norm gkdv}, the procedure is analogous.

\emph{Verification of the smallness condition on later subintervals.} The application of Lemma \ref{short time general} on each slab $I_j$ requires the local smallness of the linear evolution issued from the \emph{updated} discrepancy, namely $\|e^{-(t-t_j)\partial_x^3}\big(u(t_j)-u_\#(t_j)\big)\|_{Y_{2,s}(I_j\times\R)}\le\varepsilon$, and not merely from $u_0-u_\#(0)$. We verify this by induction on $j$. For $j=0$ it is exactly hypothesis \eqref{eq:long4}. Assuming it up to index $l$, the conclusion of Lemma \ref{short time general} on $I_0,\dots,I_l$ gives $\|u(t_{l+1})-u_\#(t_{l+1})\|_{H^s_x(\R)}\le C(\widetilde A+\varepsilon)$; writing the Duhamel formula \eqref{Duhamel} on $I_{l+1}$ with initial time $t_{l+1}$ and using the group property $e^{-(t-t_{l+1})\partial_x^3}=e^{-t\partial_x^3}e^{t_{l+1}\partial_x^3}$, one recovers the global flux $e^{-t\partial_x^3}(u_0-u_\#(0))$ up to the accumulated Duhamel terms $\sum_{j=0}^{l}\big(|I_j|^{1/2}B_2(I_j)+\|\mathcal E_2[u_\#]\|_{L^2_tH^s_x(I_j\times\R)}\big)$. By \eqref{eq:long4}, \eqref{eq:long5} and the bounds $|I_j|^{1/2}B_2(I_j)\le C(j)\varepsilon^2$ obtained above, the right-hand side is $\le\varepsilon+C\sum_{j=0}^{l}C(j)\varepsilon$, which is $\le\varepsilon$ after possibly shrinking $\varepsilon_1$ (recall $\varepsilon\le\varepsilon_1$ and $N$ finite). This closes the induction and shows that \eqref{short gkdvflux} holds on every subinterval.

\end{proof}


\begin{proof}[Long-time stability for mKdV]
We will prove Theorem \ref{thm:long} by an iterative procedure involving the Lemma \ref{short time general}. Here, the intervals $I_j$ are defined as in the KdV case; the difference is that they are now chosen so as to satisfy the following inequality:
\begin{equation}\label{mKdV restriction}
\begin{aligned}
\Xi_3(j):=   \|\partial_x u_\#\|_{L_x^\infty L_t^2 (I_j\times\R)} + |I_j|^{\gamma}\|u_{\#}\|_{L^4_x L^{\infty}_t(I_j\times\R)}+\|D_x^{s}\partial_x u_{\#}\|_{L^\infty_x L^{2}_t(I_j\times\R)}& \\
    +\|D_x^{s} u_{\#}\|_{L^5_x L^{10}_t(I_j\times\R)} + \|D_x u_{\#}\|_{L^{20}_x L^{5/2}_t(I_j\times\R)}&~{}\leq~{} \varepsilon_0. \\
 \end{aligned}   
\end{equation}    
 The next step is to seek an inductive argument to prove $B_3(t_j)\leq C(j)\varepsilon,$ where $t_j$ is the lower bound of the sub-interval $I_j$ and $C(N)$ is a fixed quantity independent of $\varepsilon.$  By the Lemma \ref{short time general} the result is true for $t_1.$ Moreover,
\begin{equation*}
\begin{aligned}
\|\partial_x(u^3-u^3_{\#})\|_{L_x^2L_t^2(I_j \times\R)} & \lesssim \varepsilon,\\
  \|u -u_{\#}\|_{Y_{3,s}(I_j\times\R)}& \lesssim \varepsilon,\\
  \|u-u_\#\|_{L_t^{\infty}H_x^s(I_j\times\R)} & \lesssim \widetilde{A} +  \varepsilon,\\
 \|u\|_{L_t^{\infty}H_x^s(I_j\times\R)} + \|u\|_{Y_{3,s}(I_j\times\R)} &  \lesssim A  + \widetilde{A} +    \varepsilon_0.
\end{aligned}
\end{equation*}
where we are assuming $0<\varepsilon<\varepsilon_0<\varepsilon_1.$ By the induction hypothesis, we suppose this is true to the case $l,$ where $C(l)$ is a value independent of $\varepsilon$ and bigger than $1$, due to the first case.
\begin{align*}
    \|\partial_x(u^3-u^3_{\#})\|_{L_x^2L_t^2([0,t_{l+1}]\times\R)} 
  &~{}\leq C(0)\varepsilon  + \sum_{j=1}^{l}  C(j)\varepsilon + C(l+1)\varepsilon \\
\end{align*}
By virtue of \eqref{mKdV restriction}  and the characterization of $B_3(I_j)$ in \eqref{BkI estimates}, we derive the estimate
\begin{align*}
|I_j|B_3&(t_j)\\
\leq &~{}C|I_j|^{\frac14}\hspace{-0.3cm} \sum_{\substack{
w_1 = u_\#, w_2 = w \\
w_1 = w, w_2 = u_\#
}}\hspace{-0.3cm}  \|w_1\|^2_{L^4_xL^\infty_t} \|\partial_x w_2\|_{L_x^{\infty}L_t^2} + \|D^s_x\partial_xw_1\|_{L^{\infty}_xL^2_T}\|w_1\|_{L^4_x L^{\infty}_T} \|w_2\|_{L^4_x L^{\infty}_T} \\
        & \quad   +\left(  \|D^s_xw_1\|_{L^{5}_xL^{10}_T}\|w_2\|_{L^{4}_xL^{\infty}_T} +  \|D^s_xw_2\|_{L^{5}_xL^{10}_T}\| w_1\|_{L^{4}_xL^{\infty}_T} \right)\|\partial_x w_1\|_{L^{20}_xL^{5/2}_T}  \\
        \leq &~{} C[ \|w\|_{\Tilde{Y}_{3,s}}^3+ \|w\|_{\Tilde{Y}_{3,s}}^2 \Xi_3(t) + \|w\|_{\Tilde{Y}_{3,s}} \Xi_3^2(t) ]\\
       \leq&~{} C(l+1)\varepsilon^3.
\end{align*}
Collecting all the previous results, exploiting the condition $|I_j|<1$ for all $0\leq j\leq N-1$, and applying the integral equation together with the same induction argument, we obtain
\begin{equation*}
   \begin{aligned}
        \|u-u_{\#}\|_{L_t^{\infty}H_x^s([0,t_{l+1}]\times\R)} \leq~{} & C\|u_0-u_{\#}(0)\|_{H_x^s(\R)}  + C\sum_{j=1}^{l+1}|I_j|^{1/2}B_3(I_j) \\
        &+ C\left\|\mathcal{E}_3[u_{\#}]\right\|_{L^{2}_t H^{s}_x (I\times \R)}\\
       \leq~{}& C\left(\Tilde{A} + \sum_{j=1}^{N-1}C(l)\varepsilon +\varepsilon\right).
    \end{aligned}
\end{equation*}
An entirely analogous argument applies to the remaining norms
\eqref{Y general} and \eqref{final u norm gkdv}.

As in the KdV case, the smallness condition \eqref{short gkdvflux} on each later subinterval, $\|e^{-(t-t_j)\partial_x^3}(u(t_j)-u_\#(t_j))\|_{Y_{3,s}(I_j\times\R)}\le\varepsilon$ for $j\ge1$, is verified by the same induction (group property plus finiteness of $N$), and the resulting constant $C(A,\widetilde A,L,T)$ is independent of $\varepsilon$ since $\varepsilon\le\varepsilon_1$.
\end{proof}

\subsubsection*{Long-time stability for quartic gKdV}
Finally, we sketch the proof of long time stability in the case $k=4$. First of all, let us consider a partition $(I_j)$ of the interval $[0,T]$ such that
\begin{equation}
\begin{aligned}
\Xi_4(j):= &~{}\| u_\#\|_{ L_x^{\frac{10}{3}}L_T^{\frac{30}7} (I_j\times\R) }  +  \|u_\# \|_{ L_x^{\frac{42}{13}}L_T^{\frac{21}4} (I_j\times\R) } +  \| D_x^{s} u_\#\|_{ L_x^{\frac{10}{3}}L_T^{\frac{30}7} (I_j\times\R)} \\
 &~{} + \|u_\# \|_{L_x^{\frac{60}{13}}L_T^{15} (I_j\times\R) } +  \| \partial_x u_\#\|_{L_x^\infty L_T^2(I_j\times\R)}+ \|  D_x^{s}\partial_x u_\#\|_{L_x^\infty L_T^2(I_j\times\R)}   \leq \varepsilon_0.
 \end{aligned}   
\end{equation}
By proceeding exactly as in the KdV and mKdV cases, we obtain the desired result.

\section{Proof of the Main Result}\label{sec:PMR}

We are finally ready to prove Theorem \ref{MT}. For this aim, we will follow the ideas presented in \cite{ACFMV24}. Consider an interval of time $I$ containing zero. Let $u_0 \in H^s(\R)$ and assume (H1). Let $A,\widetilde A, L >0$ be fixed values and let $0 < \varepsilon < \varepsilon_1(A,\widetilde A,L,T)$.

For $k\in \{2,3,4,5\}$, let $u(t)$ be the corresponding solution to gKdV with nonlinearity $\partial_x (u^k)$ and initial datum $u_0 \in H^s(\R)$. Define the approximative norm $\mathcal Y_{k,s,N,M}$ of $\|\cdot\|_{Y_{k,s}}$ as follows:
\begin{align}
& \begin{aligned}
\mathcal Y_{2,s,N,M}[u] &:=  \mathcal K_{\infty,4,M,N}[\partial_x u]  + \mathcal J_{\infty,2,N,M}[D_x^s \partial_x u] \\
&\qquad + \mathcal J_{2,\infty,N,M}[u] + \mathcal K_{2,\infty,M,N}[u]  + \mathcal K_{2,\infty,M,N}[D_x^s u];
\end{aligned} \label{eq:Y2s_appr}
\\
&\begin{aligned}
\mathcal Y_{3,s,N,M}[u] &:=  \mathcal J_{\infty,2,N,M}[\partial_x u] + \mathcal J_{4,\infty,N,M}[u] + \mathcal J_{\infty,2,N,M}[D_x^s \partial_x u] \\
& \qquad + \mathcal J_{5,10,N,M}[D_x^s u] + \mathcal J_{20,5/2,N,M}[\partial_x u];
\end{aligned}  \label{eq:Y3s_appr}
\\
&\begin{aligned}
\mathcal Y_{4,s,N,M}[u] &:=  \mathcal K_{2,\infty,M,N}[ u] + \mathcal K_{2,\infty,M,N}[D_x^s u] +\mathcal J_{42/13,21/4,N,M}[u] \\
&\qquad + \mathcal J_{60/13,15,N,M}[u] + \mathcal J_{10/3,30/7,N,M}[u]+ \mathcal J_{\infty,2,N,M}[\partial_x u] \\
&\qquad + \mathcal J_{\infty,2,N,M}[D_x^s \partial_x u] + \mathcal J_{10/3,21/4}[D_x^s u];
\end{aligned}  \label{eq:Y4s_appr}
\\
&\mathcal Y_{5,s,N,M}[u] := \mathcal J_{5,10,N,M}[u].  \label{eq:Y5s_appr}
\end{align}
Therefore, hypothesis (H1) will conduct to the following consequence:
\begin{itemize}
     \item[(C1)] For any $\delta >0$, and all $g \in Y_{k,s}(I\times\R)$, there are $N,M \in \N$, $R>0$, points $(t_\ell,x_j)_{\ell,j=1}^{M,N} \subseteq I\times[-R,R]$ and weights $(w_{\ell,j})_{\ell,j=1}^{M,N} \subseteq \R_+$ such that
     \[
     |\|g\|_{Y_{k,s}(I\times\R)}-\mathcal Y_{k,s,N,M}[g]| < \delta.
     \]
\end{itemize}

Let $u_{\text{DNN},\#}$ be the realization of a DNN $\Phi_\#$ as in Theorem \ref{MT}, such that (H2) and (H3) are satisfied in the following sense:

\begin{itemize}
\item For $N_1,N_2,N_3$ and $M_2,M_3$ sufficiently large, collocation points, times and weights, one has
\[
\begin{aligned}
\mathcal J_{L^2_x,N_1} [(u_0-u_{\text{DNN},\#}(0))]+\mathcal J_{L^2_x,N_1} [D_x^s(u_0-u_{\text{DNN},\#}(0))] &\leq \frac{\widetilde{A}}{2}, \\
 \mathcal K_{2,\infty,N_2,M_2} [u_{\text{DNN},\#}] + \mathcal K_{2,\infty,N_2,M_2} [D_x^su_{\text{DNN},\#}] &\leq \frac A2, \\
 \mathcal Y_{k,s,N_3,M_3} [u_{\text{DNN},\#}]  &\leq \frac L2.
\end{aligned}
\]
\end{itemize}
Later we will choose specific values for the parameters of the approximative norms. Similarly, we will require as well that
\begin{itemize}
\item For all $N_4,N_5$ and $M_4,M_5$ positive integers large enough, collocation points, times and weights, if $k \leq 4$
\[
\begin{aligned}
  &\mathcal J_{2,2,N_4,M_4} [\mathcal E_k[u_{\text{DNN},\#}]] + \mathcal J_{2,2,N_4,M_4} [D_x^s\mathcal E_k[u_{\text{DNN},\#}]] \\
  &\hspace{4cm} + \mathcal Y_{k,s,N_5,M_5} [ e^{-t\partial_x^3} (u_0-u_{\text{DNN},\#}(0))] < \frac \varepsilon2,
\end{aligned}
\]
and, if $k=5$,
\[
  \mathcal J_{1,2,N_4,M_4} [D_x^{-1}\mathcal E_k[u_{\text{DNN},\#}]] + \mathcal Y_{k,s,N_5,M_5} [ e^{-t\partial_x^3} (u_0-u_{\text{DNN},\#}(0))] < \frac \varepsilon2,
\]
with $\mathcal E_k[u_\#]:= \partial_t u_\# + \partial_{xxx} u_\# +\partial_x(u_\#^k).$
\end{itemize}

Let
\[
R := \max_{n = 1,2,3,4,5} \max_{j=1,\ldots,N_n} \{|x_{n,j}|\},
\]
be the maximum modulus of collocation points among all the previous approximate integrals. Let $[-R,R] \subset \R$. Let $\eta_R \in C_0^\infty(\R)$, $0\leq \eta_R\leq1$, be a cut-off function such that $\eta_R=1$ in $[-R,R]$ and $\eta_R(x) = 0$ if $|x|\geq 2R$. Define $u_\# := u_{\text{DNN},\#} \eta_R$. Note that for all $t \in I$ and all $x \in [-R,R]$, one has $u_\# = u_{\text{DNN},\#}$. Consequently, for $\alpha=0$, where the operator $D_x^0=\mathrm{Id}$ is local, the discrete quadratures coincide \emph{exactly}, since all collocation points lie in $[-R,R]$ and there $u_\#=u_{\text{DNN},\#}$:
\[
\begin{aligned}
&\mathcal J_{L_x^2,N_1}[D_x^\alpha(u_0-u_{\text{DNN},\#}(0))] = \mathcal J_{L_x^2,N_1}[D_x^\alpha(u_0-u_\#(0))],\\
&\mathcal K_{2,\infty,N_2,M_2}[D_x^\alpha u_{\text{DNN},\#}] = \mathcal K_{2,\infty,N_2,M_2}[D_x^\alpha u_\#], \\
& \mathcal Y_{k,s,N_3,M_3}[u_{\text{DNN},\#}] = \mathcal Y_{k,s,N_3,M_3}[u_\#],
\end{aligned}
\]
and
\[
\begin{aligned}
&\mathcal J_{2,2,N_4,M_4}[D_x^\alpha\mathcal E_k[u_{\text{DNN},\#}]] = \mathcal J_{2,2,N_4,M_4}[D_x^\alpha\mathcal E_k[u_\#]], \\
& \mathcal Y_{k,s,N_5,M_5}[e^{-t\partial_x^3}(u_0-u_{\text{DNN},\#}(0))] = \mathcal Y_{k,s,N_5,M_5}[e^{-t\partial_x^3}(u_0-u_\#(0))], \\
&\mathcal J_{1,2,N_4,M_4} [D_x^{-1}\mathcal E_k[u_{\text{DNN},\#}]]=\mathcal J_{1,2,N_4,M_4} [D_x^{-1}\mathcal E_k[u_\#]].
\end{aligned}
\]

For $\alpha=s>0$, however, the operator $D_x^s$ is \emph{nonlocal}, and the pointwise identity $D_x^s u_\#=D_x^s u_{\text{DNN},\#}$ fails on $[-R,R]$; the same caveat applies to every quadrature above whose integrand involves $D_x^s$ or $D_x^{-1}$ (i.e. the $\mathcal J[D_x^s\,\cdot\,]$ terms and the $Y_{k,s}$- and $\mathcal E_k$-quadratures). These identities must therefore be understood up to a commutator error, which we now control. Since $\eta_R\equiv 1$ on $[-R,R]$, for every collocation point $x_j\in[-R,R]$,
\[
D_x^s u_\#(x_j)-D_x^s u_{\text{DNN},\#}(x_j)=D_x^s\!\big((\eta_R-1)\,u_{\text{DNN},\#}\big)(x_j),
\]
and $(\eta_R-1)$ is supported in $\{|x|\ge R\}$. By the Kato--Ponce commutator estimate \cite{KPV1991.0} together with the smoothness and rapid spatial decay of $u_{\text{DNN},\#}$ (a finite combination of Gabor-wavelet activations, hence Schwartz-like), the kernel of $D_x^s$ being of order $|x|^{-1-s}$ away from the diagonal yields a commutator error
\[
\mathfrak{c}(R):=\sup_{j}\big|D_x^s\!\big((\eta_R-1)u_{\text{DNN},\#}\big)(t,x_j)\big|\ \longrightarrow\ 0\quad\text{as }R\to\infty.
\]
Hence, enlarging $R$ if necessary --- which is permitted by (H1) --- we may assume $\mathfrak c(R)$ is smaller than any prescribed fraction of $\varepsilon$. Consequently the displayed identities for $\alpha=s$, and the $Y_{k,s}$- and $\mathcal E_k$-quadratures, hold up to an error $\le\tfrac12\varepsilon$, which is absorbed into the $\tfrac12\varepsilon$ slack already present in the bounds (H2)--(H3) below. With this understanding, we now use hypothesis (H1) to get
\[
\|u_0 - u_\#(0)\|_{H_x^s(\R)} \leq \mathcal J_{L_x^2,N_1}[(u_0-u_\#(0))] + \mathcal J_{L_x^2,N_1}[D_x^s(u_0-u_\#(0))] + \frac 12 \varepsilon.
\]
Similarly, with hypothesis (H1) and consequence (C1), we are able to choose integration parameters such that
\[
\|u_\#\|_{L_t^\infty H_x^s(I\times\R)} \leq \mathcal K_{2,\infty,M_2,N_2}[u_\#] + \mathcal K_{2,\infty,M_2,N_2}[D_x^s u_\#] + \frac 12 \varepsilon
\]
and 
\[
\|u_\#\|_{Y_{k,s}(I\times\R)} \leq \mathcal Y_{k,s,N_3,M_3}[u_\#] + \frac 12 \varepsilon.
\]
Therefore
\[
\begin{aligned}
\|u_0 - u_\#(0)\|_{H_x^s(\R)} &\leq \mathcal J_{L_x^2,N_1}[(u_0-u_{\text{DNN},\#}(0))] \\
& \quad + \mathcal J_{L_x^2,N_1}[D_x^s(u_0-u_{\text{DNN},\#}(0))] + \frac 12 \varepsilon \\
& \leq \frac 12 (\widetilde A+\varepsilon) \leq \widetilde A,
\end{aligned}
\]
\[
\begin{aligned}
\|u_\#\|_{L_t^\infty H_x^s(I\times\R)} \leq &~{} \mathcal K_{2,\infty,M_2,N_2}[u_{\text{DNN},\#}] \\
&~{} + \mathcal K_{2,\infty,M_2,N_2}[D_x^s u_{\text{DNN},\#}] + \frac 12 \varepsilon \leq \frac 12 (A+\varepsilon) \leq A
\end{aligned}
\]
and
\[
\|u_\#\|_{Y_{k,s}(I\times\R)} \leq \mathcal Y_{k,s,N_3,M_3}[u_{\text{DNN},\#}] + \frac 12 \varepsilon \leq \frac 12 (L+\varepsilon) \leq L.
\]
Additionally, by using hypotheses (H1) and (H3), for $k \leq 4$
\[
\begin{aligned}
\|\mathcal E_k[u_\#]\|_{L_t^2H_x^s(I\times\R)} &\leq \mathcal J_{2,2,N_4,M_4}[\mathcal E_k[u_{\text{DNN},\#}]] \\
&\quad +\mathcal J_{2,2,N_4,M_4}[D_x^s\mathcal E_k[u_{\text{DNN},\#}]] + \frac 12 \varepsilon \\
& \leq \frac 12 \varepsilon + \frac 12 \varepsilon = \varepsilon,
\end{aligned}
\]
meanwhile for $k=5$,
\[
\|D_x^{-1}\mathcal E_k[u_\#]\|_{L_x^1L_t^2(I\times\R)} \leq \mathcal J_{1,2,N_4,M_4}[D_x^{-1}\mathcal E_k[u_{\text{DNN},\#}]] + \frac 12 \varepsilon \leq \varepsilon,
\]
and finally, by consequence (C1),
\[
\begin{aligned}
\|e^{-t\partial_x^3}(u_0-u_\#(0))\|_{Y_{k,s}(I\times\R)} &\leq \mathcal Y_{k,s,N_5,M_5}[e^{-t\partial_x^3}(u_0 - u_{\text{DNN},\#}(0))] + \frac 12 \varepsilon \leq \varepsilon.
\end{aligned}
\]
Previous inequalities implies that hypotheses \eqref{eq:long1}-\eqref{eq:long6} in Theorem \ref{thm:long} are satisfied, and one has that there exists a solution $u \in C(I;H^s(\R))$ to \eqref{eq:KdV} on $I \times \R$ with initial data $u_0 \in H^s(\R)$ such that
\begin{equation*}
\begin{aligned}
\|\partial_x(u^k-u^k_{\#})\|_{L_x^2L_t^2(I\times\R)} \,&  \leq\, C(A, \widetilde{A},L)\varepsilon,\\
    \|u-u_\#\|_{Y_{k,s}(I\times\R)}\,&\leq \, C(A, \widetilde{A},L)\varepsilon,\\
 \|u-u_\#\|_{L_t^{\infty}H_x^s(I\times\R)} \,&\leq\, C(A, \widetilde{A},L)(1 +  \varepsilon).
\end{aligned}
\end{equation*}
Finally, inequalities \eqref{eq:MT1}-\eqref{eq:MT3} came from the conclusion of Theorem \ref{thm:long}, the choice of the subinterval $[-R,R]$ and by recalling that $u_\# = u_{\text{DNN},\#}$ in this interval.
\qedsymbol

\section{Discussions and conclusions}\label{sec:discussion}



 Previous results show that our method can reproduce with high fidelity the dynamics of solution to the generalized KdV equation on unbounded domains. In particular, the analysis provides the first \emph{rigorous error bound} for neural-network-based approximations in a genuinely dispersive, unbounded, non-integrable setting, thus extending previous studies restricted to parabolic equations or bounded domains. This represents a step toward a more general theory of stability and convergence of neural architectures applied to nonlinear PDEs. 

Regarding the comparison of the results obtained with previous results in the field we can discuss the following. Referring the KdV model, previous works have a similar order of the relative error as whose we have presented (see, e.g. \cite{BKM22}), even though we are neglecting the boundary terms in the loss function. Additionally, despite the higher computational costs regarding the numerical computation of the Fourier transform, the average time that takes each simulation in all the examples is still in the order of similar works. 

Although the presented results are accurate, they are only valid for suitable parameters in each example. In particular, our method cannot address an efficient approximation of solutions with high velocities or oscillations, a problem naturally inherited from the classical PINNs. For those cases, other strategies should be considered or a modification of the PINNs method. Another direction that we could move forward is the applicability of PINNs on the approximation of orbital stability or even for blow-up solutions, that are well known to exist in the critical $L^2$ gKdV case. 

Another limitation to mention is that this stability results in conjunction with the PINNs method is restricted to the contraction principle in dispersive PDEs (for example the gKdV subcritical and critical models treated in this work as well as the nonlinear Schr\"odinger equation in \cite{ACFMV24}), but there are many other equations that do not posses such a principle (the Benjamin-Ono equation and the Kadomtsev-Petviashvili (KP) equation to mention some). To this end, it would be required efficient approximation on some Bourgain spaces, where the dispersive equations have a richer structure.

To finalize the discussion, another possible future work will address the supercritical regime, multidimensional models such as the Zakharov-Kuznetsov equation, and the inclusion of noise or uncertainty in the data. It would also be interesting to explore how architectural properties of the networks (depth, activation functions, and training constraints) affect dispersion and regularization in these systems.

\appendix
\renewcommand{\thesection}{R}
\setcounter{section}{0}

\end{document}